\documentclass[11pt,a4paper,leqno]{article}
\usepackage[T1]{fontenc}
\usepackage{lmodern}
\usepackage[margin=22mm,top=25mm,bottom=25mm]{geometry}
\usepackage{amsthm,amsmath,amsfonts,amssymb,mathrsfs}
\usepackage[numbers]{natbib}
\usepackage{booktabs,tabularx,array,enumitem,graphicx,xcolor,tikz}
\usepackage{microtype}
\usepackage[section]{placeins}
\usetikzlibrary{arrows.meta,positioning,shapes.geometric,fit,calc}
\setlist[itemize]{leftmargin=2em}
\setlist[enumerate]{leftmargin=2em}
\tikzset{
 rbox/.style={draw, rounded corners=2pt, align=center, inner sep=2pt, text width=5.3cm, font=\small},
 wbox/.style={draw, rounded corners=2pt, align=center, inner sep=2pt, text width=5.3cm, font=\small},
 keybox/.style={draw, rounded corners=2pt, align=center, inner sep=2pt, text width=5.3cm, font=\small},
 tinybox/.style={draw, rounded corners=2pt, align=center, inner sep=2pt, text width=5.3cm, font=\scriptsize},
 arr/.style={-{Stealth[length=2mm]}, thick},
 darr/.style={-{Stealth[length=2mm]}, thick, dashed}
}
\usepackage[colorlinks,citecolor=blue,urlcolor=blue]{hyperref}
\hypersetup{
 pdftitle={Uniform-in-Time Error Estimates and Invariant Measure Approximation for SDEs with SuperLinear Drift and Dissipativity at Infinity},
 pdfauthor={Shan Huang and Xiaoyue Li},
 pdfsubject={},
 pdfkeywords={65C30, 60C10, 65C40, Invariant measure, Uniform-in-time error, Truncated Euler-Maruyama scheme}
}
\theoremstyle{plain}

\newtheorem{theorem}{Theorem}[section]
\newtheorem{lemma}[theorem]{Lemma}
\newtheorem{corollary}[theorem]{corollary}
\newtheorem{proposition}{Proposition}
\theoremstyle{definition}
\newtheorem{definition}[theorem]{Definition}
\newtheorem{example}{Example}

\newtheorem{remark}{Remark}
\newtheorem{assumption}{Assumption}
\newenvironment{funding}{\section*{Funding}}{}
\title{Uniform-in-Time Error Estimates and Invariant Measure Approximation for SDEs with SuperLinear Drift and Dissipativity at Infinity}
\author{Shan~Huang\textsuperscript{1}\and Xiaoyue~Li\textsuperscript{2}\thanks{[\textbf{Corresponding author}]}}
\date{}
\begin{document}
\newcommand{\E}{\mathbb{E}}
\def\R{\mathbb{R}}
\def\p{p}
\def\1{{\bf 1}}
\def\e{\mathrm{e}}
\def\B{\mathscr{B}}
\def\L{\mathscr{L}}
\def\k{\kappa}
\def\md{\mathrm{d}}
\def\ep{\epsilon}
\def\bx{\bar{x}(t)}
\def\d{\delta}
\def\D{\Delta}
\def\tr{\triangle}
\maketitle
\begin{center}
\small
\textsuperscript{1}School of Mathematics and Statistics, Northeast Normal University\\
\href{mailto:huangs668@nenu.edu.cn}{huangs668@nenu.edu.cn}\\[0.5em]
\textsuperscript{2}School of Mathematical Sciences, Tiangong  University\\
\href{mailto:lixy@tiangong.edu.cn}{lixy@tiangong.edu.cn}
\end{center}
\begin{abstract}
For ergodic stochastic differential equations with drift coefficients that are dissipative only at infinity and grow superlinearly, can explicit numerical schemes preserve ergodicity and approximate invariant measures with quantitative convergence rates? This paper gives a positive answer to this question by proposing an explicit truncated Euler--Maruyama (TEM) scheme. Under a dissipativity-at-infinity condition, we prove the numerical ergodicity of the TEM scheme in the \(L^p\)-Wasserstein distance for all \(p\geqslant 1\). By combining truncation arguments with coupling techniques, we establish a uniform-in-time strong convergence rate of order \(1/2\) in moments. To the best of our knowledge, this is the first result to obtain such a uniform-in-time strong error estimate for numerical solutions under the framework of long-range dissipativity and superlinear growth through a purely discrete-time analysis. Moreover, using the exponential ergodicity of both the exact and numerical solutions, we derive an order \(1/2\) convergence rate for their invariant measures in the \(L^1\)-Wasserstein distance. Numerical experiments, including one for a high-dimensional nonlinear system, are presented to illustrate the theoretical results.
\end{abstract}

\medskip
\noindent\textbf{MSC:} Primary 65C30, 60C10; secondary 65C40.

\smallskip
\noindent\textbf{Keywords:} Invariant measure; Uniform-in-time error; Truncated Euler-Maruyama scheme.

\section{Introduction}
Invariant measures of ergodic stochastic differential equations (SDEs) play a central role in stochastic sampling, stochastic optimization, statistical mechanics and related areas. A prototypical example is the overdamped Langevin diffusion
\begin{equation}\label{E1}
	\mathrm d x_t=-\nabla U(x_t)\,\mathrm dt+\sqrt{2}\,\mathrm d B_t ,
\end{equation}
whose invariant distribution is the Gibbs measure associated with the potential $U$. Under suitable assumptions on $U$, the law of \eqref{E1} converges to its invariant measure as $t\to\infty$; see, for example, \cite{Meyn1993}. This long-time convergence provides the theoretical basis for many sampling and optimization algorithms. Since closed-form solutions of nonlinear SDEs are generally unavailable, numerical discretization is indispensable. A fundamental problem is therefore to design implementable schemes that approximate not only the finite-time dynamics but also the invariant measure of the underlying diffusion.

In this paper we study this problem for the more general additive-noise SDE
\begin{equation}\label{E}
	\mathrm d x_t=b(x_t)\,\mathrm dt+\Sigma\,\mathrm d B_t,\qquad x_0\in\mathbb R^d,
\end{equation}
where $B$ is a $q$-dimensional Brownian motion, $b:\mathbb R^d\to\mathbb R^d$ is locally Lipschitz and has superlinear growth, and $\Sigma\in\mathbb R^{d\times q}$ is non-degenerate in the sense that \(\Gamma:=\Sigma\Sigma^{\mathrm T}>0\).
Our interest is in explicit discretizations for \eqref{E} under a non-uniform dissipativity condition. More precisely, we assume that the drift is dissipative only at infinity: there exist constants $L,K>0$ and $R\geqslant 0$ such that
\begin{equation}\label{NGC}
	\langle x-y,b(x)-b(y)\rangle
	\leqslant L|x-y|^2\mathbf 1_{\{|x-y|\leqslant R\}}
	-K|x-y|^2\mathbf 1_{\{|x-y|>R\}},
	\qquad x,y\in\mathbb R^d .
\end{equation}
This condition is less restrictive than uniform dissipativity and also more flexible than the commonly used outside-of-sphere dissipativity condition; see Remark \ref{ree1} for details. It accommodates drifts with local oscillations and multiple-well structures, and therefore covers nonlinear models that are not included in the uniformly dissipative or outside-of-sphere frameworks. The condition \eqref{NGC} is satisfied by a number of systems arising from statistical mechanics and stochastic sampling. For instance, bistable potential-driven Langevin dynamics and stochastic overdamped Duffing-type equations \cite{McD} may fail to be uniformly dissipative, or may also fall outside the standard outside-of-sphere framework, while remaining dissipative at infinity in the sense of \eqref{NGC}. This illustrates that \eqref{NGC} offers a substantially broader framework for the long-time numerical analysis of nonlinear stochastic systems.

For linear SDEs, the Euler--Maruyama (EM) method is a standard and effective explicit discretization; when applied to the overdamped Langevin diffusion, it coincides with the unadjusted Langevin algorithm widely used in sampling \cite{Bernoulli}. However, for SDEs with superlinearly growing drift, the Euler--Maruyama method may diverge in both strong and weak senses even on finite time intervals \cite{Hutzenthaler}. Implicit schemes, such as backward Euler--Maruyama methods, stochastic $\theta$-methods and split-step backward Euler methods, provide suitable approximations for SDEs with superlinearly growing coefficients; however, their implementation requires the solution of an implicit equation at each time step. This motivates the development of explicit methods which remain stable in the long-time regime.

A large body of work has studied approximation of invariant measures under the uniform dissipativity condition
\begin{equation}\label{E2}
	\langle x-y,b(x)-b(y)\rangle\leqslant -K|x-y|^2,\qquad x,y\in\mathbb R^d ,
\end{equation}
which corresponds to strong convexity in the Langevin case. Under \eqref{E2}, convergence of numerical invariant measures has been obtained for both implicit and explicit schemes; see, for example, \cite{Brehier, ChenC, AMO, JiangY, LiuW, IM3, Talay} and \cite{IM2, IMA, SIAM}, respectively. Uniform-in-time strong error estimates have also been established in this setting \cite{Cui YP, Wu XM}. Non-asymptotic error bounds have been derived for Langevin-type dynamics in the linear setting \cite{JLMR}, for the Euler scheme with decreasing step sizes in the context of jump processes \cite{IM1}, and for tamed schemes applied to nonlinear Langevin diffusions \cite{Bro}. Further developments and related results under uniform dissipativity can be found in \cite{W.F, N.K} and the references therein. The analysis is substantially more delicate when \eqref{E2} fails, because one can no longer rely on global contraction of either the exact flow or the numerical transition.

For non-uniformly dissipative systems, coupling techniques have become an effective tool. Eberle and Majka \cite{EJP} developed contraction estimates for Markov chains on general state spaces, and Majka, Mijatovi\'c and Szpruch \cite{Majka} applied such ideas to derive non-asymptotic Wasserstein error bounds for Euler-type approximations of SDEs satisfying a non-uniform dissipativity condition. Related developments include \cite{Huang, 2025SPA, Schuh}. These works, however, are mainly formulated under global Lipschitz-type assumptions. For superlinearly growing drifts, tamed and projected algorithms have been studied in \cite{Bao, Bao1, Bao2, Bro, L.ZH, A.D. N, P.C}. In particular, \cite{Bao, Bao2, A.D. N} obtained approximation results for numerical invariant measures under outside-of-sphere dissipativity, while \cite{Bao1} established uniform-in-time strong convergence for tamed schemes in a different framework. Moreover, \cite{P.C, L.ZH} considered uniform-in-time weak error or non-asymptotic distributional bounds.

Taken together, the above results leave open the following natural problem. For SDEs with superlinearly growing drift under the dissipativity-at-infinity condition \eqref{NGC}, can one design an explicit and implementable numerical scheme which simultaneously provides a uniform-in-time strong error estimate with a quantitative convergence rate and an error estimate for invariant measures? To the best of our knowledge, existing results do not provide these two long-time estimates simultaneously in this setting. This problem is particularly relevant under \eqref{NGC}, since this condition imposes a less restrictive requirement on the drift coefficient \(b\) than the outside-of-sphere condition used in several previous works. It allows more complicated local behaviour of the drift, such as local oscillations or multiple-well structures, and therefore covers a broader class of nonlinear models.

This paper gives an affirmative answer to the above problem. We construct an explicit truncated Euler--Maruyama scheme for \eqref{E}, referred to as the TEM scheme, and investigate its long-time approximation properties in a fully discrete-time framework. First, by applying a Markov-chain ergodicity argument to the transition mechanism of the TEM scheme, we prove that, for every admissible step size \(h\), the numerical chain admits a unique invariant probability measure \(\mu_h\) and converges to \(\mu_h\) exponentially in Wasserstein distance. Thus the proposed explicit scheme preserves the ergodic structure of the exact diffusion under \eqref{NGC}. Second, under a polynomial Lipschitz condition on the drift, we establish the uniform-in-time strong error estimate
\begin{equation}\label{main-strong-intro}
	\sup_{k\geqslant 0}\mathbb E|X_k-x_{kh}|\leqslant C h^{1/2},
\end{equation}
where \(C\) is independent of \(k\) and \(h\). This gives a quantitative strong convergence rate that is uniform over the whole time horizon. Finally, combining \eqref{main-strong-intro} with the exponential ergodicity of both the exact diffusion and the numerical chain, we obtain the invariant-measure approximation
\begin{equation}\label{main-invariant-intro}
	\mathcal W_1(\mu_h,\mu)\leqslant C h^{1/2},
\end{equation}
where \(\mu\) is the invariant measure of \eqref{E}. Hence the TEM scheme approximates the invariant measure with order \(1/2\) in the \(L^1\)-Wasserstein distance.

We next explain the main difficulty and the key idea of the proof. The central issue is not merely to prove finite-time convergence, but to control the numerical error uniformly over an infinite time interval. Standard arguments based on global dissipativity or exponential contractivity of the numerical transition are not available under \eqref{NGC}, because the drift may be locally expansive even though it is dissipative at infinity. In such a situation, synchronous coupling alone cannot prevent the accumulation of errors over long times. Moreover, the superlinear growth of the drift rules out the direct use of the classical Euler--Maruyama method, while the truncation in the TEM scheme produces a numerical solution which is naturally a discrete-time Markov chain and does not naturally correspond to a continuous semimartingale to which It\^o's formula can be directly applied.

To overcome these obstacles, we develop a fully discrete coupling argument at the level of one-step transition kernels. The coupling links one step of the TEM approximation directly with the corresponding exact solution, and combines synchronous and reflection couplings according to the size and direction of the one-step separation. Since the diffusion matrix is allowed to be a non-degenerate constant matrix, the reflection is formulated in the metric induced by \(\Gamma^{-1}=(\Sigma\Sigma^{\mathrm T})^{-1}\). Together with a suitably constructed concave distance function and the truncation map, this one-step estimate can be propagated and iterated purely at the grid points. This is the mechanism that yields the uniform-in-time strong error estimate \eqref{main-strong-intro}; the invariant-measure error bound \eqref{main-invariant-intro} then follows by combining this estimate with the exponential ergodicity of the exact and numerical dynamics.

Let us finally clarify the relation with the existing literature. Compared with \cite{Bao, Bao2, A.D. N}, where invariant-measure approximation is studied under outside-of-sphere dissipativity, our analysis is carried out under the more flexible dissipativity-at-infinity condition \eqref{NGC} and also yields a uniform-in-time strong error estimate. Compared with the uniform-in-time strong convergence result in \cite{Bao1}, our explicit scheme is analyzed under dissipativity at infinity rather than outside-of-sphere dissipativity, and the proof is fully discrete rather than based on a continuous-time framework. Compared with \cite{Majka}, our result allows drifts with superlinear growth and further gives a quantitative error estimate for invariant measures. Compared with \cite{P.C, L.ZH}, which focus on weak errors, our principal estimate is a uniform-in-time strong error bound. Finally, the present formulation allows a non-degenerate constant diffusion matrix, going beyond the scalar-noise setting considered in several closely related works \cite{Bao, Bao1, Bao2, 2025SPA, Majka, A.D. N}.

The rest of the paper is organized as follows. Section 2 introduces notation, assumptions and preliminary ergodicity results. Section 3 defines the TEM scheme and proves its numerical ergodicity. Section 4 is devoted to the fully discrete coupling construction and the proof of the uniform-in-time strong error estimate; the invariant-measure approximation is then obtained as a consequence. Section 5 presents numerical experiments, including a high-dimensional example, to illustrate the theoretical results.
	
\section{Preliminaries}
\subsection{Notations and assumptions}
Throughout this paper, we use the following notations. Let $\mathbb{N}$ denote the set of all positive integers. Let $\R_+$ denote the set of all non-negative real numbers. 
Let $C$ denote a  generic positive constant whose value may change from line to line while $C_l$ is a positive constant depending on the parameter $l$. For any $a, b\in \R$, let $a\wedge b = \min\{a, b\}$, $a\vee b = \max\{a, b\}$ and \((a-b)^+ = a-a\wedge b\). Let $|\cdot|$ denote both the Euclidean norm in $\mathbb{R}^d$ and the Frobenius norm in $\mathbb{R}^{d \times m}$. 
For a vector $A\in\R^d$, denote its transpose by $A^\mathrm{T}$. For $x, y\in\mathbb{R}^d$, let $\langle x, y \rangle$ denote the inner product of $x$ and $y$. Let $\mathbf{I}$ be the $d\times d$ identity matrix. 
	
Let $(\Omega, \mathcal{F}, \{\mathcal{F}_t\}_{t\geqslant 0}, \mathbb{P})$ be a complete filtered probability space with $\{\mathcal{F}_t\}_{t\geqslant0}$ satisfying the usual conditions (that is, it is right-continuous and increasing while $\mathcal{F}_0$ contains all $\mathbb{P}$-null sets). Let $\mathbb{E}$ denote the expectation corresponding to $\mathbb{P}$.
Let $\{B_t\}_{t\geqslant 0}$ be a  $q$-dimensional Brownian motion defined on the probability space. Let $\mathcal{B}(E)$ denote the family of all Borel sets on $E$. For a set $\mathbb{D}$, let $\mathbb{D}^c$ be the complement of $\mathbb{D}$ and $\1_\mathbb{D}(x) = 1$ if $x\in\mathbb{D}$ and $0$ otherwise.  For $p\geqslant1$, let $L^p(\Omega;\R^d)$ be the family of \(\R^d\)-valued random variables $\xi$ satisfying $\|\xi\|_p:=(\E|\xi |^p)^{1/p}<\infty$. 
For a random variable $\xi\in L^p(\Omega;\R^d)$, let $\mathcal{L}(\xi)$ denote its distribution.
Let $N(a, \Gamma)$ denote the normal distribution with mean $a\in\R^d$ and positive definite covariance matrix $\Gamma\in \R^{d\times d}$. 
Let $p_{a,\Gamma}(\cdot)$ be the density function of $N(a,\Gamma)$ described by 
\begin{equation*}
p_{a,\Gamma}(x) := \frac{1}{\sqrt{(2\pi)^d\mathrm{det}(\Gamma)}}\mathrm{exp}\left(-\frac{1}{2}(x-a)^\mathrm{T}\Gamma^{-1}(x-a)\right).
\end{equation*} 

Let $\mathcal{P}(\R^{d})$ denote the set of all  probability measures on $\R^{d}$.
For $\mu\in \mathcal{P}(\R^d)$ and a measurable function $g$, define $\mu(g):=\int_{\R^d}g(x)\mu(\md x)$ whenever the integral exists. 
Let $f:[0, +\infty)\rightarrow [0, +\infty)$ be a strictly increasing function satisfying $f(0) = 0$.
Given two probability measures $\mu, \nu\in\mathcal{P}(\R^d)$ satisfying \(\mu(f(|\cdot|))<\infty\) and \(\nu(f(|\cdot|))<\infty\), define 
$$\mathcal{W}_f(\mu, \nu):= \inf_{\mathbf{\pi}\in \mathscr{C}(\mu, \nu)}\int_{\R^d\times \R^d} f(|x-y|)\mathbf{\pi}(\md x, \md y),$$
where $\mathscr{C}(\mu, \nu)$ is the collection of probability measures on $\R^d\times \R^d$ with marginals $\mu$ and $\nu$, respectively. 
For any $p\geqslant 1$, define
$$\mathcal{P}_p(\R^d):=\left\{\mu \in \mathcal{P}(\R^{d}) : \mu(|\cdot|^p)=\int_{\R^{d}}|x|^p\mu(\mathrm{d}x)<\infty\right\} $$
and the $L^p$-Wasserstein distance $\mathcal{W}_p$  for $\mu, \nu\in\mathcal{P}_p(\R^d)$
$$\mathcal{W}_p(\mu, \nu):= \inf_{\mathbf{\pi}\in \mathscr{C}(\mu, \nu)}\left(\int_{\R^d\times \R^d} |x-y|^p\mathbf{\pi}(\md x, \md y)\right)^{\frac{1}{p}}.$$
Then, $(\mathcal{P}_p(\R^d),\mathcal{W}_p)$ is a complete metric space for any $p\geqslant 1$, see \cite{Chen}.
For convenience, we impose the following assumptions.
\begin{assumption}\label{A}
For any $n>0$, there is a  constant $L_n>0$ such that 
\begin{equation}\label{LL}
|b(x)-b(y)|\leqslant L_n|x-y|
\end{equation}
for any $x, y\in\R^d$ with $|x|\vee|y|\leqslant n$.     
\end{assumption}
\begin{assumption}\label{A1}
There exist constants $L, K>0$ and $R\geqslant 0$ such that 
\begin{equation}\label{CaI}
	\langle x-y, b(x)-b(y)\rangle \leqslant L|x-y|^2 {\bf{1}}_{[0, R]}(|x-y|)-K|x-y|^2\1_{(R, +\infty)}(|x-y|) 
\end{equation}
for any $x, y\in\R^d$.
\end{assumption}
\begin{remark}\label{Re0}
It follows from \eqref{CaI} that for any $\epsilon\in(0, K)$,
\begin{equation*}
\langle x, b(x)\rangle
\leqslant(L+K)R^2-(K-\epsilon)|x|^2+\frac{|b(\mathbf{0})|^2}{4\epsilon}.
\end{equation*}
\end{remark}
\begin{remark}\label{ree1}
Some previous references \cite{Bao, Bao1, Bao2, 2025SPA, A.D. N} impose the outside-of-sphere contractivity condition. That is, there exist constants $L',K',R' > 0$ such that
\begin{equation}\label{DOS}
\langle x-y, b(x)-b(y)\rangle\leqslant L'|x-y|^2 {\bf{1}}_{[0, R']}(|x|\vee|y|)-K'|x-y|^2\1_{(R', +\infty)}(|x|\vee|y|)
\end{equation}
for any $x, y\in\R^d$. In the one-dimensional case, there are several functions $b$ satisfying (\ref{DOS}), such as the functions
$-x, x-x^3$ and $b_1(x) := -x|x|+6x$.
We observe that  
$(\ref{CaI})$ follows from $(\ref{DOS})$ with $L = L', K = K'$ and $R = 2R'$. 
However, $(\ref{DOS})$ is only sufficient but not necessary for $(\ref{CaI})$.
In fact, there exist several functions $b$ satisfying $(\ref{CaI})$ but not $(\ref{DOS})$.  For example, the functions $-x+\sin(2x), -3x/2+\ln(1+|x|)-2\sin (x), -x+(2-\cos(\pi x/4))\sin(2\pi x)$ and $b_2(x) := -x|x|+6x-2x \cos(x)$.
		
Figure \ref{Fre} geometrically visualizes the difference between the two dissipative properties described by (\ref{CaI}) and (\ref{DOS}). The functions $-x$ and $b_1(x)$ exhibit uniform dissipation outside a certain sphere. However, the class of functions satisfying dissipativity at infinity is broader: it includes not only such functions, but also functions with frequent local oscillations, such as \(-x+\sin(2x)\) and \(b_2(x)\).
\begin{figure}[!htb]
\centering
\includegraphics[width=12cm]{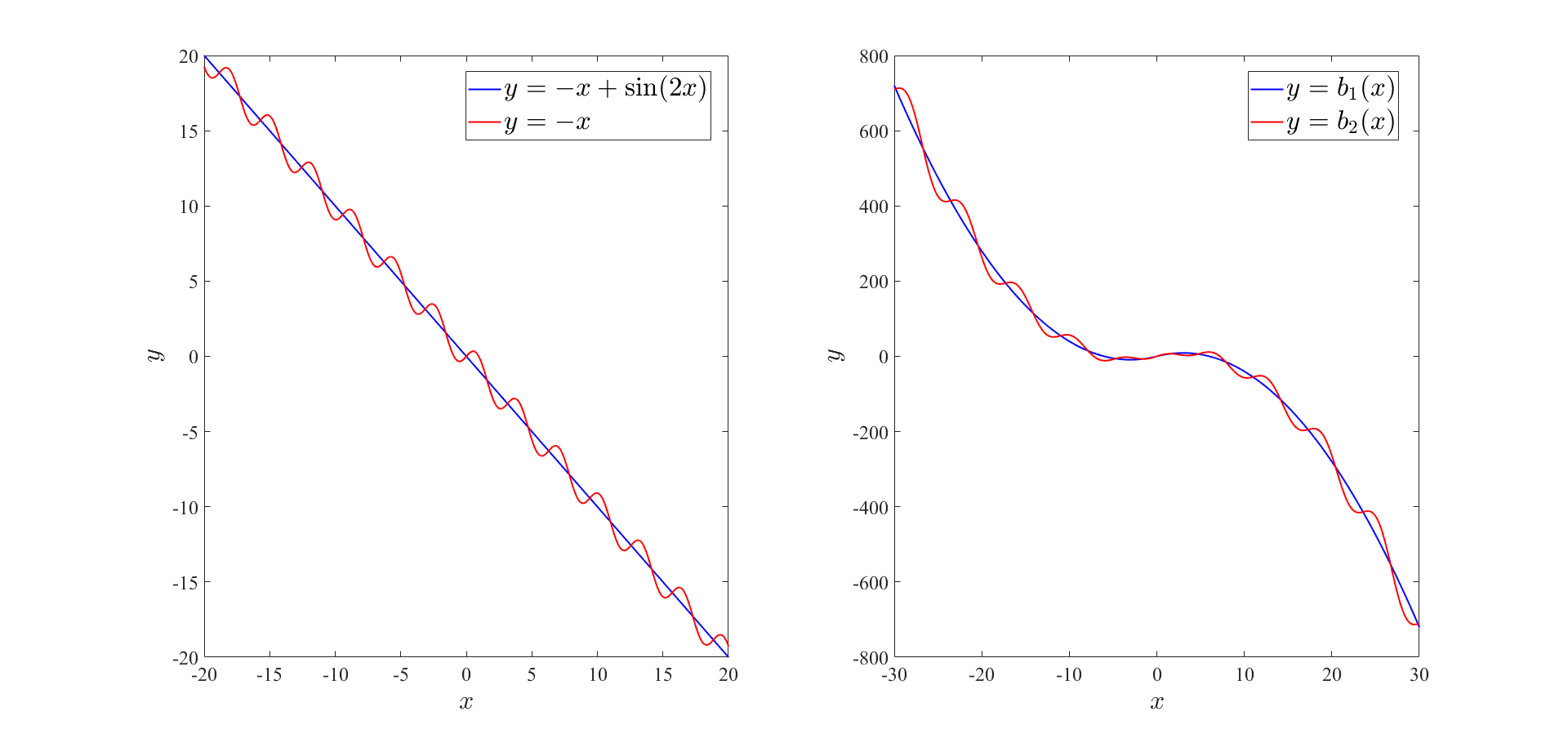}
\caption{Functions satisfying outside-of-sphere dissipativity or dissipativity at infinity.}
\label{Fre}
\end{figure}
\end{remark}
	
\subsection{Useful lemmas}
In order to study the numerical approximations, we prepare some results on the exact solutions of (\ref{E}) from \cite{IMA} and \cite{MN}. 
\begin{lemma}[{\cite[Theorems 2.3 and 5.2]{IMA}}]\label{Ju}
Under Assumption \ref{A}, the SDE \eqref{E} admits a unique global solution
\((x_t)_{t\geqslant 0}\) for every initial value \(x_0\in\mathbb R^d\).
Moreover, \((x_t)_{t\geqslant 0}\) is a time-homogeneous Markov process with
transition probability \(P_t(x,A):=\mathbb P(x_t\in A\mid x_0=x)\) for any \(t\geqslant0, A\in\mathcal B(\mathbb R^d)\).
If, in addition, Assumption \ref{A1} holds, then for every \(q>0\),
\[\sup_{0\leqslant t<\infty}\mathbb E|x_t|^q
\leqslant C_q(1+|x_0|^q),\]
where \(C_q\) is independent of \(t\).
\end{lemma}
	
\begin{lemma}[{\cite[Corollary 1.8]{MN}}]\label{LeX}
Under Assumptions \ref{A} and \ref{A1}, the SDE (\ref{E}) has a unique invariant probability measure $\mu\in \bigcap_{q\geqslant 1}\mathcal{P}_q(\R^d)$. Furthermore, there is a constant $\lambda>0$ such that for any $q\geqslant 1$ and $\nu\in\mathcal{P}_q(\R^d)$,
\begin{equation*}
\mathcal{W}_q(\nu P_t, \mu)\leqslant C\mathrm{e}^{-\lambda t/q}\mathcal{W}_{\phi_q}(\nu, \mu), \quad t\geqslant 0,
\end{equation*}
where $(P_t)_{t\geqslant 0}$ is the transition semigroup associated with (\ref{E}), the measure $\nu P_t$ is defined by
\[\nu P_t (A) = \int_{\R^d}P_t(x, A)\nu(\md x),\quad \forall A\in\mathcal{B}(\R^d),\]
and
\begin{equation*}
\phi_q(x) =\left\{
\begin{aligned}
&x^{1/q},& \text{if}~0\leqslant x< q^{-q/(q-1)},\\
&x-q^{-q/(q-1)}+q^{-1/(q-1)},& \text{if}~ x\geqslant q^{-q/(q-1)},
\end{aligned}\right.
\end{equation*}
for $q>1$ and $\phi_1(x) = x$ for all $x\geqslant 0$.
\end{lemma}	
	
To facilitate the proof of numerical ergodicity, we close this section by citing an ergodicity theorem. We first recall the definition of the transition probability of a general discrete-time homogeneous Markov chain.
\begin{definition}[{\cite{Meyn}}]
	Let \(\Phi=\{\Phi_n\}_{n\geqslant0}\) be a discrete-time homogeneous Markov chain on \((E,\mathcal{B}(E))\). For any \(m\geqslant0\), \(x\in E\), and \(A\in\mathcal{B}(E)\), define
	\[
	P^m_{\Phi}(x,A):=\mathbb{P}_x(\Phi_m\in A),
	\]
	which is called the \(m\)-step transition probability of the chain \(\Phi\). For notational simplicity, we denote the one-step transition probability \(P^1_{\Phi}(x,A)\) by \(P_{\Phi}(x,A)\).
\end{definition}
Next, we introduce an ergodicity theorem for the discrete-time homogeneous Markov chain \(\Phi=\{\Phi_n\}_{n\geqslant0}\).

\begin{lemma}[{\cite[Theorem 15.0.1]{Meyn}}]\label{GET}
	Assume that the following conditions hold:
	\begin{enumerate}
		\item[(i)] The \(E\)-valued time-homogeneous Markov chain
		\(\Phi=\{\Phi_n\}_{n\geqslant0}\) is \(\psi\)-irreducible and aperiodic, where \(\psi\) is a measure on \((E,\mathcal{B}(E))\).
		
		\item[(ii)] There exist a petite set \(A_*\in\mathcal{B}(E)\), constants
		\(\gamma<\infty\), \(\beta<1\), and a measurable function
		\(V:E\to[1,+\infty)\) such that
		\[
		\int_E V(y)P_{\Phi}(x,\md y)
		\leqslant \beta V(x)+\gamma \mathbf{1}_{A_*}(x),
		\qquad x\in E.
		\]
	\end{enumerate}
	Then the chain \(\Phi\) admits a unique invariant probability measure \(\pi\), and there exist constants \(\bar{c}>1\) and \(\tilde{c}>0\) such that, for any \(x\in E\),
	\begin{equation}\label{2.7}
		\sum_{m=1}^{\infty}\bar{c}^m
		\sup_{g\in\mathbb{G}_V}
		\left|
		\int_E g(y)P^m_{\Phi}(x,\md y)
		-
		\int_E g(y)\pi(\md y)
		\right|
		\leqslant \tilde{c}V(x),
	\end{equation}
	where
	\[
	\mathbb{G}_V
	:=
	\{g:E\to\R \text{ measurable}: |g(x)|\leqslant V(x),\ \forall x\in E\}.
	\]
\end{lemma}	
\begin{remark}\label{RE1}
	It follows from \eqref{2.7} that, for any \(x\in E\) and \(m\geqslant1\),
	\[
	\sup_{g\in\mathbb{G}_V}
	\left|
	\int_E g(y)P^m_{\Phi}(x,\md y)
	-
	\int_E g(y)\pi(\md y)
	\right|
	\leqslant
	\tilde{c}\mathrm{e}^{-m\ln\bar{c}}V(x).
	\]
\end{remark}
\section{TEM scheme and numerical ergodicity}
In this section, we propose an implementable explicit scheme and establish its numerical ergodicity.
\subsection{Construction of the numerical scheme}
Motivated by \cite{IMA}, we begin with the construction of an explicit scheme for the SDE \eqref{E}. Under the local Lipschitz condition \eqref{LL}, we choose a strictly increasing continuous function
\(\phi:\R_+\to\R_+\) such that \(\phi(u)\to\infty\) as \(u\to\infty\) and
\begin{equation}\label{Tr}
	\sup_{\substack{|x|\vee |y|\leqslant u\\ x\neq y}}
	\frac{|b(x)-b(y)|}{|x-y|}
	\leqslant \phi(u),
	\qquad u>0.
\end{equation}
\begin{remark}
The function \(\phi\) satisfying \eqref{Tr} is not unique. For example, if \eqref{Tr} holds for some function \(\phi\), then it also holds for any function in the family \(\{k\phi+l:\ k\geqslant1,\ l>0\}\).
Here we only need to choose one that works.
\end{remark}
For any $x\in\R^d$, define $\boldsymbol{e}_{x} = x/|x|$ if $x\neq \boldsymbol{0}$ and \(\boldsymbol{e}_{x} = \boldsymbol{0}\) otherwise.
Fix a parameter $M> \phi(0)\vee|b(\mathbf{0})|$ and a constant $\bar{\theta}\in(0, 1/2)$. For any step size $h\in(0,1]$ and $\theta\in(0, \bar{\theta}]$, define the truncation mapping $\pi_h:\R^d \rightarrow\R^d$ by
\begin{equation}\label{TF}
\pi_h(x) =  \Big(|x| \wedge \phi^{-1}(Mh^{-\theta})\Big) \boldsymbol{e}_{x} ,
\end{equation}
where $\phi^{-1}$ denotes the inverse function of $\phi$. 
It follows from (\ref{Tr}) and (\ref{TF}) that  
\begin{equation}\label{GL}
|b(\pi_h(x))-b(\pi_h(y))|\leqslant Mh^{-\theta}|\pi_h(x)-\pi_h(y)|,~~\forall~ x, y\in\R^d.
\end{equation}
This, together with \(M\geqslant |b(\mathbf{0})|\) and \(h\in(0,1]\), yields
\[
|b(\pi_h(x))|
\leqslant Mh^{-\theta}(1+|\pi_h(x)|),
\qquad x\in\R^d.
\]
Next, we propose an explicit scheme for (\ref{E}),  called the truncated Euler-Maruyama (TEM) scheme. For any integer $k\geqslant 0$, define
\begin{equation}\label{TEM}
\left\{
\begin{aligned}
X_0 &= \pi_h(x_0),\\
\hat{X}_{k+1} &= X_k + b(X_k)h + \Sigma Z_{k+1}, \\
X_{k+1} &= \pi_h(\hat{X}_{k+1}), 
\end{aligned}\right.
\end{equation}
where $t_k=kh$ and $Z_{k+1} := B_{t_{k+1}}-B_{t_{k}}$. 

\begin{remark}
	Assume that the drift coefficient \(b\) is globally Lipschitz continuous; that is, there exists a constant \(\tilde{L}>0\) such that $|b(x)-b(y)|\leqslant \tilde{L} |x-y|$ for all $x, y\in\R^d$.
	In this case, no truncation is needed. Equivalently, one may formally take
	\(\pi_h(x)=x\) for all \(x\in\R^d\). Then the TEM scheme \eqref{TEM}
	reduces to the classical Euler--Maruyama scheme.
	Therefore, the TEM scheme can be regarded as a generalization of the EM scheme to SDEs with nonlinear drift coefficients.
\end{remark}

We cite the time-homogeneous Markov property and the uniform boundedness of the numerical solution \((X_k)_{k\geqslant0}\) from \cite[p.882, Lemma 7.8; p.869, Theorem 5.5]{IMA}.
	
\begin{lemma}[{\cite[p.882, Lemma 7.8; p.869, Theorem 5.5]{IMA}}]\label{Le2} Let Assumptions \ref{A} and \ref{A1} hold.  
$(X_k)_{k\geqslant 0}$ is a time-homogeneous Markov chain with the $m$-step transition probability $P^m_X(x, A) = \mathbb{P}(X_m\in A \mid X_0 = x)$ and has the property
\begin{equation*}
\sup_{h\in(0, \bar{h}]}\sup_{k\geqslant 0}\E |X_k|^q\leqslant \sup_{h\in(0, \bar{h}]}\sup_{k\geqslant 1}\E |\hat{X}_k|^q\leqslant C_q(1+|x_0|^q), 
\end{equation*}
for any $q>0$, where $\bar{h}:= (KM^{-2})^{1/(1-2\bar{\theta})}\wedge\left(2K^{-1}\right)$ provided that $M\geqslant \sqrt{K}$.
\end{lemma}
	
To close this subsection, we prove the Feller continuity of the numerical solution.
\begin{lemma}\label{Le3} Let Assumption \ref{A} hold.
\((X_k)_{k\geqslant0}\) is a Feller Markov chain.
\end{lemma}	

\begin{proof}
	Let \((X_k^{x_0})_{k\geqslant0}\) and
	\((X_k^{x'_0})_{k\geqslant0}\) be two numerical solutions of \eqref{TEM}
	with initial values \(x_0,x'_0\in\R^d\), respectively, driven by the same
	Brownian increments. By the argument in \cite[p.142, Lemma 8.1.4]{Ok},
	it is sufficient to prove that, for any fixed \(k\geqslant1\),
	\[
	\lim_{x_0\to x'_0}\E|X_k^{x_0}-X_k^{x'_0}|^2=0.
	\]	
	Since \(\pi_h\) is the projection onto the closed ball centered at the origin
	with radius \(\phi^{-1}(Mh^{-\theta})\), it is non-expansive, that is, for any \(x,y\in\R^d\), \(|\pi_h(x)-\pi_h(y)|\leqslant |x-y|\). 
	Hence,
	using \eqref{GL} and \eqref{TEM}, for any fixed integer \(k\geqslant1\),
	\begin{align*}
		\E|X_k^{x_0}-X_k^{x'_0}|^2
		&\leqslant
		\E\left|
		X_{k-1}^{x_0}-X_{k-1}^{x'_0}
		+h\left(b(X_{k-1}^{x_0})-b(X_{k-1}^{x'_0})\right)
		\right|^2\\
		&\leqslant
		2(1+M^2h^{2-2\theta})
		\E|X_{k-1}^{x_0}-X_{k-1}^{x'_0}|^2.
	\end{align*}
	Iterating this inequality gives
	\[
	\E|X_k^{x_0}-X_k^{x'_0}|^2
	\leqslant
	2^k(1+M^2h^{2-2\theta})^k
	|\pi_h(x_0)-\pi_h(x'_0)|^2.
	\]
	Since \(\pi_h\) is non-expansive, we further obtain
	\[
	\E|X_k^{x_0}-X_k^{x'_0}|^2
	\leqslant
	2^k(1+M^2h^{2-2\theta})^k
	|x_0-x'_0|^2.
	\]
Letting $x_0\rightarrow x'_0$, the desired inequality \label{Fe} follows. The proof is complete.
\end{proof}
	
\subsection{Numerical ergodicity}
In this subsection, we aim to prove the existence and uniqueness of the numerical invariant measure. Furthermore, we show that the distribution of the numerical solution $(X_k)_{k\geqslant 0}$ converges exponentially to the numerical invariant measure in the $L^p(p\geqslant 1)$-Wasserstein distance.

We next establish the ergodicity of the discrete-time Markov chain induced by
the TEM transition kernel. The proof relies on the Markov-chain ergodicity
criterion recalled in Lemma \ref{GET}. More precisely, the criterion requires irreducibility, aperiodicity, a suitable petite set and a Lyapunov drift condition, and then yields the existence and uniqueness of an invariant probability measure together with geometric ergodicity in the weighted \(V\)-norm. When such criteria are applied to numerical approximations of SDEs, however, their assumptions are not automatic and must be verified for the resulting discrete-time Markov chain. The essential points are to prove that the numerical chain is irreducible and satisfies a suitable Lyapunov drift condition, together with the required petite-set
condition. For nonlinear SDEs with superlinear coefficients, this strategy has been used in \cite{LiuWu} to establish the numerical ergodicity of tamed Euler--Maruyama schemes. Since the present TEM scheme is different from the tamed scheme, and in order to make the argument self-contained, we verify below the corresponding irreducibility, petite-set and Lyapunov conditions for the TEM chain.

\begin{theorem}\label{Pro1}
Let Assumptions \ref{A} and \ref{A1} hold. For any integer $p\geqslant 1$, $\theta\in(0, \bar{\theta}]$ and $h\in(0, \bar{h}]$, the TEM chain $(X_k)_{k\geqslant 0}$  admits a unique invariant measure $\mu_h$ on \(\mathbb{B}^h:=\{x\in\R^d: |x|\leqslant \phi^{-1}(Mh^{-\theta})\}\). Moreover, there exist constants $\bar{c}>1$, $\tilde{c}>0$, possibly
depending on \(p,h,\theta\), such that  
\begin{equation*}
\sum_{k = 1}^{\infty}\bar{c}^k\sup_{g\in\mathbb{G}_{V_p}}|\E g(X_k)-\mu_h(g)|\leqslant \tilde{c}(1+|x_0|^{2p}), 
\end{equation*}
where $\bar{h}$ is defined in Lemma \ref{Le2}, $V_p(x) = 1+|x|^{2p}$ and \[\mathbb{G}_{V_p}: = \{g:\R^d\rightarrow\R:~|g(x)|\leqslant V_p(x), ~\forall x\in \R^d\}.\]
\end{theorem}
	
\begin{proof}
Fix $p\geqslant 1$, $\theta\in(0, \bar{\theta}]$ and $h\in(0, \bar{h}]$. Let \(\mathcal{E}_h:=\mathcal{B}(\R^d)\cap \mathbb{B}^h\),
which is a $\Sigma$-algebra by \cite[p. 145, Lemma 3]{Pro}.
For the numerical ergodicity, we only verify that conditions $(\mathrm{i})$ and $(\mathrm{ii})$ hold in Lemma \ref{GET} with $E = \mathbb{B}^h$. Since the proof is rather technical, we divide it into two steps.
		
{\bf Step 1: Irreducibility and aperiodicity.} Let \(\mu^{\mathrm{Leb}}\) denote the Lebesgue measure restricted to \(\mathbb{B}^h\). For \(x\in \mathbb B^h\), the random variable \(x+b(x)h+\Sigma Z_1\) has a strictly positive Gaussian density on \(\R^d\). Hence, for any \(x\in\mathbb{B}^h\) and nonempty set
\(A\in\mathcal{E}_h\) with \(\mu^{\mathrm{Leb}}(A)>0\),
\[
\begin{aligned}
P_X(x,A)
=
	\mathbb{P}\big(\pi_h(x+b(x)h+\Sigma Z_1)\in A\big)  &\geqslant
	\mathbb{P}\big(x+b(x)h+\Sigma Z_1\in A\big) \\
	&=
	\int_A p_{x+b(x)h,h\Sigma \Sigma^{\mathrm T}}(y)\,\md y
	>0.
\end{aligned}
\]
Therefore, according to \cite[Proposition 4.2.1]{Meyn}, \((X_k)_{k\geqslant0}\) is  $\mu^{\mathrm{Leb}}$-irreducible.

We next show that it is aperiodic. Suppose, to the contrary, that the chain has
period \(d_0\geqslant2\). Then, according to the definition of a periodic Markov chain \cite[p.122-123]{Meyn}, there exist disjoint nonempty sets
\(D_1,\ldots,D_{d_0}\in\mathcal{E}_h\) such that
\[
\mathbb{B}^h=\bigcup_{i=1}^{d_0}D_i,
\qquad
P_X(x,D_{i+1})=1,\quad x\in D_i,
\]
with the convention \(D_{d_0+1}=D_1\). Since $\mu^{\mathrm{Leb}}(\mathbb{B}^h)>0$, there exists an $i_0\in \{1, 2, \cdots, d_0\}$  such that $\mu^{\mathrm{Leb}}(D_{i_0} )>0$. By the irreducibility argument above,
\[
P_X(x,D_{i_0})>0,
\qquad x\in E_h.
\]
In particular, this holds for \(x\in D_{i_0}\). However, the cyclic decomposition
implies
\[
P_X(x,D_{i_0})=0,
\qquad x\in D_{i_0},
\]
because the chain must move from \(D_{i_0}\) to \(D_{i_0+1}\) in one step. This is
a contradiction. Hence the chain is aperiodic.
		
{\bf Step 2: Lyapunov drift condition.} For any $x\in \mathbb{B}^h$, by \(|\pi_h(z)|\leqslant |z|\), we have
\begin{align*}
\mathbb{E}V_p(\pi_h(x+b(x)h+\Sigma Z_1)) &= 1+\E|\pi_h(x+b(x)h+\Sigma Z_1)|^{2p}\leqslant 1+\E|x+b(x)h+\Sigma Z_1|^{2p},
\end{align*}
which yields 
\begin{equation}\label{Proeq1}
\begin{aligned}
&\mathbb{E}V_p(\pi_h(x+b(x)h+\Sigma Z_1)) \\
\leqslant &1+\E\left(|x|^2+2h\langle x, b(x) \rangle+h^2|b(x)|^2+|\Sigma Z_1|^2+2\langle x+b(x)h,\Sigma Z_1 \rangle\right)^p.
\end{aligned}
\end{equation}
Using Remark \ref{Re0} with $\epsilon = K/4$ implies that for any $x\in\mathbb{B}^h$
\begin{equation}\label{Reeq}
\langle x, b(x)\rangle \leqslant -\frac{3K}{4}|x|^2+R^2(L+K)+\frac{1}{K}|b(\mathbf{0})|^2.
\end{equation}
By (\ref{GL}), (\ref{Reeq}) and $M^2h^{2-2\theta}\leqslant Kh$, we further have for any $x\in\mathbb{B}^h$
\begin{align*}
2h\langle x, b(x) \rangle+h^2|b(x)|^2 &\leqslant-\frac{3Kh}{2}|x|^2+2hR^2(L+K)+\frac{2h}{K}|b(\mathbf{0})|^2+M^2h^{2-2\theta} \left|x\right|^2\\
&\leqslant -\frac{Kh}{2}|x|^2+q_1h,
\end{align*}
where $q_1 := 2R^2(L+K)+2|b(\mathbf{0})|^2/K$.
Substituting this into (\ref{Proeq1}), we obtain that
\begin{equation}\label{EV}
\mathbb{E}V_p(\pi_h(x+b(x)h+\Sigma Z_1))\!\leqslant\! 1\!+\!\E\!\left(\!\left(1\!-\!Kh/2\right)|x|^2+q_1h\!+\!|\Sigma Z_1|^2\!+\!2\langle x\!+\!b(x)h, \!\Sigma Z_1 \rangle\!\right)^p\!.
\end{equation}
For $p =1$, since $\mathbb{E}\langle x+b(x)h, \Sigma Z_1 \rangle = 0$ and $\E|\Sigma Z_1|^{2p} \leqslant C_dh^p$ (see \cite{BMSC}), we have
\begin{equation}\label{EV1}
\begin{aligned}
\mathbb{E}V_1(\pi_h(x+b(x)h+\Sigma Z_1))\leqslant 1+\left(1-\frac{Kh}{4}\right)|x|^2+\gamma_1 h,
\end{aligned}
\end{equation}
where $\gamma_1:= q_1+C_d$. Next, we consider the case $p\geqslant 2$. For any $h\in(0, 2/K]$, we have \(1-Kh/2\in[0,1)\). Together with (\ref{EV}) and the binomial theorem, this gives that for any $x\in\mathbb{B}^h$
\begin{align*}
\mathbb{E}V_p(\pi_h(x+b(x)h+\Sigma Z_1)) 
\leqslant 1+\left(1-\frac{Kh}{2}\right)|x|^{2p}+\mathrm{I}_1+\mathrm{I}_2,
\end{align*}
where
$$\mathrm{I}_1 := p|x|^{2p-2}\left(q_1h+\E|\Sigma Z_1|^2+2\E\langle x+b(x)h, \Sigma Z_1 \rangle\right),$$
$$\mathrm{I}_2:=\sum_{k =0}^{p-2}C_p^k|x|^{2k}\E\left(q_1h+|\Sigma Z_1|^2+2\langle x+b(x)h, \Sigma Z_1 \rangle\right)^{p-k},$$
and $C_p^k := p!/(k!(p-k)!)$.
Next, we estimate $\mathrm{I}_1$ and $\mathrm{I}_2$. 
Similar to the argument for \eqref{EV1}, we have	
\begin{equation}\label{I1}
\mathrm{I}_1\leqslant p|x|^{2p-2}\left(q_1h+C_d\Sigma^2 h^p\right)\leqslant C_d|x|^{2p-2}h\leqslant Kh|x|^{2p}/8+C_dh.	
\end{equation}
Since $C_p^k\leqslant C_p^{\lfloor p/2 \rfloor}$ for any $k\in[0, p]\cap \mathbb{N}$, $\E|Z_1|^{p-k} \leqslant C_dh^{(p-k)/2}\leqslant C_dh$ for any $k\in[0,p-2]\cap\mathbb{N}$, and $|b(x)|h\leqslant Mh^{1-\theta}(1+|x|)\leqslant M\sqrt{h}(1+|x|)$ for any $x\in\mathbb{B}^h$, we have
\begin{equation}\label{I2}
\begin{aligned}
\mathrm{I}_2\leqslant C_p^{\lfloor p/2 \rfloor}\sum_{k =0}^{p-2}|x|^{2k}\E\left(q_1h+|\Sigma Z_1|^2+2\langle x+b(x)h, \Sigma Z_1 \rangle\right)^{p-k}\!\leqslant\! C_d\sum_{k =0}^{p-2}\left(|x|^{2k}h^2+|x|^{p+k}h\right).	
\end{aligned}
\end{equation}
Noticing that $p/(p-k)\in[1, p/2]$ and $k/(p-k)\in[0, p/2-1]$ for any $k\in [0, p-2]$, by Young's inequality and \eqref{I2}, we have
\begin{equation}\label{I2-}
\mathrm{I}_2\leqslant \frac{Kh}{8}|x|^{2p} + C_dh.
\end{equation}
This, together with (\ref{I1}) and (\ref{I2-}), implies that for any $x\in\mathbb{B}^h$
\begin{equation}\label{Proeq5}
\mathbb{E}V_p(\pi_h(x+b(x)h+\Sigma Z_1)) \leqslant \left(1-\frac{Kh}{4}\right)|x|^{2p}+\gamma_p,
\end{equation}
where $\gamma_p = 1+C_dh$. Therefore, combining (\ref{EV1}) with (\ref{Proeq5}), we obtain that for any $p\in[1, +\infty)\cap\mathbb{N}$ and $x\in\mathbb{B}^h$,
$$\mathbb{E}V_p(\pi_h(x+b(x)h+\Sigma Z_1)) \leqslant \left(1-\frac{Kh}{4}\right)|x|^{2p}+\gamma_p = \left(1-\frac{Kh}{8}\right)|x|^{2p}+\left(\gamma_p -\frac{Kh}{8}|x|^{2p}\right).$$
Hence, for any $p\in\mathbb{N}$ and $x\in\mathbb{B}^h$,
\begin{equation}\label{Ev}
			\mathbb{E}V_p(\pi_h(x+b(x)h+\Sigma Z_1)) \leqslant \gamma_p\1_{A_p}(x)+\left(1-\frac{Kh}{8}\right) V_p(x),
		\end{equation}
		where the compact set %$A_p$ is defined by
		$A_p := \left\{x\in\R^d: |x|^{2p}\leqslant 8\gamma_p/(Kh)\right\}\cap\mathbb{B}^h.$ 
		Since $(X_k)_{k\geqslant 0}$ is a Feller process and $(X_k)_{k\geqslant 0}$ is $\mu^{\mathrm{Leb}}$-irreducible, $A_p$ is a petite set for $(X_k)_{k\geqslant 0}$ by using \cite[Proposition 6.2.8]{Meyn}.
		This, together with (\ref{Ev}), implies that condition $(\mathrm{ii})$ in Lemma \ref{GET} holds with $V_p(x) = 1+|x|^{2p}$, $\beta = 1-Kh/8$, $\gamma = \gamma_p$ and petite set $A_p$. 
		According to Lemma \ref{GET} and Remark \ref{RE1}, there exists a unique invariant probability $\mu_h$ such that for any $k\geqslant 1$,
		$$
		\sup_{g\in \mathbb{G}_{V_p}}|\E g(X_k)-\mu_h(g)|\leqslant \tilde{c} \mathrm{e}^{-k\ln \bar{c}}V_p(X_0)$$
		for suitable constants $\tilde{c},\bar{c}>0$. This completes the proof.
	\end{proof}
\begin{remark}
	In \cite{IMA}, the existence and uniqueness of the numerical invariant measure for the TEM scheme were obtained under the uniformly dissipative condition. Since Assumption \ref{A1} reduces to the uniformly dissipative case when \(R=0\) in \eqref{CaI}, the present result extends the corresponding result in \cite{IMA} to the weaker long-distance dissipativity setting.
\end{remark}

Next, we show that the numerical solution \((X_k)_{k\geqslant0}\) of the TEM scheme \eqref{TEM} is exponentially ergodic in the \(L^q\)-Wasserstein distance for any \(q\geqslant1\).

\begin{theorem}\label{Co}
	Let Assumptions \ref{A} and \ref{A1} hold. For any \(h\in(0,\bar h]\), where \(\bar h\) is given in Lemma \ref{Le2}, for any \(q\geqslant1\) and \(k\geqslant1\), the numerical solution \((X_k)_{k\geqslant0}\) satisfies
	\begin{equation*}
		\mathcal{W}_q(\mathcal{L}(X_k), \mu_h)\leqslant C_q (1+|x_0|)\mathrm{e}^{-\bar{c}'k},
	\end{equation*}
	where  $\bar{c}' = \ln \bar{c}/(q+2)$ and $\bar{c}>1$ are given in Theorem \ref{Pro1}. Here \(C_q>0\) is a constant independent of \(k\).
\end{theorem}
\begin{proof}
	We first recall the following standard estimate connecting the Wasserstein distance with a weighted total variation norm \cite[Theorem 6.15]{CV}. For any
	\(\mu,\nu\in\mathcal{P}_{2p}(\R^d)\) with $\int_{\R^d}|x|^{2p}\md |\mu-\nu|(x)<\infty$,
	\[\mathcal{W}_{2p}^{2p}(\mu,\nu)
	\leqslant 2^{2p}\int_{\R^d}(1+|x|^{2p})|\mu-\nu|(\md x).\]
	By the dual representation of the weighted total variation norm \cite[p.110-111]{Hairer},
	\[
	\int_{\R^d}(1+|x|^{2p})|\mu-\nu|(\md x)
	=
	\sup_{g\in\mathbb{G}_{V_p}}
	|\mu(g)-\nu(g)|,
	\]
	where \(V_p(x)=1+|x|^{2p}\).	
	Taking \(\mu=\mathcal{L}(X_k)\) and \(\nu=\mu_h\), and using Theorem \ref{Pro1}, we get
	\[
	\sup_{g\in\mathbb{G}_{V_p}}
	|\E g(X_k)-\mu_h(g)|
	\leqslant
	\tilde c \bar c^{-k}(1+|x_0|^{2p}).
	\]
	Therefore,
	\[\mathcal{W}_{2p}(\mathcal{L}(X_k),\mu_h)
		\leqslant
		2\left(
		\sup_{g\in\mathbb{G}_{V_p}}
		|\E g(X_k)-\mu_h(g)|
		\right)^{1/(2p)}  \leqslant
		2\tilde c^{1/(2p)}
		\exp\left(-\frac{k\ln\bar c}{2p}\right)
		(1+|x_0|).\]	
	Now let \(q\geqslant1\). Choose \(p\in\mathbb{N}\) such that \(q\in[2(p-1),2p)\).
	Then, by the monotonicity of Wasserstein distances,
	\[
	\mathcal{W}_q(\mathcal{L}(X_k),\mu_h)
	\leqslant
	\mathcal{W}_{2p}(\mathcal{L}(X_k),\mu_h).
	\]
	Therefore,
	\[
	\mathcal{W}_q(\mathcal{L}(X_k),\mu_h)
	\leqslant
	C_q(1+|x_0|)
	\exp\left(-\frac{k\ln\bar c}{q+2}\right).
	\]
	This completes the proof.
\end{proof}

Theorem \ref{Co} shows that the TEM scheme preserves the exponential ergodicity of the exact solution in the Wasserstein space
\((\mathcal{P}_p(\R^d),\mathcal{W}_p)\) for any \(p\geqslant1\), in agreement with the ergodicity property of the exact solution stated in Lemma \ref{LeX}.

\section{Uniform-in-time error analysis}
This section is devoted to establishing the uniform-in-time \(1/2\)-order convergence rate of the numerical solution \((X_k)_{k\geqslant0}\), generated by the TEM scheme \eqref{TEM}, to the exact solution \((x_t)_{t\geqslant0}\) under a condition slightly stronger than the local Lipschitz continuity of the drift.
Moreover, we derive the \(1/2\)-order convergence rate of the numerical invariant measure \(\mu_h\) to the underlying invariant measure \(\mu\) in the \(L^1\)-Wasserstein distance.
We impose the following polynomial Lipschitz condition on the drift.
	\begin{assumption}\label{A2}
		There exist constants $L^*>0$ and $\ell>0$ such that for any $x, y\in\R^d$,
		\begin{equation*}
			|b(x)-b(y)|\leqslant L^*(1+|x|^\ell+|y|^\ell)|x-y|.
		\end{equation*}
	\end{assumption}
	
\begin{remark}\label{PG}
			Under Assumption \ref{A2}, by Young's inequality, one notes that for any $x\in\R^d$
			\[|b(x)|\leq \frac{L^*(\ell+2)}{\ell+1}|x|^{\ell+1}+\frac{L^*\ell}{\ell+1}+|b(\mathbf{0})|.\]
\end{remark}

\begin{remark}\label{mapping form} 
Under Assumption \ref{A2}, we observe  that any function $\phi$ with $\phi(u)\geqslant L^*(1+2u^\ell)$ for any $u>0$ satisfies \eqref{Tr}. In particular, we choose
$\phi(u) = L^*(1+2u^\ell),~\forall~u\geqslant0$. Then, $\phi^{-1}(u) = ((u-L^*)/(2L^*))^{1/\ell}$ for all $u\geqslant L^*$. 
\end{remark}
	
For convenience, we introduce the auxiliary one-step process $(\hat{S}_{k})_{k=1}^\infty$ by
\begin{equation}\label{S}
\hat{S}_{k}:= \pi_h(x_{t_{k-1}}) + hb(\pi_h(x_{t_{k-1}})) + \Sigma Z_{k},\quad k\geqslant 1.
\end{equation}
By the triangle inequality, we obtain the following error decomposition: 
\begin{align*}
\E\left|X_{k}-x_{t_{k}}\right|\leqslant\E|X_{k}-\hat{X}_{k}|+\E |\hat{X}_{k}-\hat{S}_{k}|+ \E |\hat{S}_{k}-x_{t_{k}}|.
\end{align*}
Therefore, it remains to estimate the three terms $\E|X_{k}-\hat{X}_{k}|$, $\E|\hat{X}_{k}-\hat{S}_{k}|$ and $\E|\hat{S}_{k}-x_{t_{k}}|$. 
	
\subsection{Upper bounds of $\E|X_{k}-\hat{X}_{k}|$ and $\E |\hat{S}_{k}-x_{t_{k}}|$}
In this subsection, we estimate the uniform upper bounds of
\(\E |X_k-\hat X_k|\) and \(\E|\hat S_k-x_{t_k}|\) for \(k\geqslant0\).
These estimates are the basic ingredients for proving the uniform-in-time
strong convergence rate of the TEM scheme. The following lemma quantifies the error introduced by the projection mapping
\(\pi_h\). It will be used to control both the truncation error of the exact
solution and the projection error of the numerical solution.
\begin{lemma}\label{Le4}
	For any \(M\geqslant 3L^*\), \(p,q>1\) with \(1/p+1/q=1\), and
	\(s,n\geqslant1\), if
	\[
	\zeta\in L^{np}(\Omega;\R^d)\cap L^{sq}(\Omega;\R^d),
	\]
	then
	\[
	\E \left|\pi_h(\zeta)-\zeta\right|^n
	\leqslant
	2^n\|\zeta\|^n_{np}\|\zeta\|^s_{sq}
	\left(3L^*M^{-1}\right)^{s/\ell}h^{s\theta/\ell}.
	\]
\end{lemma}

	\begin{proof}
		Define
		$$\mathbb{A}_h := \left\{\omega\in\Omega: |\zeta(\omega)|\leqslant \left(\frac{Mh^{-\theta}-L^*}{2L^*}\right)^{1/\ell}\right\}.$$
		On \(\mathbb A_h\), we have \(\pi_h(\zeta)=\zeta\). Hence, by H{\"o}lder's inequality \cite[p.5]{Mao} and the fact $|\pi_h(\zeta)|\leqslant |\zeta|$, we obtain that for any $n\geqslant 1$,
		\begin{equation}\label{Leeq1}
			\E|\zeta\!-\!\pi_h(\zeta)|^n \!=\! \E\!\left(|\zeta\!-\!\pi_h(\zeta)|^n\1_{\mathbb{A}_h^c}\!\right) \leqslant \left(\E |\zeta-\!\pi_h(\zeta)|^{np}\right)^\frac{1}{p}\!\left(\mathbb{P}(\mathbb{A}_h^c)\right)^{\frac{1}{q}}\leqslant 2^n \|\zeta\|^n_{np}\left(\mathbb{P}(\mathbb{A}_h^c)\right)^{\frac{1}{q}}\!.
		\end{equation}
		By Chebyshev's inequality \cite[p.5]{Mao}, 
		$$\left(\mathbb{P}(\mathbb{A}_h^c)\right)^{\frac{1}{q}}\leqslant \|\zeta\|_{sq}^s\left(\frac{Mh^{-\theta}-L^*}{2L^*}\right)^{-s/\ell}.$$
		Since \(M\geqslant3L^*\) and \(h\in(0,1]\), we have
		$(Mh^{-\theta}-L^*)/(2L^*)\geqslant Mh^{-\theta}/(3L^*).$
		Consequently, for any $s\geqslant 1$, $\left((Mh^{-\theta}-L^*)/(2L^*)\right)^{-s/\ell}\leqslant \left(3L^*M^{-1}\right)^{s/\ell}h^{s\theta/\ell}$.
		Substituting this estimate into \eqref{Leeq1} gives the desired result.
	\end{proof}
	Applying Lemma \ref{Le4} together with the uniform moment bounds of the exact
	and numerical solutions, we obtain the following uniform truncation estimates.
	These estimates will be used repeatedly in the local error analysis.
	\begin{corollary}\label{co2}
		Let Assumptions \ref{A1} and \ref{A2} hold. For any
		\(M\geqslant |b(\mathbf{0})|\vee 3L^*\), \(\theta\in(0,\bar{\theta}]\), and
		\(h\in(0,\bar h]\), the exact solution \((x_t)_{t\geqslant0}\) and the numerical
		solution \((X_k)_{k\geqslant0}\) defined by \eqref{TEM} satisfy that, for any
		\(s,n\geqslant1\),
		\[
		\sup_{t\geqslant0}\E|\pi_h(x_t)-x_t|^n
		\leqslant
		Ch^{s\theta/\ell},
		\qquad
		\sup_{k\geqslant0}\E|X_k-\hat X_k|^n
		\leqslant
		Ch^{s\theta/\ell}.
		\]
	\end{corollary}
	
The next lemma estimates the local consistency error between the auxiliary
one-step approximation \(\hat S_k\) and the exact solution \(x_{t_k}\).
 \begin{lemma}\label{Le5}
 	Let Assumption \ref{A2} hold. Then \(\sup_{k\geqslant1}\E|\hat S_k-x_{t_k}|
 	\leqslant
 	Ch^{3/2}\).
 \end{lemma}
 
 \begin{proof}
 	Combining \eqref{E} with \eqref{S}, we have
 	\begin{equation}\label{LEeq}
 			\E |\hat S_{k+1}-x_{t_{k+1}}|
 			=
 			\E\left|
 			\pi_h(x_{t_k})
 			+h b(\pi_h(x_{t_k}))
 			-x_{t_k}
 			-\int_{t_k}^{t_{k+1}} b(x_t)\md t
 			\right|  \leqslant
 			J_1+J_2+J_3,
 	\end{equation}
 	where
 	\[
 	J_1:=\E|\pi_h(x_{t_k})-x_{t_k}|,
 	\qquad
 	J_2:=h\E|b(\pi_h(x_{t_k}))-b(x_{t_k})|,
 	\]
 	and
 	\[
 	J_3:=
 	\E\left|
 	\int_{t_k}^{t_{k+1}}
 	\bigl[b(x_{t_k})-b(x_t)\bigr]\md t
 	\right|.
 	\]	
 	By Corollary \ref{co2}, choosing
 	\(s\geqslant 1\vee 3\ell/(2\theta)\), we get
 	\begin{equation}\label{LEeq1}
 		J_1
 		\leqslant
 		Ch^{s\theta/\ell}
 		\leqslant
 		Ch^{3/2}.
 	\end{equation}	
 	Using Assumption \ref{A2}, Hölder's inequality, Lemma \ref{Ju}, and
 	Corollary \ref{co2}, for
 	\(s\geqslant 1\vee \ell/(2\theta)\), we obtain
 	\begin{equation}\label{LEeq2}
 		\begin{aligned}
 			J_2
 			&\leqslant
 			hL^*
 			\E\left[
 			(1+|\pi_h(x_{t_k})|^\ell+|x_{t_k}|^\ell)
 			|\pi_h(x_{t_k})-x_{t_k}|
 			\right]  \\
 			&\leqslant
 			hL^*
 			\left(
 			\E(1+|\pi_h(x_{t_k})|^\ell+|x_{t_k}|^\ell)^2
 			\right)^{1/2}
 			\left(
 			\E|\pi_h(x_{t_k})-x_{t_k}|^2
 			\right)^{1/2}  \\
 			&\leqslant
 			Ch^{1+s\theta/\ell}
 			\leqslant
 			Ch^{3/2}.
 		\end{aligned}
 	\end{equation} 	
 	For \(J_3\), by Assumption \ref{A2}, Hölder's inequality, and Lemma \ref{Ju},
 	\begin{equation}\label{LEeq3}
 		\begin{aligned}
 			J_3
 			&\leqslant
 			L^*
 			\int_{t_k}^{t_{k+1}}
 			\E\left[
 			(1+|x_{t_k}|^\ell+|x_t|^\ell)
 			|x_{t_k}-x_t|
 			\right]\md t  \\
 			&\leqslant
 			L^*
 			\left(
 			\int_{t_k}^{t_{k+1}}
 			\E(1+|x_{t_k}|^\ell+|x_t|^\ell)^2\md t
 			\right)^{1/2}
 			\left(
 			\int_{t_k}^{t_{k+1}}
 			\E|x_{t_k}-x_t|^2\md t
 			\right)^{1/2}  \\
 			&\leqslant
 			C\sqrt h
 			\left(
 			\int_{t_k}^{t_{k+1}}
 			\E|x_t-x_{t_k}|^2\md t
 			\right)^{1/2}.
 		\end{aligned}
 	\end{equation}
 	Moreover, by Hölder's inequality, Itô's isometry, Remark \ref{PG}, and
 	Lemma \ref{Ju}, for any \(t\in(t_k,t_{k+1}]\),
 	\begin{equation}\label{Lqq4}
 		\begin{aligned}
 			\E|x_t-x_{t_k}|^2
 			&=
 			\E\left|
 			\int_{t_k}^t b(x_u)\md u
 			+
 			\Sigma(B_t-B_{t_k})
 			\right|^2  \\
 			&\leqslant
 			2(t-t_k)\int_{t_k}^t \E|b(x_u)|^2\md u
 			+
 			2|\Sigma|^2(t-t_k)  \\
 			&\leqslant
 			C\bigl[(t-t_k)^2+(t-t_k)\bigr].
 		\end{aligned}
 	\end{equation}
 	Substituting \eqref{Lqq4} into \eqref{LEeq3} yields
 	\begin{equation}\label{LEeq4}
 			J_3
 			\leqslant
 			C\sqrt h
 			\left(
 			\int_{t_k}^{t_{k+1}}
 			\bigl[(t-t_k)^2+(t-t_k)\bigr]\md t
 			\right)^{1/2}  \leqslant
 			C\sqrt h\,(h^3+h^2)^{1/2}
 			\leqslant
 			Ch^{3/2}.
 	\end{equation}
 	Combining \eqref{LEeq}, \eqref{LEeq1}, \eqref{LEeq2}, and \eqref{LEeq4}, we obtain \(\E|\hat S_{k+1}-x_{t_{k+1}}|
 	\leqslant
 	Ch^{3/2}\).
 	Taking the supremum over \(k\geqslant0\) completes the proof.
 \end{proof}

\subsection{Construction of coupling}
It remains to estimate \(\E|\hat X_k-\hat S_k|\). Since the drift term in
\eqref{E} is dissipative only at infinity, a uniform-in-time estimate of
\(\E|\hat X_k-\hat S_k|\) cannot be obtained directly by the classical
synchronous coupling method. Inspired by \cite{Majka}, we construct in this
subsection a one-step coupling for the TEM scheme by combining synchronous
coupling with reflection coupling. This coupling is designed to handle the
non-uniform dissipativity of the drift and will be used to derive a
uniform-in-time bound for \(\E|\hat X_k-\hat S_k|\).

We first formulate the coupling in a form which is valid for a general
constant non-degenerate diffusion matrix. For fixed \(h>0\), consider one step of the TEM scheme. For any \(x\in\R^d\),
let \(Z\sim N(\mathbf{0},h\mathbf{I}_q)\) and define
\[
\hat X_h(x,Z):=u_h(x)+\Sigma Z,
\]
where \(u_h(x):=\pi_h(x)+h b(\pi_h(x))\) and \(\pi_h\) is defined by \eqref{TF}. Note that \(\Sigma Z\sim N(\mathbf{0},h\Gamma)\), where \(\Gamma:=\Sigma\Sigma^{\mathrm T}\). For \(x,y\in\R^d\), set
\[
\hat r_h(x,y):=u_h(x)-u_h(y).
\]
We shall couple \(\hat X_h(x, Z)\) and \(\hat X_h(y, Z)\).
For \(\mathbf{0}\neq r\in \R^d\) and \(a, b\in \R^d\), define
\[
\langle a,b\rangle_\Gamma:=a^{\mathrm T}\Gamma^{-1}b,
\qquad
|a|_\Gamma:=\langle a,a\rangle_\Gamma^{1/2},
\]
and define the \(\Gamma^{-1}\)-reflection operator \(\mathcal R_r^\Gamma\) by
\[
\mathcal R_r^\Gamma w
:=
w-2\frac{\langle r,w\rangle_\Gamma}{|r|_\Gamma^2}r,
\qquad w\in\mathbb R^d.
\]
Furthermore, define the one-dimensional projection
\[
\lambda_r(w)
:=
\frac{|r|}{|r|_\Gamma^2}\langle r,w\rangle_\Gamma,
\qquad w\in\mathbb R^d,
\]
and, for \(m>0\),
\[
A_r^m
:=
\{w\in\mathbb R^d:|\lambda_r(w)|\leqslant m\}.
\]
Finally, let
\[
\varphi_\Gamma(w)
:=
\frac{1}{\sqrt{(2\pi h)^d\det\Gamma}}
\exp\left\{-\frac{1}{2h}|w|_\Gamma^2\right\}
\]
be the density of \(N(0,h\Gamma)\), and set
\[
v_r^m(w)
:=
\left[
\frac{\varphi_\Gamma(r+w)}{\varphi_\Gamma(w)}
\mathbf 1_{A_r^m}(w)\mathbf 1_{A_r^m}(r+w)
\right]
\wedge
\mathbf 1_{A_r^m}(w).
\]
When \(r=\mathbf{0}\), we set \(v_{\mathbf{0}}^m\equiv1\).
Let \(\zeta\) be uniformly distributed on \([0,1]\), independently of \(Z\).
For \(\hat r_h(x,y)=\mathbf{0}\), set \(\hat Y_h(x,y, Z):=u_h(y)+\Sigma Z\).
For \(\hat r_h(x,y)\neq\mathbf{0}\), define
\begin{equation}\label{hY}
\hat Y_h(x,y, Z)
:=
\begin{cases}
	u_h(x)+\Sigma Z,
	& \zeta\leqslant v_{\hat r_h(x,y)}^m(\Sigma Z),\ \Sigma Z\in A_{\hat r_h(x,y)}^m,\ |\hat r_h(x,y)|\leqslant 2R,\\[1mm]
	u_h(y)+\mathcal R_{\hat r_h(x,y)}^\Gamma\Sigma Z,
	& \zeta> v_{\hat r_h(x,y)}^m(\Sigma Z),\ \Sigma Z\in A_{\hat r_h(x,y)}^m,\ |\hat r_h(x,y)|\leqslant 2R,\\[1mm]
	u_h(y)+\Sigma Z,
	& \Sigma Z\notin A_{\hat r_h(x,y)}^m\ \text{or}\ |\hat r_h(x,y)|>2R .
\end{cases}
\end{equation}

We shall prove in Lemma \ref{coup} that \(\mathcal{L}(\hat Y_h(x,y,Z))=\mathcal{L}(\hat X_h(y,Z))\). Consequently, the joint distribution of \((\hat X_h(x,Z),\hat Y_h(x,y,Z))\) is a coupling of the laws of \(\hat X_h(x,Z)\) and \(\hat X_h(y,Z)\). We briefly explain the construction. If \(|\hat r_h(x,y)|\leqslant 2R\), we
first compare the projected noise \(\lambda_{\hat r_h(x,y)}(\Sigma Z)\) with the threshold \(m\). If this quantity is no larger than \(m\), we then
compare \(\zeta\) with \(v_{\hat r_h(x,y)}^m(\Sigma Z)\). When \(\zeta\leqslant v_{\hat r_h(x,y)}^m(\Sigma Z)\), we set \(\hat Y_h(x,y,Z)=\hat X_h(x,Z)\),
so that the two components coincide on this event.
Otherwise, we use reflection coupling and set \(\hat Y_h(x,y,Z)
=
u_h(y)+R_{\hat r_h(x,y)}^\Gamma \Sigma Z\).
If \(|\lambda_{\hat r_h(x,y)}(\Sigma Z)|>m\) or \(|\hat r_h(x,y)|>2R\),
we use synchronous coupling, namely, \(\hat Y_h(x,y,Z)=\hat X_h(y,Z)\).
\begin{remark}
	The coupling constructed above is formulated entirely in a discrete-time framework. Our construction is tailored to the one-step transition of the TEM scheme. In particular, the coupling is applied directly to the nonlinear numerical map
	\[x\mapsto u_h(x)=\pi_h(x)+h b(\pi_h(x)),\]
	rather than to the exact continuous-time flow. This is essential for the analysis of the numerical scheme, since the projection operator \(\pi_h\) and the nonlinear drift both enter the transition mechanism. Moreover, the present construction uses a directional control of the noise. Instead of restricting the full noise magnitude \(|\Sigma Z|\) as in \cite{Majka}, we control only its projection along the direction of the one-step separation, \(|\lambda_{\hat r_h(x,y)}(\Sigma Z)|\).
	This directional truncation is more refined, since only the component of the noise that directly affects the coupling distance is restricted. It also reduces the relevant estimates to one-dimensional Gaussian integrals along the coupling direction, which is useful for treating high-dimensional problems in the subsequent analysis.
\end{remark}
\begin{lemma}[Properties of the \(\Gamma^{-1}\)-reflection]\label{reflprop}
For \(\mathbf 0 \neq r\in\mathbb R^d\), the following properties hold.
\begin{equation}\label{R1}
	\mathcal R_r^\Gamma r=-r,
	\qquad
	\bigl(\mathcal R_r^\Gamma\bigr)^2=\mathbf I .
\end{equation}	
For every \(w\in\mathbb R^d\),
\begin{equation}\label{R2}
	\langle r,\mathcal R_r^\Gamma w\rangle_\Gamma
	=
	-\langle r,w\rangle_\Gamma,\qquad |\mathcal R_r^\Gamma w|_\Gamma=|w|_\Gamma.
\end{equation}	
Moreover,
\begin{equation}\label{R3}
	\det(\mathcal R_r^\Gamma)=-1,
	\qquad
	|\det(\mathcal R_r^\Gamma)|=1.
\end{equation}
Consequently, for the Gaussian density \(\varphi_\Gamma\),
\begin{equation}\label{R4}
	\varphi_\Gamma(\mathcal R_r^\Gamma w)=\varphi_\Gamma(w),
	\qquad w\in\mathbb R^d .
\end{equation}
Furthermore,
\begin{equation}\label{R5}
	\lambda_r(\mathcal R_r^\Gamma w)
	=
	-\lambda_r(w),
	\qquad w\in\mathbb R^d,
\end{equation}
and hence, for every \(m>0\) and \(w\in A_r^m\)
\begin{equation}\label{R6}
w\in A_r^m\Longleftrightarrow
	\mathcal R_r^\Gamma w\in A_r^m.
\end{equation}	
Finally, for every \(u\in\mathbb R^d\),
\begin{equation}\label{R7}
r+\mathcal R_r^\Gamma u
	=
	\mathcal R_r^\Gamma(u-r).
\end{equation}
\end{lemma}

\begin{proof}
	By the definition of \(\mathcal R_r^\Gamma\), \(\mathcal R_r^\Gamma r = -r\).
	Also,
	\[
	\begin{aligned}
		\langle r,\mathcal R_r^\Gamma w\rangle_\Gamma
		=
		\left\langle
		r,
		w-2\frac{\langle r,w\rangle_\Gamma}{|r|_\Gamma^2}r
		\right\rangle_\Gamma     =
		\langle r,w\rangle_\Gamma
		-
		2\frac{\langle r,w\rangle_\Gamma}{|r|_\Gamma^2}
		\langle r,r\rangle_\Gamma  =
		-\langle r,w\rangle_\Gamma .
	\end{aligned}
	\]
Applying this identity once more gives
	\[\bigl(\mathcal R_r^\Gamma\bigr)^2w=
\mathcal R_r^\Gamma w-2\frac{\langle r,\mathcal R_r^\Gamma w\rangle_\Gamma}{|r|_\Gamma^2}r =w-2\frac{\langle r,w\rangle_\Gamma}{|r|_\Gamma^2}r+
2\frac{\langle r,w\rangle_\Gamma}{|r|_\Gamma^2}r=w.\]
	Hence \eqref{R1} follows.	
	Next,
	\[|\mathcal R_r^\Gamma w|_\Gamma^2
=\left|w-2\frac{\langle r,w\rangle_\Gamma}{|r|_\Gamma^2}r\right|_\Gamma^2                    =|w|_\Gamma^2-4\frac{\langle r,w\rangle_\Gamma^2}{|r|_\Gamma^2}
+4\frac{\langle r,w\rangle_\Gamma^2}{|r|_\Gamma^2}=|w|_\Gamma^2 .\]
This proves \eqref{R2}. Since
\[\det(\mathcal R_r^\Gamma)
=1-2\frac{r^{\mathrm T}\Gamma^{-1}r}{|r|_\Gamma^2}=-1,\]
\eqref{R3} holds. By the definition of \(\varphi_\Gamma(w)\),
\eqref{R3} immediately implies \eqref{R4}. By the definition of \(\lambda_r\) and \eqref{R2},
\[\lambda_r(\mathcal R_r^\Gamma w)
=\frac{|r|}{|r|_\Gamma^2}\langle r,\mathcal R_r^\Gamma w\rangle_\Gamma
=-\frac{|r|}{|r|_\Gamma^2}\langle r,w\rangle_\Gamma
=-\lambda_r(w).\]
This proves \eqref{R5}, which implies \eqref{R6}. Finally, \eqref{R1}
immediately implies \eqref{R7}. The proof is complete.
\end{proof}
\begin{remark}
	For a fixed nonzero vector \(r\), the \(\Gamma^{-1}\)-reflection \(\mathcal R_r^\Gamma w\)
	is the reflection with respect to the hyperplane \(H_\Gamma
	=
	\{x\in \R^d:\langle r,x\rangle_\Gamma=0\}\).
	In other words, the reflecting hyperplane is determined by the notion of
	orthogonality induced by the \(\Gamma^{-1}\)-inner product.
	In the example shown in Figure~\ref{sigma-ref}, we take \(	r=(1,1)^{\mathrm T},
	w=(2,1)^{\mathrm T}\).
	When \(\Gamma=I\), the \(\Gamma^{-1}\)-inner product \(\langle \cdot, \cdot\rangle_\Gamma\) is the usual Euclidean
	inner product. Hence the reflecting hyperplane is \(H_{\mathbf I}=\{(x_1,x_2):x_1+x_2=0\}\).
	Thus \(\mathcal R_r^{\mathbf I}\) maps the point \(w\) to \(\mathcal R_r^{\mathbf I} w=(-1,-2)^{\mathrm T}\).
	The corresponding unit ball is \(\{x = (x_1, x_2):|x|_{\mathbf I}\leqslant 1\}
	=
	\{(x_1,x_2):x_1^2+x_2^2\leqslant 1\}\),
	whose boundary is the usual unit circle. By contrast, when \(\Gamma=\operatorname{diag}(4,1)\),
	distances, angles, and orthogonality are changed by
	\(\Gamma^{-1}\), rather than the standard dot product. The reflecting hyperplane becomes \(H_\Gamma
	=
	\{(x_1,x_2): x_1/4+x_2=0\}\).
	Therefore the same point \(w\) is reflected to a different point, \(\mathcal R_r^\Gamma w
	=(-2/5,-7/5)\). In this geometry, the unit ball is \(\{x=(x_1, x_2):|x|_\Gamma\leqslant 1\}
	=
	\{(x_1, x_2):x_1^2/4+x_2^2\leqslant 1\}\). So its boundary is an ellipse rather than a circle. Thus, as illustrated in Figure~\ref{sigma-ref}, the difference
	between the two cases is that when \(\Gamma=I\), reflection is the usual
	Euclidean mirror reflection, whereas when \(\Gamma\neq I\), the mirror
	hyperplane, the unit ball, and the reflected point are all determined by
	the weighted geometry induced by \(\Gamma^{-1}\).
\end{remark}

\begin{figure}[h!]
	\centering
	\resizebox{\textwidth}{!}{		
	\begin{tikzpicture}[scale=0.8, >=Stealth]
		
		% -------------------- First panel --------------------
		\begin{scope}[shift={(-4.8,0)}]
			
			% Axes
			\draw[->] (-4.2,0) -- (4.2,0) node[right] {$x_1$};
			\draw[->] (0,-4.2) -- (0,4.2) node[above] {$x_2$};
			
			% Mirror line y = -x
			\draw[thick, blue] (-2.8,2.8) -- (2.8,-2.8)
			node[pos=0.75, below right] {$x_1+x_2=0$};
			
			% Unit circle
			\draw[thick, green!60!black] (0,0) circle (1)node[pos=0.75, above left]{$x_1^2+x_2^2 = 1$~~~~~~~~~~~~~~};
			
			% Vector r
			\draw[->, thick] (0,0) -- (1,1) node[above right] {$r$};
			
			% Points
			\filldraw[blue] (2,1) circle (1.5pt)
			node[right] {$w=(2,1)$};
			
			\filldraw[blue] (-1,-2) circle (1.5pt)
			node[left] {$\mathcal R_r^{\bf I} w=(-1,-2)$};
			
			% Dashed segment
			\draw[dashed, blue] (2,1) -- (-1,-2);
			
			% Label
			\node at (0,-4.8) {\(\Gamma={\bf I}\): ordinary Euclidean reflection};
			
		\end{scope}
		
		% -------------------- Second panel --------------------
		\begin{scope}[shift={(5.2,0)}]
			
			% Axes
			\draw[->] (-4.2,0) -- (4.2,0) node[right] {$x_1$};
			\draw[->] (0,-4.2) -- (0,4.2) node[above] {$x_2$};
			
			% Mirror line y = -x/4
			\draw[thick, blue] (-4,1) -- (4,-1)
			node[pos=0.8, below right] {$\frac14 x_1+x_2=0$};
			
			% Unit ellipse x^2/4 + y^2 = 1
			\draw[thick, green!60!black] (0,0) ellipse (2 and 1) node[pos=0.75, above left = 1.1cm]{$\frac{x_1^2}{4}+x_2^2 = 1$};
			
			% Vector r
			\draw[->, thick] (0,0) -- (1,1) node[above right] {$r$};
			
			% Points
			\filldraw[blue] (2,1) circle (1.5pt)
			node[right] {$w=(2,1)$};
			
			\filldraw[blue] (-0.4,-1.4) circle (1.5pt)
			node[left] {$\mathcal R_r^\Gamma w=\left(-\frac25,-\frac75\right)$};
			
			% Dashed segment
			\draw[dashed, blue] (2,1) -- (-0.4,-1.4);
			
			% Label
			\node at (0,-4.8) {\(\Gamma=\operatorname{diag}(4,1)\): \(\Gamma^{-1}\)-reflection};
			
		\end{scope}			
	\end{tikzpicture}}		
	\caption{Comparison between the ordinary Euclidean reflection and the
		\(\Gamma^{-1}\)-reflection.}
		\label{sigma-ref}
\end{figure}

\begin{lemma}\label{coup}
For every \(h, m>0\) and \(x,y\in\R^d\), the joint distribution of the random vector \((\hat X_h(x,Z),\hat Y_h(x,y,Z))\) defined by \eqref{hY} is a coupling of the laws of \(\hat X_h(x,Z)\) and \(\hat X_h(y,Z)\).
\end{lemma}
\begin{proof}
We need to show that
\[\mathcal L(\hat Y_h(x,y,Z))=\mathcal L(u_h(y)+\Sigma Z).\]
We write $\hat{r}_h = \hat{r}_h(x, y)$ and $|\hat{r}_h| = |\hat{r}_h(x, y)|$ for short. If \(\hat r_h=\mathbf 0\) or \(|\hat r_h|>2R\), then by construction of \(\hat Y_h(x,y,Z)\), the conclusion is immediate.

Assume therefore that \(\hat r_h\neq\mathbf 0\) and \(|\hat r_h|\leqslant 2R\).Let \(g:\mathbb R^d\to\mathbb R\) be bounded and measurable. We write $\hat{r}_h = \hat{r}_h(x, y)$ and $|\hat{r}_h| = |\hat{r}_h(x, y)|$ for short.
According to the definition of $\hat{Y}_h(x, y, Z)$ given in (\ref{hY}) and using the independence of \(\zeta\), we immediately get 
	\[
	\begin{aligned}
		\mathbb E g(\hat Y_h(x,y,Z))
		&=
		\int_{(A_{\hat r_h}^m)^c}
		g(u_h(y)+w)\varphi_\Gamma(w)\,\mathrm dw +
		\int_{A_{\hat r_h}^m}
		g(u_h(x)+w)
		v_{\hat r_h}^m(w)
		\varphi_\Gamma(w)\,\mathrm dw                 \\
		&\quad+
		\int_{A_{\hat r_h}^m}
		g\bigl(
		u_h(y)+\mathcal R_{\hat r_h}^\Gamma w
		\bigr)
		\bigl(1-v_{\hat r_h}^m(w)\bigr)
		\varphi_\Gamma(w)\,\mathrm dw .
	\end{aligned}
	\]
	By the definitions of \(u_h(x)\) and \(v_{\hat r_h}^m\),
	\[v_{\hat r_h}^m(w)\varphi_\Gamma(w)=
		\bigl[
		\mathbf 1_{A_{\hat r_h}^m}(w)
		\varphi_\Gamma(w)
		\bigr] \wedge
		\bigl[
		\mathbf 1_{A_{\hat r_h}^m}
		\bigl(\hat r_h+w\bigr)
		\varphi_\Gamma\bigl(\hat r_h+w\bigr)
		\bigr],\]
	we obtain
	\[
	\begin{aligned}
		&\mathbb E g(\hat Y_h(x,y,Z))
		\\
		=&\int_{(A_{\hat r_h}^m)^c}
		g(u_h(y)+w)\varphi_\Gamma(w)\,\mathrm dw \\
		&+\int_{\mathbb R^d}
		g\bigl(u_h(y)+\hat r_h+w\bigr)
		\Bigl(\bigl[
		\mathbf 1_{A_{\hat r_h}^m}(w)
		\varphi_\Gamma(w)
		\bigr] \wedge
		\bigl[
		\mathbf 1_{A_{\hat r_h(x,y)}^m}
		\bigl(\hat r_h+w\bigr)
		\varphi_\Gamma\bigl(\hat r_h+w\bigr)
		\bigr]
		\Bigr)\,\mathrm dw                                         \\
		&+
		\int_{\mathbb R^d}
		g\bigl(u_h(y)+\mathcal R_{\hat r_h}^\Gamma w\bigr)
		\Bigl(\mathbf 1_{A_{\hat r_h}^m}(w)
		\varphi_\Gamma(w)-\Bigl(
		\bigl[
		\mathbf 1_{A_{\hat r_h}^m}(w)
		\varphi_\Gamma(w)
		\bigr]\wedge
		\bigl[
		\mathbf 1_{A_{\hat r_h}^m}
		\bigl(\hat r_h+w\bigr)
		\varphi_\Gamma\bigl(\hat r_h+w\bigr)
		\bigr]
		\Bigr)
		\Bigr)\,\mathrm dw .
	\end{aligned}
	\]
	In the second integral, use \(s=\hat r_h+w\), which gives
	\[\begin{aligned}
		&\int_{\mathbb R^d}
		g\bigl(u_h(y)+\hat r_h+w\bigr)\Bigl(
		\bigl[
		\mathbf 1_{A_{\hat r_h}^m}(w)
		\varphi_\Gamma(w)
		\bigr]\wedge
		\bigl[
		\mathbf 1_{A_{\hat r_h}^m}
		\bigl(\hat r_h+w\bigr)
		\varphi_\Gamma\bigl(\hat r_h+w\bigr)
		\bigr]
		\Bigr)\,\mathrm dw \\
		=&
		\int_{\mathbb R^d}
		g(u_h(y)+s)\Bigl(
		\bigl[
		\mathbf 1_{A_{\hat r_h}^m}
		\bigl(s-\hat r_h\bigr)
		\varphi_\Gamma\bigl(s-\hat r_h\bigr)
		\bigr]\wedge
		\bigl[
		\mathbf 1_{A_{\hat r_h}^m}(s)
		\varphi_\Gamma(s)
		\bigr]
		\Bigr)\,\mathrm ds .
	\end{aligned}\]	
	In the third integral, use \(v=\mathcal R_{\hat r_h}^\Gamma w\). By \eqref{R3}, the Jacobian is one. By \eqref{R4} and \eqref{R6},
	\[
	\mathbf 1_{A_{\hat r_h}^m}
	\bigl(\mathcal R_{\hat r_h}^\Gamma v\bigr)
	\varphi_\Gamma
	\bigl(\mathcal R_{\hat r_h}^\Gamma v\bigr)
	=
	\mathbf 1_{A_{\hat r_h}^m}(v)\varphi_\Gamma(v).
	\]
	Moreover, by \eqref{R7}, \eqref{R4}, and \eqref{R6},
	\[\mathbf 1_{A_{\hat r_h}^m}
		\bigl(
		\hat r_h+
		\mathcal R_{\hat r_h}^\Gamma v
		\bigr)
		\varphi_\Gamma
		\bigl(
		\hat r_h+
		\mathcal R_{\hat r_h}^\Gamma v
		\bigr) =
		\mathbf 1_{A_{\hat r_h}^m}
		\bigl(v-\hat r_h\bigr)
		\varphi_\Gamma
		\bigl(v-\hat r_h\bigr).\]
	Therefore the third integral becomes
	\[\int_{\mathbb R^d}
		g(u_h(y)+v)
		\Bigl(
		\mathbf 1_{A_{\hat r_h}^m}(v)\varphi_\Gamma(v)
		-
		\Bigl(
		\bigl[
		\mathbf 1_{A_{\hat r_h}^m}(v)
		\varphi_\Gamma(v)
		\bigr]
		\wedge
		\bigl[
		\mathbf 1_{A_{\hat r_h}^m}
		\bigl(v-\hat r_h\bigr)
		\varphi_\Gamma
		\bigl(v-\hat r_h\bigr)
		\bigr]
		\Bigr)
		\Bigr)
		\,\mathrm dv.\]
	Hence we obtain
	\[\mathbb E g(\hat Y_h(x,y,Z))=
		\int_{(A_{\hat r_h}^m)^c}
		g(u_h(y)+w)\varphi_\Gamma(w)\,\mathrm dw +
		\int_{A_{\hat r_h}^m}
		g(u_h(y)+v)\varphi_\Gamma(v)\,\mathrm dv =
		\mathbb E g(u_h(y)+\Sigma Z).\]
	Since this holds for every bounded measurable \(g\),
	\[\mathcal L(\hat Y_h(x,y,Z))=\mathcal L(u_h(y)+\Sigma Z)=\mathcal L(\hat X_h(y,Z)).\]
	The proof is complete.
\end{proof}

	\begin{remark}
			The proof of Lemma \ref{coup} relies on the orthogonality of the reflection
			operator and the invariance of the Gaussian law \(N(0, h\Gamma)\) under the \(\Gamma^{-1}\)-reflection. This is why
			the reflected noise \(R_{\hat r_h}^\Gamma \Sigma Z\) has the same distribution as
			\(\Sigma Z\).
		\end{remark}

\begin{remark}
	The coupling constructed above is intrinsically discrete in time and is
	adapted to the one-step transition kernel of the TEM scheme. It is therefore
	not obtained by directly transferring a continuous-time reflection coupling to
	the numerical level. Instead, the coupling is applied after the nonlinear
	one-step map \(x\mapsto u_h(x)=\pi_h(x)+h b(\pi_h(x))\)
	has been performed. This point is important for the present analysis, because
	both the projection \(\pi_h\) and the nonlinear drift \(b\) enter the transition
	mechanism of the scheme. In this respect, the construction is specifically
	designed for discretizations of nonlinear SDEs, rather than for
	the standard Euler--Maruyama transition.
	
	The construction is also distinct from the reflection-type coupling developed
	in \cite{Majka}, which is formulated in a linear discrete-time setting with a
	scalar constant diffusion coefficient. That framework provides an important
	prototype for discrete reflection couplings with noise truncation. In contrast,
	the present coupling is adapted to the nonlinear one-step map of the TEM scheme
	and allows for a general non-degenerate constant diffusion matrix. In particular, the noise covariance is \(h\Gamma\), where \(\Gamma=\Sigma\Sigma^{\mathrm T}>0\),
	and the reflection is taken with respect to the
	\(\Gamma^{-1}\)-geometry. The relevant reflection direction and projection are
	therefore determined by the metric induced by \(\Gamma^{-1}\), rather than by
	the Euclidean geometry associated with a scalar or isotropic diffusion
	coefficient. This matrix-valued setting is essential for preserving the
	Gaussian marginal law under reflection. Moreover,
	while the truncation in \cite{Majka} is imposed on the magnitude of the full
	noise increment, here the truncation acts only on the component of \(\Sigma Z\)
	along the one-step separation direction, namely on \(|\lambda_{\hat r_h}(\Sigma Z)|\).	The directional truncation has two advantages for the subsequent estimates.
	First, it is more local to the coupling mechanism, since the Euclidean distance
	between the two coupled components depends only on the noise projected along
	\(\hat r_h\). Secondly, it reduces the key estimates to one-dimensional
	Gaussian integrals along this coupling direction, while the full noise remains
	\(d\)-dimensional and has covariance \(h\Gamma\). This feature is particularly
	useful in high-dimensional settings and in the analysis of nonlinear projected
	numerical schemes. 
\end{remark}
	
We establish the equality $\E|\hat{X}_h(x,Z)-\hat{Y}_h(x, y,Z)| = |u_h(x)-u_h(y)|$ in the following lemma—a property pivotal to subsequent analysis.
\begin{lemma}\label{1mom}
	For any $h, m>0$ and $x, y\in\R^d$, the coupling $(\hat{X}_h(x,Z), \hat{Y}_h(x, y,Z))$ defined by (\ref{hY}) satisfies $\E|\hat{X}_h(x,Z)-\hat{Y}_h(x, y,Z)| = |u_h(x)-u_h(y)|$.
\end{lemma}

\begin{proof}
For simplicity, we define \(\hat{R}_h(x, y):=\hat{X}_h(x, Z)-\hat{Y}_h(x, y,Z)\).
For short, we write $\hat{R}_h = \hat{R}_h(x, y)$ and $|\hat{R}_h| = |\hat{R}_h(x, y)|$.
If \(\hat r_h=\mathbf 0\) or \(|\hat r_h|>2R\), then the construction is synchronous. Hence \(\mathbb E|\hat{R}_h|=|\hat r_h|\).	
	
It remains to consider the case \(\hat r_h\neq\mathbf 0\) and \(|\hat r_h|\leqslant 2R\). On the event
\[\left\{\zeta\leqslant v_{\hat r_h}^m(\Sigma Z),
	\ \Sigma Z\in A_{\hat r_h}^m\right\},\]
the two components coincide, and hence \(|\hat{R}_h|=0\).
On the event
\[\left\{\zeta
	>v_{\hat r_h}^m(\Sigma Z),\ \Sigma Z\in A_{\hat r_h}^m
	\right\},\]
we have
\[\hat{R}_h
=\hat r_h+\Sigma Z-\mathcal R_{\hat r_h}^\Gamma\Sigma Z=
\left(1+2\frac{\langle \hat r_h,\Sigma Z\rangle_\Gamma}{
|\hat r_h|_\Gamma^2}\right)\hat r_h.\]
Therefore
\begin{equation}\label{refR}
|\hat{R}_h|=\left|
	|\hat r_h|+2\lambda_{\hat r_h}(\Sigma Z)\right|.
\end{equation}
On the event \(\{\Sigma Z\notin A_{\hat r_h}^m\}\), the coupling is
synchronous and \(|\hat{R}_h|
=|\hat r_h|\).
Therefore, we get
	\[\mathbb E|\hat{R}_h|=
		\int_{(A_{\hat r_h}^m)^c}
		|\hat r_h|\varphi_\Gamma(w)\,\mathrm dw+
		\int_{A_{\hat r_h}^m}
		\left|
		|\hat r_h(x,y)|
		+
		2\lambda_{\hat r_h}(w)
		\right|
		\left(
		1-v_{\hat r_h}^m(w)
		\right)
		\varphi_\Gamma(w)\,\mathrm dw .\]
	Subtracting \(|\hat r_h|\) from both sides gives
	\[
	\begin{aligned}
	\mathbb E|\hat{R}_h|
		-|\hat r_h|&=
		\int_{A_{\hat r_h}^m}
		\left(
		\left|
		|\hat r_h|
		+
		2\lambda_{\hat r_h}(w)
		\right|
		-
		|\hat r_h|
		\right)
		\varphi_\Gamma(w)\,\mathrm dw-
		\int_{A_{\hat r_h}^m}
		\left|
		|\hat r_h|
		+
		2\lambda_{\hat r_h}(w)
		\right|
		v_{\hat r_h}^m(w)
		\varphi_\Gamma(w)\,\mathrm dw .
	\end{aligned}
	\]	
	Since \(\Sigma Z\sim N(\mathbf 0,h\Gamma)\), \(\lambda_{\hat r_h}(\Sigma Z)
	\sim
	N\left(
	0,
	h|\hat r_h|^2/|\hat r_h|_\Gamma^2
	\right)\).
	Moreover, \(\lambda_{\hat r_h}
	\bigl(\hat r_h+w\bigr)
	=
	|\hat r_h|
	+
	\lambda_{\hat r_h}(w)\),
	and
	\[\frac{\varphi_\Gamma\bigl(\hat r_h+w\bigr)}{\varphi_\Gamma(w)}=
		\frac{p_{0,h|\hat r_h|^2/|\hat r_h|_\Gamma^2}
			\bigl(|\hat r_h|+\lambda_{\hat r_h}(w)
			\bigr)}{p_{0,h|\hat r_h|^2/|\hat r_h|_\Gamma^2}
			\bigl(\lambda_{\hat r_h}(w)\bigr)},\]
	where \(p_{0,\tau^2}\) denotes the density of \(N(0,\tau^2)\).
	Therefore the \(d\)-dimensional integral can be evaluated with respect to this
	one-dimensional Gaussian law, and we obtain
	\[
	\begin{aligned}
		&\mathbb E|\hat{R}_h|
		-|\hat r_h|\\
		=&
		\int_{\R}
		\left(
		\left|
		|\hat r_h|+2s
		\right|
		-
		|\hat r_h|
		\right)1_{[-m,m]}
		p_{0,h|\hat r_h|^2/|\hat r_h|_\Gamma^2}(s)
		\,\mathrm ds                                      \\
		&-
		\int_{\mathbb R}
		\left|
		|\hat r_h|+2s
		\right|\Bigl(
		\bigl[
		\mathbf 1_{[-m,m]}(s)
		p_{0,h|\hat r_h|^2/|\hat r_h|_\Gamma^2}(s)
		\bigr]\wedge
		\bigl[
		\mathbf 1_{[-m,m]}
		\bigl(|\hat r_h|+s\bigr)
		p_{0,h|\hat r_h|^2/|\hat r_h|_\Gamma^2}
		\bigl(|\hat r_h|+s\bigr)
		\bigr]
		\Bigr)
		\,\mathrm ds .
	\end{aligned}
	\]	
	The density \(p_{0,h|\hat r_h|^2/|\hat r_h|_\Gamma^2}\)
	is even and decreasing on \([0,\infty)\). Hence
	\[p_{0,h|\hat r_h|^2/|\hat r_h|_\Gamma^2}
		\bigl(|\hat r_h|+s\bigr)
		\leqslant
		p_{0,h|\hat r_h|^2/|\hat r_h|_\Gamma^2}(s)\Longleftrightarrow
		s\geqslant -\frac{|\hat r_h|}{2}.\]
	It follows that the minimum in the last display equals
	\[
	\begin{cases}
		\mathbf 1_{[-m,m]}(s)
		p_{0,h|\hat r_h|^2/|\hat r_h|_\Gamma^2}(s),
		&
		s<-\dfrac{|\hat r_h|}{2},\\[2mm]
		\mathbf 1_{[-m,m]}\bigl(|\hat r_h|+s\bigr)
		p_{0,h|\hat r_h|^2/|\hat r_h|_\Gamma^2}
		\bigl(|\hat r_h|+s\bigr),
		&
		s\geqslant -\dfrac{|\hat r_h|}{2}.
	\end{cases}
	\]
	Therefore
	\[
	\begin{aligned}
		&\mathbb E|\hat{R}_h|
		-|\hat r_h| \\
		=&
		-
		|\hat r_h|
		\int_{-\infty}^{-|\hat r_h|/2}
		\mathbf 1_{[-m,m]}(s)
		p_{0,h|\hat r_h|^2/|\hat r_h|_\Gamma^2}(s)
		\,\mathrm ds+
		2
		\int_{-|\hat r_h|/2}^{\infty}
		s\,\mathbf 1_{[-m,m]}(s)
		p_{0,h|\hat r_h|^2/|\hat r_h|_\Gamma^2}(s)
		\,\mathrm ds\\
		&-
		\int_{-|\hat r_h|/2}^{\infty}
		\bigl(|\hat r_h|+2s\bigr)
		\mathbf 1_{[-m,m]}
		\bigl(|\hat r_h|+s\bigr)
		p_{0,h|\hat r_h|^2/|\hat r_h|_\Gamma^2}
		\bigl(|\hat r_h|+s\bigr)
		\,\mathrm ds .
	\end{aligned}
	\]
	In the last integral, set \(u=|\hat r_h|+s\).
	Then
	\[
	\begin{aligned}
		&\int_{-|\hat r_h|/2}^{\infty}
		\bigl(|\hat r_h|+2s\bigr)
		\mathbf 1_{[-m,m]}
		\bigl(|\hat r_h|+s\bigr)
		p_{0,h|\hat r_h|^2/|\hat r_h|_\Gamma^2}
		\bigl(|\hat r_h|+s\bigr)
		\,\mathrm ds \\
		=&
		\int_{|\hat r_h|/2}^{\infty}
		\bigl(2u-|\hat r_h|\bigr)
		\mathbf 1_{[-m,m]}(u)
		p_{0,h|\hat r_h|^2/|\hat r_h|_\Gamma^2}(u)
		\,\mathrm du .
	\end{aligned}
	\]
	Hence, we obtain
	\[
	\begin{aligned}
		\mathbb E|\hat{R}_h|
		-|\hat r_h| =&
		2
		\int_{-|\hat r_h|/2}^{|\hat r_h|/2}
		s\,\mathbf 1_{[-m,m]}(s)
		p_{0,h|\hat r_h|^2/|\hat r_h|_\Gamma^2}(s)
		\,\mathrm ds+
		|\hat r_h|
		\int_{|\hat r_h|/2}^{\infty}
		\mathbf 1_{[-m,m]}(s)
		p_{0,h|\hat r_h|^2/|\hat r_h|_\Gamma^2}(s)
		\,\mathrm ds                                     \\
		&-
		|\hat r_h|
		\int_{-\infty}^{-|\hat r_h|/2}
		\mathbf 1_{[-m,m]}(s)
		p_{0,h|\hat r_h|^2/|\hat r_h|_\Gamma^2}(s)
		\,\mathrm ds .
	\end{aligned}
	\]
	The first integral vanishes because the integrand is odd. The last two
	integrals cancel because both \(\mathbf 1_{[-m,m]}(s)\) and \(p_{0,h|\hat r_h|^2/|\hat r_h|_\Gamma^2}(s)\) are even functions. Hence \(\mathbb E|\hat{R}_h|
	-|\hat r_h| = 0\). The proof is complete.
\end{proof}
			
			\subsection{Contractivity of the coupling}\label{Sec4.3}
			In this section, our primary objective is to analyze the contractivity of the coupling random vector $(\hat{X}_h(x,Z), \hat{Y}_h(x, y,Z))$ defined by (\ref{hY}), with respect to the initial positions $(\pi_h(x), \pi_h(y))$ with an appropriate function.
			
			For any $h>0$ and  $x, y\in\R^d$, we recall $|\hat{r}_h|$, $|\hat{R}_h|$ and define $|r_h|$ by
			\begin{subequations}\label{eq3}
				\begin{align}
					|r_h|: &= |r_h(x, y)| = |\pi_h(x)-\pi_h(y)|,  \label{Nata}\\
					|\hat{r}_h|: &= |\hat{r}_h(x, y)| = |u_h(x)-u_h(y)|, \label{Natb}\\
					|\hat{R}_h|: &= |\hat{R}_h(x, y)| = |\hat{X}_h(x,Z)-\hat{Y}_h(x, y,Z)|.\label{Natc}
				\end{align}
			\end{subequations}	
			According to Lemma \ref{1mom}, for any $h>0$ and $x, y\in\R^d$, we have
			\begin{equation}\label{cc1}
				\mathbb{E}(|\hat{R}_h| - |r_h|) = |\hat{r}_h| - |r_h|.
			\end{equation}
			If $\pi_h(x) = \pi_h(y)$, we find that $|\hat{r}_h| = 0$. Thus, according to (\ref{cc1}), we know that \(\mathbb{E}(|\hat{R}_h|) = 0\). For this case, it is evident that the contractivity of the coupling random vector $(\hat{X}_h(x, Z), \hat{Y}_h(x, y,Z))$ is trivial. Therefore, we  primarily focus on the contractivity of coupling $(\hat{X}_h(x, Z), \hat{Y}_h(x, y,Z))$ for $\pi_h(x) \neq \pi_h(y)$.  
We are going to construct an appropriate function $f:[0, +\infty) \to [0, +\infty)$ satisfying $f' > 0$ and $f'' < 0$ and then obtain the contractivity of $(\hat{X}_h(x, Z), \hat{Y}_h(x, y,Z))$ under the function $f$, that is, there exist constants $c, h^*>0$ such that for any $x, y\in\R^d$ and $h \in (0, h^*)$,
\begin{equation*}
\mathbb{E}f(|\hat{R}_h|) \leqslant (1 - ch)f(|r_h|).
\end{equation*}

Let \(\vartheta_-:=\sqrt{\lambda_{\min}(\Gamma)}\) and \(\vartheta_+:=\sqrt{\lambda_{\max}(\Gamma)}\). Define
\begin{equation}\label{c0}
\alpha_{1, \Gamma} := \frac1{16}
\left(1-\exp\left\{-\frac{1}{8\vartheta_+^2}\right\}\right)\inf_{\vartheta\in[\vartheta_-,\vartheta_+]}\int_{1/4}^{3/8}p_{0,\vartheta^2}(u)\,\mathrm du .
\end{equation}
Define 				
\begin{equation}\label{varphi}
\varphi (u) = \mathrm{exp}\left(-\frac{2L+1}{2\alpha_{1, \Gamma}}\left(u^2+2u\right)\right), \quad\quad\Phi(u)= \int_0^u \varphi(s) \md s, \quad u\geqslant 0.
\end{equation}	
Clearly, the function $\varphi$ is decreasing  
while $\Phi(u)$ is increasing for $u>0$ due to the nonnegativity of $\varphi$. 
Define 			
\begin{equation}\label{Ac} 
r_1 = 2R+16 ,~~~~
\alpha_{2, \Gamma}= \frac{1}{2}\left(\int_0^{r_1}\frac{\Phi(s+1)}{\alpha_{1, \Gamma}\varphi(s)}\md s\right)^{-1}.
\end{equation} 
Hence, owing to the monotonicity of $\Phi$ and $\varphi$, we have 
\begin{equation}\label{dd6}
\frac{\alpha_{1, \Gamma}\varphi(r_1)}{2r_1(r_1+1)}\leqslant \frac{1}{2}\left(\int_0^{r_1}\frac{r_1+1}{\alpha_{1, \Gamma}\varphi(r_1)}\md s\right)^{-1}\leqslant \alpha_{2, \Gamma}\leqslant \frac{1}{2}\left(\int_0^{r_1}\frac{\Phi(1)}{\alpha_{1, \Gamma}}\md s\right)^{-1} =
 \frac{\alpha_{1, \Gamma}}{2r_1\Phi(1)}.
\end{equation}
We now proceed to construct the appropriate function $f$. For any $u\geqslant 0$, we define 
\begin{equation}\label{f}
\rho(u) = 1- \alpha_{2, \Gamma}\int_0^u \frac{\Phi(s+1)}{\alpha_{1, \Gamma}\varphi(s)}\md s,  \quad f(u) = \int_0^u \varphi (s\wedge r_1)\rho(s\wedge r_1) \md s.
\end{equation}
Noticing that for any $u\in[0, r_1]$,  $1/2\leqslant \rho(u) \leqslant 1$, we know that  
\begin{equation}\label{f1}
\frac{1}{2}\varphi(r_1)u\leqslant \frac{1}{2}\int_0^u\varphi(s\wedge r_1) \md s \leqslant f(u)\leqslant \Phi(u)\leqslant u, ~~~~\forall~ u\geqslant 0.
\end{equation}
To preserve the ergodicity of our numerical scheme, we select the suitable parameter $M$ satisfying
\begin{equation}\label{M1}
M \geqslant
\max\left\{|b(\mathbf{0})|, \frac{\sqrt{\varphi(r_1)RK}}{8\sqrt{r_1}}, \frac{\sqrt{\alpha_{1, \Gamma}}}{4\sqrt{2r_1\Phi(1)}}, 1\right\} 
\end{equation} in the truncation mapping $\pi_h$ defined by \eqref{TF}.

\begin{remark} It is worth noting that   \(M\geqslant \max\left\{|b(\mathbf{0})|, \sqrt{K}, \sqrt{L+1}, 1\right\}\) implies that \eqref{M1} holds.
In fact, noting that $\varphi(u)\leqslant 1$ for all $u>0$ and \eqref{Ac}, we know that
\(\sqrt{\varphi(r_1)RK}/(8\sqrt{r_1})<\sqrt{K}\).
By using inequality \(\mathrm{e}^{-u}\leqslant 1/(1+u)\) for all $u>0$, we have
\begin{align*}
\Phi(1) &= \int_0^1 \mathrm{exp}\left(-\frac{2L+1}{2\alpha_{1, \Gamma}}\left(u^2+2u\right)\right)\md u\geqslant \int_0^1 \mathrm{exp}\left(-\frac{6L+3}{2\alpha_{1, \Gamma}}u\right)\md u> \frac{2\alpha_{1, \Gamma}}{2\alpha_{1, \Gamma}+6L+3},
\end{align*}
which together with \eqref{c0} and \eqref{Ac} implies that
\(\sqrt{\alpha_{1, \Gamma}}/(4\sqrt{2r_1\Phi(1)}) < \sqrt{L+1}\). Then the desired assertion follows.
\end{remark}

Define
\begin{equation}\label{c1}
	\alpha_{3,\Gamma}
	:=
	\inf_{\vartheta\in[\vartheta_-,\vartheta_+]}
	\frac{5}{\vartheta^3\sqrt{2\pi}}
	\min\left\{
	\mathrm e^{-1/(2\vartheta^2)},
	8\mathrm e^{-32/\vartheta^2}
	\right\}
\end{equation}
and 
\begin{equation}\label{c2}
	\alpha_{4,\Gamma}
	:=
	\inf_{\vartheta\in[\vartheta_-,\vartheta_+]}
	4\vartheta^3
	\left(1-\mathrm e^{-1/\vartheta^2}\right)
	\int_0^{1/(2\vartheta)}
	u^3p_{0,1}(u)\,\mathrm du .
\end{equation}
Since \(\Gamma>0\), one has \(0<\vartheta_-\leqslant \vartheta_+<\infty\), and hence \(\alpha_{1, \Gamma}, \alpha_{3,\Gamma},\alpha_{4,\Gamma}>0\).
Moreover, when \(\Gamma=\Sigma^2\mathbf I\), we have
\(\vartheta_-=\vartheta_+=|\Sigma|\), and \(\alpha_{1, \Gamma}, \alpha_{3,\Gamma}\), \(\alpha_{4,\Gamma}\) reduce to
\[\frac1{16}
\left(1-\mathrm e^{-1/(8\Sigma^2)}\right)\int_{1/4}^{3/8}p_{0,\Sigma^2}(u)\,\mathrm du<\frac{3\sqrt{3}}{16\mathrm{e}^{3/2}\sqrt{2\pi}}\approx 0.03,\]
\[\frac{5}{|\Sigma|^3\sqrt{2\pi}}
\min\left\{
\mathrm e^{-1/(2\Sigma^2)},
8\mathrm e^{-32/\Sigma^2}
\right\}<\frac{15\sqrt{3}}{\mathrm{e}^{3/2}\sqrt{2\pi}}\approx 2.31,\]
and
\[4|\Sigma|^3
\left(1-\mathrm e^{-1/\Sigma^2}\right)
\int_0^{1/(2|\Sigma|)}
u^3p_{0,1}(u)\,\mathrm du<\frac{1}{\sqrt{2\pi}}\approx 0.4.\]
Now we establish the main result that the coupling random vector $(\hat{X}_h(x,Z), \hat{Y}_h(x, y,Z))$ defined by \eqref{hY} is globally contractive. Define \begin{equation*}
h_0 = \min\left\{\left(KM^{-2}\right)^{1/(1-2\bar{\theta})}, K^{-1}, L^{-1}\right\}.
\end{equation*}

We first derive a simple deterministic estimate for the drift-updated distance
\(|\hat r_h|\). It captures the contraction induced by the dissipativity condition outside the radius \(R\), and provides a uniform control inside this
radius.

\begin{lemma}\label{Le4.8}
		Let Assumption \ref{A} hold. Then,
		for any \(h\in(0,h_0)\), if \(|r_h|>R\), then
		\begin{equation}\label{Sub4eq1}
			|\hat r_h|\leqslant |r_h|\left(1-\frac{hK}{2}\right).
		\end{equation}
		If, in addition, Assumption \ref{A1} holds, then for any \(h\in(0,h_0)\),
		whenever \(|r_h|\leqslant R\), one has
		\begin{equation}\label{Sub4eq2}
			|\hat r_h|\leqslant 2R.
		\end{equation}
	\end{lemma}
\begin{proof} 
	For any $h\in(0, h_0)$, we have
	\begin{equation}\label{cc2}
		M^2h^{2-2\theta}\leqslant M^2h^{2-2\bar{\theta}}\leqslant hK\leqslant 1,\quad \forall \theta\in(0, \bar{\theta}]\quad \text{and} \quad 1+2hL+M^2h^{2-2\theta}\leqslant 4.
	\end{equation}
	Using  (\ref{GL}) one notices from the definition of $|\hat{r}_h|$ that 
	\begin{equation}\label{eq+6}
		\begin{aligned}
			|\hat{r}_h|
			&=\sqrt{|r_h|^2+\!2h\left\langle\pi_h(x)-\pi_h(y), b(\pi_h(x))-b(\pi_h(y))\right\rangle+h^2|b(\pi_h(x))-b(\pi_h(y))|^2}\\
			&\leqslant \sqrt{(1+M^2h^{2-2\theta})|r_h|^2+\!2h\left\langle\pi_h(x)-\pi_h(y), b(\pi_h(x))-b(\pi_h(y))\right\rangle }.
		\end{aligned}
	\end{equation}
	If $|r_h|>R$, by virtue of \eqref{CaI} and \eqref{cc2}, it follows from \eqref{eq+6} that
	\begin{equation*}
		|\hat{r}_h| \leqslant  |r_h|\sqrt{1+M^2h^{2-2\theta}-2hK}\leqslant |r_h|\sqrt{1-hK}.
	\end{equation*}
	Then, this together with the fact $\sqrt{1+x}\leqslant 1+x/2$ for all $x\geqslant-1$   implies that \eqref{Sub4eq1} holds. On the other hand, if $|r_h|\leqslant R$, by virtue of  (\ref{CaI}) and \eqref{cc2}, it follows from \eqref{eq+6} that
	\begin{equation}\label{ss}
		\begin{aligned}
			|\hat{r}_h|	\leqslant |r_h|\sqrt{1+M^2h^{2-2\theta}+2hL}\leqslant 2R.
		\end{aligned}
	\end{equation}
	The proof is complete.
\end{proof}

With the concave distance function (f), the constants introduced above, and the deterministic estimate in Lemma \ref{Le4.8} at hand, we now state the global one-step contractivity estimate for the coupling defined in \eqref{hY}.

\begin{theorem}\label{Th1}
Let Assumptions \ref{A} and \ref{A1} hold.
For any $x, y\in\R^d$ and \(h \in (0, \bar{h}]\), consider the coupling $(\hat{X}_h(x,Z), \hat{Y}_h(x,y,Z))$ defined by \eqref{hY}, with the parameter $m=8$. Then, we have
\begin{equation*}
\E f(|\hat{R}_h|)\leqslant (1-ch)f(|r_h|),
\end{equation*}
where $|r_h|, |\hat{R}_h|$ are defined in (\ref{eq3}), function $f$ is defined by \eqref{f}, 
the constant $\bar{h}$ is defined by
\begin{equation}\label{bar}
\bar{h}=\min\left\{h_0, \left(\frac{\alpha_{2, \Gamma}\alpha_{3, \Gamma}\Phi(1)}{8M\alpha_{1, \Gamma}}\right)^{2/(1-2\bar{\theta})}, \left(\frac{r_1}{18}\right)^2, R^2, (4M)^{-2}, (4Mr_1)^{2/(2\bar{\theta}-1)}, M^{1/(2\bar{\theta}-1)}\right\},
\end{equation}
and the contractive rate $c$ is given by 
\begin{equation*}
c = \min\left\{M, \frac{\varphi(r_1)KR }{4r_1}, \frac{\alpha_{2, \Gamma}\alpha_{4, \Gamma}\Phi(1)}{4\alpha_{1, \Gamma}}, \alpha_{2, \Gamma}\right\}.
\end{equation*}
\end{theorem}

Theorem \ref{Th1} extends the one-step contractivity estimate of \cite[Theorem 2.1]{Majka} to the present nonlinear setting with a superlinearly growing drift and a constant non-degenerate diffusion matrix. Its proof has to accommodate both the local non-contractivity and the dissipativity at infinity of the drift in \eqref{E}, as well as the truncation mechanism of the TEM scheme. The estimate obtained in Theorem \ref{Th1} will serve as a key ingredient in the subsequent derivation of the uniform-in-time consistency estimates for the TEM approximation.

Because the proof of Theorem \ref{Th1} is rather technical, we begin by
outlining its structure; see Figure \ref{Figtu} for a roadmap. The proof is built on two
main preparatory ingredients: the construction of a suitable concave distance
function \(f\), and the one-step coupling estimate given in Lemma \ref{1mom}.
With these tools in hand, we divide the argument according to the size of the
current distance \(|r_h|\). More precisely, we consider three cases:
\(0<|r_h|\leqslant \sqrt h\), \(\sqrt h<|r_h|\leqslant r_1\), and
\(|r_h|>r_1\). The first and third cases are treated in
Propositions \ref{PR1} and \ref{PR5}, respectively. The intermediate regime
\(\sqrt h<|r_h|\leqslant r_1\) is the most delicate one, and it is further
split according to the size of the one-step updated distance
\(|\widehat r_h|\), which leads to three subcases handled separately in the
proof of Proposition \ref{PR4}. More specifically, the subcases
\(|\widehat r_h|\leqslant \sqrt h\) and
\(\sqrt h<|\widehat r_h|\leqslant 2R\) are treated via
Propositions \ref{PR2} and \ref{PR3}, while the case
\(|\widehat r_h|>2R\) is controlled by synchronous coupling together with
inequality \eqref{I-3}. Since these three propositions cover all possible
values of \(|r_h|\), their combination yields the desired estimate in Theorem
\ref{Th1}.
\begin{figure}[htbp]
	\centering
	\resizebox{\textwidth}{!}{%
		\begin{tikzpicture}[node distance=8mm and 8mm]
			\node[keybox] (goal) {Goal: Theorem \ref{Th1}\\\(\mathbb{E}f(|\widehat R_h|)\leqslant (1-ch)f(|r_h|)\)};
			\node[wbox, above left=6mm and 5mm of goal] (prep1) {Construct a concave distance function \(f\)};
			\node[wbox, above right=6mm and 5mm of goal] (prep2) {one-step coupling\\Lemma \ref{1mom}: first-order identity};
			
			\node[rbox, below left=8mm and 10mm of goal] (case1) {Case 1: \(0<|r_h|\leqslant \sqrt h\)\\Lemma \ref{Le9}\\ \(\Rightarrow\) Proposition \ref{PR1}};
			\node[rbox, below=8mm of goal] (case2) {Case 2: \(\sqrt h<|r_h|\leqslant r_1\)\\further split by \(|\widehat r_h|\)\\\(\Rightarrow\) Proposition \ref{PR4}};
			\node[rbox, below right=8mm and 10mm of goal] (case3) {Case 3: \(|r_h|>r_1\)\\Lemma \ref{Le4.8} + concavity of \(f\)\\\(\Rightarrow\) Proposition \ref{PR5}};
			
			\node[tinybox, below left=8mm and 10mm of case2] (sub1) {2-I: \(|\widehat r_h|\leqslant \sqrt h\)\\Lemma \ref{Le10} \(\Rightarrow\) Proposition \ref{PR2}};
			\node[tinybox, below=8mm of case2] (sub2) {2-II: \(\sqrt h<|\widehat r_h|\leqslant 2R\)\\Lemma \ref{2mom} \(\Rightarrow\) Proposition \ref{PR3}};
			\node[tinybox, below right=8mm and 10mm of case2] (sub3) {2-III: \(|\widehat r_h|>2R\)\\synchronous coupling + inequality \eqref{I-3}};
			
			\node[keybox, below=25mm of case2] (combine) {Propositions \ref{PR1}, \ref{PR4}, and \ref{PR5} cover all values of \(|r_h|\).\\Together they imply Theorem \ref{Th1}.};
			
			\draw[arr] (prep1) -- (goal);
			\draw[arr] (prep2) -- (goal);
			\draw[arr] (goal) -- (case1);
			\draw[arr] (goal) -- (case2);
			\draw[arr] (goal) -- (case3);
			\draw[arr] (case2) -- (sub1);
			\draw[arr] (case2) -- (sub2);
			\draw[arr] (case2) -- (sub3);
			\draw[arr] (sub1) -- (combine);
			\draw[arr] (sub2) -- (combine);
			\draw[arr] (sub3) -- (combine);
			\draw[arr] (case1) -- (combine);
			\draw[arr] (case3) -- (combine);
		\end{tikzpicture}%
	}
	\caption{Proof roadmap for Theorem \ref{Th1}.}\label{Figtu}
\end{figure}

To prove that Theorem \ref{Th1} holds in \underline{\bf{Case 1. \textit{ $|r_h|\in(0, \sqrt{h}]$}}}, we prepare the following lemma.

\begin{lemma}\label{Le9}
	Let Assumption \ref{A} hold. For any \(x,y\in\mathbb R^d\) and
	\(h\in(0,M^{-2}]\), the coupling \((\hat X_h(x,Z),\hat Y_h(x,y,Z))\) defined by \eqref{hY} has the property that 
	\[\mathbb E\left[
		\bigl(|\hat R_h|-|r_h|\bigr)^2
		\mathbf 1_{(|r_h|+\sqrt h,\,|r_h|+17\sqrt h)}
		(|\hat R_h|)
		\right]\geqslant
		\alpha_{3,\Gamma}
		|\hat r_h|\sqrt h\,
		\mathbf 1_{(0,\sqrt h]}(|r_h|)
		\mathbf 1_{[0,2R]}(|\hat r_h|),\]
		where $|r_h|, |\hat{r}_h|, |\hat{R}_h|$ are defined in \eqref{eq3} and $\alpha_{3, \Gamma}$ is defined by \eqref{c1}.
\end{lemma}

\begin{proof}
	If \(\hat r_h=\mathbf 0\), then both sides are zero, and there is nothing
	to prove. We therefore assume that \(\hat r_h\neq \mathbf 0\).
	
	Since on the event \(|\hat R_h|\in(|r_h|+\sqrt h,\,|r_h|+17\sqrt h)\), one has \((|\hat R_h|-|r_h|)^2\geqslant h\),
	we get
	\begin{equation}\label{matrix-J6-start}
			\mathbb E\left[
			\bigl(|\hat R_h|-|r_h|\bigr)^2
			\mathbf 1_{(|r_h|+\sqrt h,\,|r_h|+17\sqrt h)}
			(|\hat R_h|)
			\right]\geqslant
			h\,
			\mathbb E\left[
			\mathbf 1_{(|r_h|+\sqrt h,\,|r_h|+17\sqrt h)}
			(|\hat R_h|)
			\right].
	\end{equation}
	It is enough to estimate from below the probability on the right-hand side
	coming from the reflection part of the coupling.
	On the reflection event, by \eqref{refR}, \(|\hat R_h|
	=
	\left|
	|\hat r_h|
	+
	2\lambda_{\hat r_h}(\Sigma Z)
	\right|\). By the definition of \(v_{\hat r_h}^m\), we obtain
	\begin{equation}\label{matrix-J6}
		\begin{aligned}
			&\mathbb E\left[
			\mathbf 1_{(|r_h|+\sqrt h,\,|r_h|+17\sqrt h)}
			(|\hat R_h|)
			\right]                                            \\
			\geqslant&
			\int_{\mathbb R^d}
			\left(
			\mathbf 1_{A_{\hat r_h}^m}(w)\varphi_\Gamma(w)
			-
			\mathbf 1_{A_{\hat r_h}^m}
			\bigl(\hat r_h+w\bigr)
			\varphi_\Gamma
			\bigl(\hat r_h+w\bigr)
			\right)^+                                          \\
			&\qquad\qquad\times
			\mathbf 1_{[0,2R]}(|\hat r_h|)
			\mathbf 1_{(|r_h|+\sqrt h,\,|r_h|+17\sqrt h)}
			\left(
			\left|
			|\hat r_h|
			+
			2\lambda_{\hat r_h}(w)
			\right|
			\right)
			\,\mathrm dw .
		\end{aligned}
	\end{equation}
	Since \(\lambda_{\hat r_h}(\Sigma Z)
	\sim N(0,h|\hat r_h|^2/|\hat r_h|_\Gamma^2)\) and \(\lambda_{\hat r_h}
	\bigl(\hat r_h+w\bigr)
	=
	|\hat r_h|+\lambda_{\hat r_h}(w)\),
	the right-hand side of \eqref{matrix-J6} is bounded from below by
	\begin{equation}\label{matrix-J6-one-dim}
		\begin{aligned}
			&\int_{\mathbb R}
			\Bigl(
			\mathbf 1_{[-m,m]}(s)
			p_{0,h|\hat r_h|^2/|\hat r_h|_\Gamma^2}(s)
			-
			\mathbf 1_{[-m,m]}(|\hat r_h|+s)
			p_{0,h|\hat r_h|^2/|\hat r_h|_\Gamma^2}(|\hat r_h|+s)
			\Bigr)^+                                           \\
			&\qquad\qquad\times
			\mathbf 1_{[0,2R]}(|\hat r_h|)
			\mathbf 1_{(|r_h|+\sqrt h,\,|r_h|+17\sqrt h)}
			\bigl(|\,|\hat r_h|+2s\,|\bigr)
			\,\mathrm ds ,
		\end{aligned}
	\end{equation}
	where \(p_{0,h|\hat r_h|^2/|\hat r_h|_\Gamma^2}\) denotes the density of \(N(0,h|\hat r_h|^2/|\hat r_h|_\Gamma^2)\). We now estimate a lower bound for the indicator function $\1_{(|r_h|+\sqrt{h}, |r_h|+17\sqrt{h})}\left(|\hat{r}_h|+2s\right)$ in \eqref{matrix-J6-one-dim}. Using (\ref{GL}), we have
	\begin{align*}
		|\hat{r}_h| = |\pi_h(x)-\pi_h(y)+h(b(\pi_h(x))-b(\pi_h(y)))|\leqslant (1+Mh^{1-\theta})|r_h|,
	\end{align*}
	and similarly \(|r_h|= |\pi_h(x)-\pi_h(y)|
	\leqslant |\hat{r}_h|+Mh^{1-\theta}|r_h|\).
	Since \(h\in(0,M^{-2}]\) and \(M\geqslant 1\), this yields 
	\begin{equation}\label{rr}
		\left||r_h|-|\hat{r}_h|\right|\leqslant Mh^{1-\theta}|r_h|\leqslant M\sqrt{h}|r_h|\leqslant |r_h|.
	\end{equation}
	This implies that on the event \(\{|r_h|\in(0,\sqrt h]\}\),  
	$$\frac{1}{2}\left(|r_h|-|\hat{r}_h|+\sqrt{h}\right)\leqslant\sqrt{h}, \qquad \frac{1}{2}\left(|r_h|-|\hat{r}_h|+17\sqrt{h}\right) \geqslant8\sqrt{h}.$$
	Hence, on \(\{|r_h|\in(0,\sqrt h]\}\),
	\begin{equation}\label{matrix-indicator-lower}
		\begin{aligned}
			\mathbf 1_{(|r_h|+\sqrt h,\,|r_h|+17\sqrt h)}
			\bigl(|\,|\hat r_h|+2s\,|\bigr)\geqslant
			\mathbf 1_{(\sqrt h,\,8\sqrt h)}(s).
		\end{aligned}
	\end{equation}	
	Using \eqref{matrix-indicator-lower} in \eqref{matrix-J6-one-dim}, and taking
	\(m=8\), we obtain
	\begin{equation}\label{matrix-J6-lower1}
		\begin{aligned}
			\mathbb E\left[
			\mathbf 1_{(|r_h|+\sqrt h,\,|r_h|+17\sqrt h)}
			(|\hat R_h|)
			\right]\geqslant&
			\int_{\sqrt h}^{8\sqrt h}
			\Bigl(
			p_{0,h|\hat r_h|^2/|\hat r_h|_\Gamma^2}(s)
			-
			p_{0,h|\hat r_h|^2/|\hat r_h|_\Gamma^2}(|\hat r_h|+s)
			\Bigr)
			\,\mathrm ds\\
			&\qquad\qquad \times
			\mathbf 1_{(0,\sqrt h]}(|r_h|)
			\mathbf 1_{[0,2R]}(|\hat r_h|).
		\end{aligned}
	\end{equation}
	Here the positive part has been removed because, for \(s>0\), \(|\hat r_h|+s\geqslant s\), and the centered Gaussian density is decreasing on \([0,\infty)\). 
	On \(\{|r_h|\in(0,\sqrt h]\}\), \eqref{rr} gives \(0\leqslant |\hat r_h|/\sqrt{h}\leqslant2\).
	By the change of variables \(s=\sqrt h\,u\), \eqref{matrix-J6-lower1} yields
	\begin{equation*}
		\begin{aligned}
			&\mathbb E\left[
			\mathbf 1_{(|r_h|+\sqrt h,\,|r_h|+17\sqrt h)}
			(|\hat R_h|)
			\right]  \\
			\geqslant&
			\frac{|\hat r_h|_\Gamma}{|\hat r_h|\sqrt{2\pi}}
			\int_1^8
			\left[
			\exp\left(-\frac{|\hat r_h|_\Gamma^2u^2}{2|\hat r_h|^2}\right)
			-
			\exp\left(-\frac{|\hat r_h|_\Gamma^2(u+|\hat r_h|/\sqrt{h})^2}{2|\hat r_h|^2}\right)
			\right]\mathrm du
			\mathbf 1_{(0,\sqrt h]}(|r_h|)
			\mathbf 1_{[0,2R]}(|\hat r_h|).
		\end{aligned}
	\end{equation*}
	Since \(0\leqslant |\hat r_h|/\sqrt{h}\leqslant2\), the interval \((3-|\hat r_h|/\sqrt{h},8-|\hat r_h|/\sqrt{h})\) is contained in \((1,8)\).
	Therefore,
	\[
	\begin{aligned}
		&\int_1^8
		\left[
		\exp\left(-\frac{|\hat r_h|_\Gamma^2u^2}{2|\hat r_h|^2}\right)
		-
		\exp\left(-\frac{|\hat r_h|_\Gamma^2(u+|\hat r_h|/\sqrt{h})^2}{2|\hat r_h|^2}\right)
		\right]\mathrm du\\
		\geqslant&
		\int_{3-|\hat r_h|/\sqrt{h}}^{8-|\hat r_h|/\sqrt{h}}
		\int_u^{u+|\hat r_h|/\sqrt{h}}
		\frac{|\hat r_h|_\Gamma^2v}{|\hat r_h|^2}
		\exp\left(-\frac{|\hat r_h|_\Gamma^2v^2}{2|\hat r_h|^2}\right)
		\,\mathrm dv\,\mathrm du .
	\end{aligned}
	\]
	For \(u\in(3-|\hat r_h|/\sqrt{h},8-|\hat r_h|/\sqrt{h})\) and \(v\in[u,u+|\hat r_h|/\sqrt{h}]\), we have \(v\in[1,8]\). Hence
	\[
	v\exp\left(-\frac{|\hat r_h|_\Gamma^2v^2}{2|\hat r_h|^2}\right)
	\geqslant
	\min\left\{
	\mathrm e^{-|\hat r_h|_\Gamma^2/(2|\hat r_h|^2)},
	8\mathrm e^{-32|\hat r_h|_\Gamma^2/|\hat r_h|^2}
	\right\}.
	\]
	It follows that
	\[
	\begin{aligned}
		&\mathbb E\left[
		\mathbf 1_{(|r_h|+\sqrt h,\,|r_h|+17\sqrt h)}
		(|\hat R_h|)
		\right]\\
		\geqslant&
		\frac{5|\hat r_h|_\Gamma^3}
		{|\hat r_h|^2\sqrt{2\pi h}}
		\min\left\{
		\mathrm e^{-|\hat r_h|_\Gamma^2/(2|\hat r_h|^2)},
		8\mathrm e^{-32|\hat r_h|_\Gamma^2/|\hat r_h|^2}
		\right\}\mathbf 1_{(0,\sqrt h]}(|r_h|)
		\mathbf 1_{[0,2R]}(|\hat r_h|).
	\end{aligned}
	\]
	By the definition of \(\alpha_{3,\Gamma}\), and since \(|\hat r_h|/|\hat r_h|_\Gamma\in[\vartheta_-,\vartheta_+]\),
	we obtain
	\begin{equation}\label{matrix-J6-final}
\mathbb E\left[
			\mathbf 1_{(|r_h|+\sqrt h,\,|r_h|+17\sqrt h)}
			(|\hat R_h|)
			\right]\geqslant
			\frac{\alpha_{3,\Gamma}|\hat r_h|}{\sqrt h}
			\mathbf 1_{(0,\sqrt h]}(|r_h|)
			\mathbf 1_{[0,2R]}(|\hat r_h|).
	\end{equation}
	Combining \eqref{matrix-J6-start} and \eqref{matrix-J6-final}, we obtain the desired result. The proof is complete.
\end{proof}
			
We now invoke Lemmas \ref{1mom} and \ref{Le9} to establish Theorem \ref{Th1} for \underline{\bf{Case 1. \textit{ $|r_h|\in(0, \sqrt{h}]$}}}.
\begin{proposition}\label{PR1}
Let Assumptions \ref{A} and \ref{A1} hold. For any $x, y\in\R^d$ and \(h \in (0, \bar{h}]\), where $\bar{h}$ is defined by \eqref{bar}, consider the coupling $(\hat{X}_h(x,Z), \hat{Y}_h(x,y,Z))$ defined by (\ref{hY}), with $m=8$. Then, if $|r_h|\in(0, \sqrt{h}]$, we have
\begin{equation*}
\E f(|\hat{R}_h|)\leqslant (1-Mh)f(|r_h|),
\end{equation*}
where $|r_h|$ and $|\hat{R}_h|$ are defined by \eqref{Nata} and \eqref{Natc} respectively, function $f$ is given by \eqref{f}.
\end{proposition}
\begin{proof}
We use Taylor's formula with
integral remainder:
\begin{equation}\label{cc4}
\begin{aligned}
	\mathbb E\left[f(|\hat R_h|)-f(|r_h|)\right]=&
	f'(|r_h|)
	\mathbb E\bigl(|\hat R_h|-|r_h|\bigr)               \\
	&+
	\mathbb E\left[
	\bigl(|\hat R_h|-|r_h|\bigr)^2
	\int_0^1
	(1-\tau)
	f''\bigl(|r_h|+\tau(|\hat R_h|-|r_h|))
	\,\mathrm d\tau
	\right].
\end{aligned}	
\end{equation}
Here \(f''\) is understood in the a.e. sense, which is sufficient since \(f'\) is absolutely continuous.
By Lemma \ref{1mom} and $0<f'\leqslant 1$, we get
\begin{equation}\label{cc5}
f'(|r_h|)\E(|\hat{R}_h|-|r_h|) = f'(|r_h|)(|\hat{r}_h|-|r_h|)\leqslant ||\hat{r}_h|-|r_h||.
\end{equation} 
We next estimate the second-order term of \eqref{cc4}. Since  $h\leqslant R^2$ and $|r_h|\in(0, \sqrt{h}]$, we have $|r_h|\leqslant R$. Then, by (\ref{Sub4eq2}), we get $|\hat{r}_h|\leqslant 2R$. 
Moreover, since \(h\leqslant (r_1/18)^2\), \(|r_h|+17\sqrt h
\leqslant 18\sqrt h
\leqslant r_1\).
On the event \(|\hat R_h|\in(|r_h|+\sqrt h,\,|r_h|+17\sqrt h)\),
all points between \(|r_h|\) and \(|\hat R_h|\) lie in
\((0,18\sqrt h)\subset(0,r_1)\). Hence, using \(f''<0\) and Lemma \ref{Le9}, we may discard the
complement of this event and obtain
\begin{equation}\label{Tqq1}
\begin{aligned}
	&\mathbb E\left[
	\bigl(|\hat R_h|-|r_h|\bigr)^2
	\int_0^1
	(1-\tau)
	f''\bigl(|r_h|+\tau(|\hat R_h|-|r_h|)\bigr)
	\,\mathrm d\tau
	\right] \\
	\leqslant&
	\frac12
	\sup_{u\in(0,18\sqrt h)} f''(u)\,
	\mathbb E\left[
	\bigl(|\hat R_h|-|r_h|\bigr)^2
	\mathbf 1_{(|r_h|+\sqrt h,\,|r_h|+17\sqrt h)}
	(|\hat R_h|)
	\right]\\
	\leqslant& \frac{\alpha_{3,\Gamma}}{2}\sup_{u\in(0,18\sqrt h)} f''(u)\,
	|\hat r_h|\sqrt h\,
	\mathbf 1_{(0,\sqrt h]}(|r_h|)
	\mathbf 1_{[0,2R]}(|\hat r_h|).
\end{aligned}
\end{equation}
Since we are working under \(|r_h|\in(0,\sqrt h]\) and have already shown \(|\hat r_h|\leqslant 2R\), it remains to estimate \(f''\). Since \(18\sqrt h\leqslant r_1\), the definitions of
\(\varphi,\rho\), and \(f\) give, for \(u\in(0,18\sqrt h)\),
\[f''(u)=\varphi'(u)\rho(u)+\varphi(u)\rho'(u).\]
Since \(\varphi'(u)\rho(u)\leqslant 0\), we have 
\begin{equation}\label{Theq12}
	\sup_{u\in(0,18\sqrt h)}f''(u)\leqslant
	\sup_{u\in(0,18\sqrt h)}
	\varphi(u)\rho'(u)=
	-\frac{\alpha_{2, \Gamma}}{\alpha_{1, \Gamma}}
	\inf_{u\in(0,18\sqrt h)}\Phi(u+1)=
	-\frac{\alpha_{2, \Gamma}}{\alpha_{1, \Gamma}}\Phi(1),
\end{equation}
because \(\Phi\) is increasing. Substituting (\ref{Theq12}) into \eqref{Tqq1} and combining with \eqref{cc4}, (\ref{cc5}) gives
\begin{equation}\label{cc6}
\mathbb E\left[f(|\hat R_h|)-f(|r_h|)\right]\leqslant
\bigl||\hat r_h|-|r_h|\bigr|
-\frac{\alpha_{3,\Gamma}\alpha_{2, \Gamma}}{2\alpha_{1, \Gamma}}
|\hat r_h|\sqrt h\,\Phi(1).
\end{equation}
We now compare \(|r_h|\) and \(|\hat r_h|\). For any $0<h\leqslant (2M)^{-2}$, we see that $M\sqrt{h}\leqslant 1/2$. Then, by \eqref{rr}, we find that \(\left||r_h|-|\hat{r}_h|\right|\leqslant M\sqrt{h}|r_h|\leqslant |r_h|/2\).
Therefore,
\begin{equation}\label{tqq2}
|r_h|\leqslant 2|\hat{r}_h|\quad \text{and then,} \quad \left||r_h|-|\hat{r}_h|\right|\leqslant 2Mh^{1-\bar{\theta}}|\hat{r}_h|.
\end{equation}
This, together with (\ref{cc6}), yields
\begin{equation}\label{tqq3}
\E \left[f(|\hat{R}_h|)-f(|r_h|)\right]\leqslant
|\hat r_h|\sqrt h
\left(2Mh^{1/2-\bar\theta}
-\frac{\alpha_{3,\Gamma}\alpha_{2, \Gamma}\Phi(1)}
{2\alpha_{1, \Gamma}}
\right).
\end{equation}
By the definition of \(\bar h\), \(2Mh^{1/2-\bar\theta}
\leqslant 
\alpha_{3,\Gamma}\alpha_{2, \Gamma}\Phi(1)/(4\alpha_{1, \Gamma})\).
This implies that
\begin{equation*}
2Mh^{1/2-\bar\theta}
-\frac{\alpha_{3,\Gamma}\alpha_{2, \Gamma}\Phi(1)}
{2\alpha_{1, \Gamma}}
\leqslant-\frac{\alpha_{3,\Gamma}\alpha_{2, \Gamma}\Phi(1)}
{4\alpha_{1, \Gamma}}\leqslant -2Mh^{1/2-\bar{\theta}}\leqslant -2M\sqrt{h}.
\end{equation*}
Inserting this into \eqref{tqq3} gives 
\begin{equation}\label{P4.9}
\E \left[f(|\hat{R}_h|)-f(|r_h|)\right]\leqslant -2Mh|\hat{r}_h|.
\end{equation}
Moreover, as shown in (\ref{f1}), we have $f(u)\leqslant u$ for all $u\geqslant 0$. Then, using (\ref{tqq2}), we get $f(|r_h|)\leqslant |r_h|\leqslant 2|\hat{r}_h|$.
Therefore, combining this with \eqref{P4.9}, we arrive at the desired conclusion. Note that $1-Mh\in(0, 1)$ for any $h\in(0, h_1]$. The proof is complete.
\end{proof}
			
We establish the contractivity of the coupling for \underline{\bf{Case 2. \textit{$|r_h|\in(\sqrt{h}, r_1]$}}}. Define 
\begin{equation}\label{hatc}
\alpha_{5, \Gamma}=\mathrm{min}\left\{\frac{\varphi(r_1)R K}{4r_1},\frac{\alpha_{2, \Gamma}\alpha_{4, \Gamma}\Phi(1)}{4\alpha_{1, \Gamma}}\right\},
\end{equation}
where \(R, K\) are given in Assumption \ref{A1}, and constants $\alpha_{1, \Gamma}, \alpha_{2, \Gamma}, \alpha_{4, \Gamma}$ are defined in (\ref{c0}), (\ref{Ac}), and (\ref{c2}), respectively.
\begin{proposition}\label{PR4}
Let Assumptions \ref{A} and \ref{A1} hold.
For any $x, y\in\R^d$ and \(h \in (0, \bar{h}]\), where $\bar{h}$ is defined by \eqref{bar}, consider the coupling $(\hat{X}_h(x,Z), \hat{Y}_h(x,y,Z))$ defined by (\ref{hY}), with $m=8$. Then, if $|r_h|\in(\sqrt{h}, r_1]$, we have
\begin{equation*}
\E f(|\hat{R}_h|)\leqslant (1-\left(\alpha_{2, \Gamma}\wedge\alpha_{5, \Gamma}\right)h)f(|r_h|),
\end{equation*}
where $|r_h|, |\hat{r}_h|, |\hat{R}_h|$ are defined in (\ref{eq3}), function $f$, constants $\alpha_{2, \Gamma}$ and $\alpha_{5, \Gamma}$ are given by (\ref{f}), (\ref{Ac}) and (\ref{hatc}) respectively.
\end{proposition}
			
The proof of Proposition \ref{PR4} is rather technical, so we further subdivide this case  into three distinct scenarios:  \underline{$(\mathrm{I})~|\hat{r}_h|\in[0, \sqrt{h}]$}, \underline{$(\mathrm{II})~|\hat{r}_h|\in(\sqrt{h}, 2R]$}, and \underline{$(\mathrm{III})~|\hat{r}_h|\in(2R, +\infty)$}, based on the value of $|\hat{r}_h|$. In each of these cases, the conclusion of Proposition \ref{PR4} is established by Propositions \ref{PR2}, \ref{PR3}, and inequality \eqref{I-3}, respectively. The following lemma is prepared for handling case \underline{$(\mathrm{I})~|\hat{r}_h|\in[0, \sqrt{h}]$}.
\begin{lemma}\label{Le10}
Let Assumption \ref{A} hold. For any
\(x,y\in\mathbb R^d\) and \(h\in(0,1\wedge 4R^2]\),	consider the coupling \((\hat X_h(x,Z),\hat Y_h(x,y,Z))\) defined by \eqref{hY}, with \(m=8\). Then, we have
	\[\mathbb E\left[
		\bigl(|\hat R_h|-|\hat r_h|\bigr)^2
		\mathbf 1_{(0,\,|\hat r_h|+\sqrt h)}
		(|\hat R_h|)
		\right]\geqslant
		\alpha_{4,\Gamma}
		|\hat r_h|\sqrt h\,
		\mathbf 1_{[0,\sqrt h]}(|\hat r_h|),\]
	where $|r_h|, |\hat{r}_h|, |\hat{R}_h|$ are defined in (\ref{eq3}).
\end{lemma}

\begin{proof}
	If \(|\hat r_h|=0\) or \(|\hat r_h|>\sqrt h\), then the right-hand side is zero, and the assertion is
	immediate. Hence it remains to consider the case \(0<|\hat r_h|\leqslant \sqrt h\).
	Since \(h\leqslant 4R^2\), this implies \(|\hat r_h|\leqslant2R\),
	so the reflection part of the coupling is active. 
	Since the relevant
	terms in the reflection contribution depend on the noise only through
	\(\lambda_{\hat r_h}(\Sigma Z)\), by \eqref{refR}, the reflection part of the expectation gives
	\[
	\begin{aligned}
		&\mathbb E\left[
		\bigl(|\hat R_h|-|\hat r_h|\bigr)^2
		\mathbf 1_{(0,\,|\hat r_h|+\sqrt h)}
		(|\hat R_h|)
		\right] \\
		\geqslant&
		\int_{\mathbb R}
		\left(
		\left|
		|\hat r_h|+2s
		\right|
		-
		|\hat r_h|
		\right)^2
		\mathbf 1_{(0,\,|\hat r_h|+\sqrt h)}
		\left(
		\left|
		|\hat r_h|+2s
		\right|
		\right)                                                \\
		&\qquad\times
		\Bigl(
		\mathbf 1_{[-m,m]}(s)
		p_{0,h|\hat r_h|^2/|\hat r_h|_\Gamma^2}(s)
		-
		\mathbf 1_{[-m,m]}(|\hat r_h|+s)
		p_{0,h|\hat r_h|^2/|\hat r_h|_\Gamma^2}(|\hat r_h|+s)
		\Bigr)^+
		\,\mathrm ds .
	\end{aligned}
	\]
	We now restrict the last integral to \(s\in(0, \sqrt{h}/2)\).
	Since \(h\leqslant 1\), \(m=8\), and \(|\hat r_h|\leqslant\sqrt h\), we have \(0<s<m\) and \(0<|\hat r_h|+s<3\sqrt h/2<m\). Thus both truncation indicators are equal to one on this interval. Moreover, \(|\hat r_h|+s>s>0\), and hence, by monotonicity of the centred Gaussian density on \([0,\infty)\), \(p_{0,h|\hat r_h|^2/|\hat r_h|_\Gamma^2}(s)-
	p_{0,h|\hat r_h|^2/|\hat r_h|_\Gamma^2}(|\hat r_h|+s)
	\geqslant 0\).
	Furthermore, for \(0<s<\sqrt h/2\), \(0<
	|\hat r_h|+2s
	<
	|\hat r_h|+\sqrt h\) and therefore
	\[
	\left|
	|\hat r_h|+2s
	\right|
	=
	|\hat r_h|+2s,
	\qquad
	\left(
	\left|
	|\hat r_h|+2s
	\right|
	-
	|\hat r_h|
	\right)^2
	=
	4s^2 .
	\]
	Consequently,
	\begin{equation}\label{Le10-matrix-lower1}
			\mathbb E\left[
			\bigl(|\hat R_h|-|\hat r_h|\bigr)^2
			\mathbf 1_{(0,\,|\hat r_h|+\sqrt h)}
			(|\hat R_h|)
			\right]\geqslant
			4
			\int_0^{\sqrt h/2}
			s^2
			\left[
			p_{0,h|\hat r_h|^2/|\hat r_h|_\Gamma^2}(s)-
			p_{0,h|\hat r_h|^2/|\hat r_h|_\Gamma^2}(|\hat r_h|+s)
			\right]
			\,\mathrm ds .
	\end{equation}
Using the change of variables \(s=\sqrt h\,u\), the right-hand side of \eqref{Le10-matrix-lower1} becomes
	\[
	\begin{aligned}
		&4h
		\int_0^{1/2}
		u^2
		p_{0,|\hat r_h|^2/|\hat r_h|_\Gamma^2}(u)
		\left[
		1-
		\exp\left\{
		-\frac{|\hat r_h|_\Gamma^2(|\hat r_h|/\sqrt h+2u)}{2\sqrt h|\hat r_h|}
		\right\}
		\right]
		\,\mathrm du .
	\end{aligned}
	\]
	For \(u\in[0,1/2]\) and \(|\hat r_h|/\sqrt h\in(0,1]\), the exponent \(-|\hat r_h|_\Gamma^2(|\hat r_h|/\sqrt h+2u)/(2\sqrt h|\hat r_h|)\) belongs to \([-|\hat r_h|_\Gamma^2/|\hat r_h|^2, 0]\). By the concavity of \(x\mapsto1-\mathrm e^x\), we have, for every
	\(x\in[-|\hat r_h|_\Gamma^2/|\hat r_h|^2,0]\),
	\[
	1-\mathrm e^x
	\geqslant
	-\frac{|\hat r_h|^2}{|\hat r_h|_\Gamma^2}
	\left(1-\mathrm e^{-|\hat r_h|_\Gamma^2/|\hat r_h|^2}\right)x .
	\]
	Applying this, we get
	\[
	\begin{aligned}
		1-
		\exp\left\{
		-\frac{|\hat r_h|_\Gamma^2(|\hat r_h|/\sqrt h+2u)}{2\sqrt h|\hat r_h|}\right\}
		&\geqslant\frac{|\hat r_h|^2}{|\hat r_h|_\Gamma^2}
		\left(1-\mathrm e^{-|\hat r_h|_\Gamma^2/|\hat r_h|^2}\right)
		\frac{|\hat r_h|_\Gamma^2(|\hat r_h|/\sqrt h+2u)}{2\sqrt h|\hat r_h|}    \\
		&\geqslant 
		\left(1-\mathrm e^{-|\hat r_h|_\Gamma^2/|\hat r_h|^2}\right)
		\frac{|\hat r_h|}{\sqrt{h}}u .
	\end{aligned}
	\]
	Therefore
	\[\mathbb E\left[
		\bigl(|\hat R_h|-|\hat r_h|\bigr)^2
		\mathbf 1_{(0,\,|\hat r_h|+\sqrt h)}
		(|\hat R_h|)
		\right]\geqslant
		4|\hat r_h|\sqrt h
		\left(1-\mathrm e^{-|\hat r_h|_\Gamma^2/|\hat r_h|^2}\right)
		\int_0^{1/2}
		u^3
		p_{0,|\hat r_h|^2/|\hat r_h|_\Gamma^2}(u)\,\mathrm du .\]
	Since \(|\hat r_h|/|\hat r_h|_\Gamma\in[\vartheta_-,\vartheta_+]\) and
	\[\int_0^{1/2}
	u^3p_{0,|\hat r_h|^2/|\hat r_h|_\Gamma^2}(u)\,\mathrm du
	=\frac{|\hat r_h|^3}{|\hat r_h|_\Gamma^3}
	\int_0^{|\hat r_h|_\Gamma^2/(2|\hat r_h|^2)}
	v^3p_{0,1}(v)\,\mathrm dv ,\]
	we obtain
	\[\begin{aligned}
		&\mathbb E\left[
		\bigl(|\hat R_h|-|\hat r_h|\bigr)^2
		\mathbf 1_{(0,\,|\hat r_h|+\sqrt h)}
		(|\hat R_h|)
		\right] \\
		\geqslant&
		4|\hat r_h|\sqrt h
		\left(1-\mathrm e^{-|\hat r_h|_\Gamma^2/|\hat r_h|^2}\right)
		\frac{|\hat r_h|^3}{|\hat r_h|_\Gamma^3}
		\int_0^{|\hat r_h|_\Gamma^2/(2|\hat r_h|^2)}
		v^3p_{0,1}(v)\,\mathrm dv.
	\end{aligned}\]
By the definition of \(\alpha_{4,\Gamma}\), we get the desired conclusion. The proof is complete.
\end{proof}
			
We now invoke Lemmas \ref{1mom} and \ref{Le10} to establish Proposition \ref{PR4} for \underline{\bf{Case 2. \textit{$|r_h|\in(\sqrt{h}, r_1]$}}}  and \underline{$(\mathrm{I})~|\hat{r}_h|\in[0, \sqrt{h}]$}. 
\begin{proposition}\label{PR2}
Let Assumptions \ref{A} and \ref{A1} hold. For any $x, y\in\R^d$ and \( h\in(0, h_0\wedge(4M)^{-2}\wedge (4R^2)]\), consider the coupling $(\hat{X}_h(x,Z), \hat{Y}_h(x,y,Z))$ defined by (\ref{hY}), with $m=8$. Then, if $|r_h|\in(\sqrt{h}, r_1]$ and $|\hat{r}_h|\in [0, \sqrt{h}]$, we have
\begin{equation*}
\E f(|\hat{R}_h|)\leqslant (1-\alpha_{5, \Gamma}h)f(|r_h|),
\end{equation*}
where $|r_h|, |\hat{r}_h|, |\hat{R}_h|$ are defined in (\ref{eq3}), function $f$ is defined by (\ref{f}), constant $\alpha_{5, \Gamma}$ is defined by \eqref{hatc}.
\end{proposition}
\begin{proof}
We decompose \(\E[f(|\hat{R}_h|)-f(|r_h|)] = \E[f(|\hat{R}_h|)-f(|\hat{r}_h|)] + f(|\hat{r}_h|)-f(|r_h|)\). Therefore, it suffices to prove either 
\begin{equation}\label{Sc1}
\E[f(|\hat{R}_h|)-f(|\hat{r}_h|)]\leqslant 0,\quad f(|\hat{r}_h|)-f(|r_h|)\leqslant -\alpha_{5, \Gamma} hf(|r_h|)
\end{equation}
or
\begin{equation}\label{Sc2}
\E[f(|\hat{R}_h|)-f(|\hat{r}_h|)]\leqslant-\alpha_{5, \Gamma} hf(|r_h|), \quad f(|\hat{r}_h|)-f(|r_h|)\leqslant 0.	
\end{equation}
				
{\bf Step I. The case $|r_h|\in(R,r_1]$.} We first show that \(\mathbb E\bigl[f(|\hat R_h|)-f(|\hat r_h|)\bigr]\leqslant 0\).
According to \eqref{hY} and \eqref{Ac}, we have 
\[0\leqslant |\hat{R}_h|\leqslant |\hat{r}_h|+2|\lambda_{\hat r_h}(\Sigma Z)|\leqslant 2R+2m=r_1.\] 
By Lemma \ref{1mom}, \(\mathbb E|\hat R_h|=|\hat r_h|\).
Since $f$ is concave on $(0,r_1)$, Jensen's inequality yields
\[\mathbb E f(|\hat R_h|)
\leqslant f(\mathbb E|\hat R_h|)=f(|\hat r_h|).\]
Consequently, \(\mathbb E\bigl[f(|\hat R_h|)-f(|\hat r_h|)\bigr]\leqslant 0\). On the other hand, together with the fact that $f''(u)< 0$ on $u\in(0, r_1)$ and \eqref{Sub4eq1}, we have
\begin{align*}
f(|\hat{r}_h|)-f(|r_h|)\leqslant f'(|r_h|)(|\hat{r}_h|-|r_h|)\leqslant -\frac{hK}{2}f'(|r_h|)|r_h|.
\end{align*}
Since $f'(u)= \varphi(u)\rho(u)\geqslant \varphi(u)/2$ for all $u\in(0, r_1)$, $\varphi$ is decreasing on $[0, +\infty)$ and $f(u)\leqslant u$ for all $u\geqslant 0$, we obtain
$$\frac{hK}{2}f'(|r_h|)|r_h|\geqslant \frac{hK}{4}\varphi(|r_h|)|r_h|\geqslant \frac{hK}{4}\varphi(r_1)R\geqslant \frac{\varphi(r_1)R K}{4r_1}hf(r_1)\geqslant \frac{\varphi(r_1)R K}{4r_1}hf(|r_h|).$$ 
This, together with the definition of $\alpha_{5, \Gamma}$ in \eqref{hatc}, yields
\begin{equation*}
f(|\hat{r}_h|)-f(|r_h|) \leqslant -\frac{\varphi(r_1)R K}{4r_1}hf(|r_h|)\leqslant -\alpha_{5, \Gamma}hf(|r_h|).
\end{equation*}
Hence \eqref{Sc1} holds for $|r_h|\in(R,r_1]$.
				
{\bf Step II. The case $|r_h|\in(\sqrt h,R]$.} Since $|\hat r_h|\leqslant \sqrt h<|r_h|$ and $f$ is increasing, we have \(f(|\hat r_h|)-f(|r_h|)\leqslant 0\). For $|\hat{r}_h|\in [0, \sqrt{h}]$, by using Lemmas \ref{1mom} and \ref{Le10}, we have 
\begin{align*}
\E f(|\hat{R}_h|)-f(|\hat{r}_h|)
&\leqslant \frac{1}{2}\E\left[(|\hat{R}_h|-|\hat{r}_h|)^2\1_{(0, |\hat{r}_h|+\sqrt{h})}(|\hat{R}_h|)\right]\sup_{u\in(0, |\hat{r}_h|+\sqrt{h})}f''(u)\\
&\leqslant \frac{\alpha_{4, \Gamma}}{2}|\hat{r}_h|\sqrt{h}\sup_{u\in(0,2\sqrt{h})}f''(u).
\end{align*}
Similarly to (\ref{Theq12}), we get \(\sup_{u\in(0, 2\sqrt{h})}f''(u)\leqslant -\frac{\alpha_{2, \Gamma}}{\alpha_{1, \Gamma}}\Phi(1)\),
which yields
\begin{equation}\label{ccc}
\E f(|\hat{R}_h|)-f(|\hat{r}_h|)\leqslant -\frac{\alpha_{2, \Gamma}\alpha_{4, \Gamma}}{2\alpha_{1, \Gamma}}\Phi(1)|\hat{r}_h|\sqrt{h}.
\end{equation}
This, together with $f(|r_h|)\leqslant |r_h|$ and \eqref{tqq2} implies
$-|\hat{r}_h|\sqrt{h}\leqslant -|r_h| h/2\leqslant -h f(|r_h|)/2$. Substituting this into \eqref{ccc} and using \eqref{hatc} gives \(\E[f(|\hat{R}_h|)-f(|\hat{r}_h|)]\leqslant-\alpha_{5, \Gamma}hf(|r_h|)\). Hence \eqref{Sc2} holds for $|r_h|\in(\sqrt{h}, R]$. 

Therefore, combining the above results leads to the desired conclusion. 
Note that $\bar{\theta}\in(0, 1/2)$ and $M\geqslant 1\vee [\sqrt{\varphi(r_1)RK}/(8\sqrt{r_1})]$. Thus, by \eqref{hatc}, \[(4M)^{-2}\leqslant \frac{4r_1}{\varphi(r_1)RK}\leqslant \alpha_{5, \Gamma}^{-1}.\]  
Therefore, we have $1-\alpha_{5, \Gamma}h>0$ for any $h\in(0, (4M)^{-2}]$. The proof is complete.
\end{proof}
			
For \underline{\bf{Case 2. \textit{$|r_h|\in(\sqrt{h}, r_1]$}}} and \underline{$(\mathrm{III})~|\hat{r}_h|>2R$}, the synchronous coupling occurs. Hence, $|\hat{R}_h| = |\hat{r}_h|$. This implies $\E f(|\hat{R}_h|)-f(|\hat{r}_h|) = 0$. Hence, similarly to {\bf Step I} in Proposition \ref{PR2}, using \eqref{Sub4eq1} and \eqref{f1}, we have
\begin{equation}\label{I-3}
\E f(|\hat{R}_h|)-f(|r_h|) =f(|\hat{r}_h|)-f(|r_h|) \leqslant -\frac{\varphi(r_1)R K}{4r_1}hf(|r_h|).
\end{equation}
Next we consider the remaining \underline{\bf{Case 2. \textit{ $|r_h|\in(\sqrt{h}, r_1]$}}} and \underline{$(\mathrm{II})~\sqrt{h}<|\hat{r}_h|\leqslant 2R$}. To treat this case, we prove a useful lemma.

\begin{lemma}\label{2mom}
Let Assumption \ref{A} hold. For any $x, y\in\R^d$ and $h\in (0, (4M)^{-2}\wedge (4R^2)]$, consider the coupling $(\hat{X}_h(x,Z), \hat{Y}_h(x, y,Z))$ defined by (\ref{hY}), with \(m=8\). Then, we have
$$\E\left[(|\hat{R}_h|-|r_h|)^2\1_{\left(|r_h|-\sqrt{h},|r_h|\right)}(|\hat{R}_h|)\right]\geqslant \alpha_{1, \Gamma}h\1_{\left[\sqrt{h}, 2R\right]}(|\hat{r}_h|)\1_{\left[\sqrt{h}, h^{\bar{\theta}-1/2}/(4M)\right]}(|r_h|),$$
where $|r_h|, |\hat{r}_h|, |\hat{R}_h|$ are given by (\ref{eq3}), constant \(\alpha_{1,\Gamma}\) is defined by \eqref{c0}.
\end{lemma}
\begin{proof}
If either \(|\hat r_h|\notin[\sqrt h,2R]\) or \(|r_h|\notin[\sqrt h,h^{\bar\theta-1/2}/(4M)]\),
then the right-hand side is zero and the assertion is immediate. Hence it remains to consider the case \(|\hat r_h|\in[\sqrt h,2R]\) and \(|r_h|\in[\sqrt h,h^{\bar\theta-1/2}/(4M)]\).
By \(\1_{\left(|r_h|-\sqrt{h}, |r_h|\right)}(|\hat{R}_h|)\geqslant \1_{\left(|r_h|-\sqrt{h}, |r_h|-\sqrt{h}/4\right)}(|\hat{R}_h|)\), we observe that
\begin{equation}\label{ST1} 
\begin{aligned}
\E \left[(|\hat{R}_h|-|r_h|)^2\1_{\left(|r_h|-\sqrt{h}, |r_h|\right)}(|\hat{R}_h|)\right]\geqslant& \E \left[(|\hat{R}_h|-|r_h|)^2\1_{\left(|r_h|-\sqrt{h}, |r_h|-\sqrt{h}/4\right)}(|\hat{R}_h|)\right]\\
\geqslant &\frac{h}{16}\E \left[\1_{\left(|r_h|-\sqrt{h}, |r_h|-\sqrt{h}/4\right)}(|\hat{R}_h|)\right].
\end{aligned}
\end{equation}			 
We next establish a lower bound for \(\1_{\left(|r_h|-\sqrt{h}, |r_h|-\sqrt{h}/4\right)}(|\hat{R}_h|)\). Since $h\in(0, (4M)^{-2}]$, we have $\sqrt{h}\leqslant h^{\bar{\theta}-1/2}/(4M)$. Therefore, under the condition $|r_h|\in[\sqrt{h}, h^{\bar{\theta}-1/2}/(4M)]$, we obtain $Mh^{1-\bar{\theta}}|r_h|\leqslant \sqrt{h}/4$. This, together with \eqref{rr}, yields \(\left||r_h|-|\hat{r}_h|\right|\leqslant Mh^{1-\bar{\theta}}|r_h|\leqslant \sqrt{h}/4\),
which implies \(|r_h|-\sqrt{h}\leqslant|\hat{r}_h|-3\sqrt{h}/4\) and \(|\hat{r}_h|-\sqrt{h}/2\leqslant |r_h|-\sqrt{h}/4\).
Hence, 
\begin{equation}\label{ST2}
\1_{\left(|r_h|-\sqrt{h}, |r_h|-\sqrt{h}/4\right)}(|\hat{R}_h|)\geqslant \1_{\left(|\hat{r}_h|-3\sqrt{h}/4, |\hat{r}_h|-\sqrt{h}/2\right)}(|\hat{R}_h|)\1_{\left[\sqrt{h}, h^{\bar{\theta}-1/2}/(4M)\right]}(|r_h|).
\end{equation}
Then, inserting \eqref{ST2} into \eqref{ST1} gives 
\begin{equation*}
\begin{aligned}
&\E\left[(|\hat{R}_h|-|r_h|)^2\1_{\left(|r_h|-\sqrt{h}, |r_h|\right)}(|\hat{R}_h|)\right]\\
\geqslant&\frac{h}{16}\E \left[\1_{\left(|\hat{r}_h|-3\sqrt{h}/4, |\hat{r}_h|-\sqrt{h}/2\right)}(|\hat{R}_h|)\right]\1_{\left[\sqrt{h}, h^{\bar{\theta}-1/2}/(4M)\right]}(|r_h|).
\end{aligned}
\end{equation*}
It remains to derive a uniform lower bound for \(\E[\1_{\left(|\hat{r}_h|-3\sqrt{h}/4, |\hat{r}_h|-\sqrt{h}/2\right)}(|\hat{R}_h|)]\).
As in the proof of Lemma \ref{Le10}, by \eqref{refR}, the reflection contribution gives
\[
\begin{aligned}
	&\E \left[\1_{\left(|\hat{r}_h|-3\sqrt{h}/4, |\hat{r}_h|-\sqrt{h}/2\right)}(|\hat{R}_h|)\right]\\
	\geqslant&
	\int_{\mathbb R}
	\mathbf 1_{\left(|\hat r_h|-3\sqrt h/4,\,|\hat r_h|-\sqrt h/2\right)}
	\left(\left||\hat r_h|+2s\right|\right)  \\
	&\qquad\times
	\left(
	\mathbf 1_{[-m,m]}(s)
	p_{0,h|\hat r_h|^2/|\hat r_h|_\Gamma^2}(s)
	-
	\mathbf 1_{[-m,m]}(|\hat r_h|+s)
	p_{0,h|\hat r_h|^2/|\hat r_h|_\Gamma^2}(|\hat r_h|+s)
	\right)^+
	\,\mathrm ds .
\end{aligned}
\]
We now restrict the integral to \(s\in(-3\sqrt h/8,-\sqrt h/4)\).
Since \(M\geqslant 1\) and \(h\leqslant (4M)^{-2}\), we have \(|s|<3\sqrt h/8<m\). Moreover, for \(|\hat r_h|\in[\sqrt h,2R]\),
we have \(\left||\hat r_h|+2s\right|<m\) and 
\[
\left||\hat r_h|+2s\right|
=
|\hat r_h|+2s
\in
\left(|\hat r_h|-\frac{3\sqrt h}{4},
|\hat r_h|-\frac{\sqrt h}{2}\right).
\]
Furthermore, \(|\hat r_h|+2s>0\) implies \(|s|<|\hat r_h|+s\).
By monotonicity of the centred Gaussian density on \([0,\infty)\), \(p_{0,h|\hat r_h|^2/|\hat r_h|_\Gamma^2}(s)
\geq
p_{0,h|\hat r_h|^2/|\hat r_h|_\Gamma^2}(|\hat r_h|+s)\).
Consequently,
\begin{equation}\label{J8-lower1}
	\begin{aligned}
		&\E \left[\1_{\left(|\hat{r}_h|-3\sqrt{h}/4, |\hat{r}_h|-\sqrt{h}/2\right)}(|\hat{R}_h|)\right]\\
		\geqslant&
		\mathbf 1_{[\sqrt h,2R]}(|\hat r_h|)
		\int_{-3\sqrt h/8}^{-\sqrt h/4}
		\left[
		p_{0,h|\hat r_h|^2/|\hat r_h|_\Gamma^2}(s)
		-
		p_{0,h|\hat r_h|^2/|\hat r_h|_\Gamma^2}(|\hat r_h|+s)
		\right]
		\,\mathrm ds .
	\end{aligned}
\end{equation}
Making the change of variables \(s=\sqrt h\,u\), the right-hand side of \eqref{J8-lower1} becomes
\begin{equation}\label{J8-lower2}
	\mathbf 1_{[\sqrt h,2R]}(|\hat r_h|)
	\int_{-3/8}^{-1/4}
	p_{0,|\hat r_h|^2/|\hat r_h|_\Gamma^2}(u)
	\left[
	1-
	\exp\left\{
	-\frac{|\hat r_h|_\Gamma^2}{|\hat r_h|^2}
	\frac{|\hat r_h|}{\sqrt h}
	\left(
	\frac{|\hat r_h|}{2\sqrt h}+u
	\right)
	\right\}
	\right]
	\,\mathrm du .
\end{equation}
On the present event, \(|\hat r_h|\geqslant \sqrt h\). For fixed \(u\in (-3/8,-1/4)\) define
\[
H_u(t)
:=
1-
\exp\left\{
-\frac{|\hat r_h|_\Gamma^2}{|\hat r_h|^2}
t\left(\frac{t}{2}+u\right)
\right\}.
\]
Observe that \(H_u(t)\geqslant H_u(1)\) for all \(t\geqslant 1\).
Moreover, for \(u\in(-3/8,-1/4)\), \(1/2+u\in (1/8, 1/4)\).
Thus
\[
H_u(1)
=
1-
\exp\left\{
-\frac{|\hat r_h|_\Gamma^2}{|\hat r_h|^2}
\left(\frac12+u\right)
\right\}
\geqslant
1-
\exp\left\{
-\frac{|\hat r_h|_\Gamma^2}{8|\hat r_h|^2}
\right\}.
\]
Consequently, since \(p_{0,a}(\cdot)\) is even, \eqref{J8-lower2} implies
\[\begin{aligned}
	&\E \left[\1_{\left(|\hat{r}_h|-3\sqrt{h}/4, |\hat{r}_h|-\sqrt{h}/2\right)}(|\hat{R}_h|)\right]\\
	\geqslant&
	\mathbf 1_{[\sqrt h,2R]}(|\hat r_h|)
	\left(
	1-
	\exp\left\{
	-\frac{|\hat r_h|_\Gamma^2}{8|\hat r_h|^2}
	\right\}
	\right)
	\int^{3/8}_{1/4}
	p_{0,|\hat r_h|^2/|\hat r_h|_\Gamma^2}(u)
	\,\mathrm du .
\end{aligned}\]	
By the uniform ellipticity bounds, \(|\hat r_h|/|\hat r_h|_\Gamma\in[\vartheta_-,\vartheta_+]\), and by the definition of \(\alpha_{1,\Gamma}\) in \eqref{c0}, we obtain the desired conclusion. The proof is complete.
\end{proof}
			
We now invoke Lemmas \ref{1mom} and \ref{2mom} to establish Proposition \ref{PR4} for \underline{\bf{Case 2. \textit{ $|r_h|\in(\sqrt{h}, r_1]$}}} and \underline{$(\mathrm{II})~|\hat{r}_h|\in(\sqrt{h}, 2R]$}.
\begin{proposition}\label{PR3}
Let Assumptions \ref{A} and \ref{A1} hold. For any $x, y\in\R^d$ and $h \in (0, \bar{h}]$, where $\bar{h}$ is defined by \eqref{bar}, consider the coupling $(\hat{X}_h(x,Z), \hat{Y}_h(x,y,Z))$ defined by (\ref{hY}), with $m=8$. Then, if $|r_h|\in(\sqrt{h}, r_1]$ and $|\hat{r}_h|\in (\sqrt{h}, 2R]$, we have
\begin{equation*}
\E f(|\hat{R}_h|)\leqslant (1-\alpha_{2, \Gamma}h)f(|r_h|),
\end{equation*}
where $|r_h|, |\hat{r}_h|, |\hat{R}_h|$ are defined in (\ref{eq3}), the function
$f$ is defined in \eqref{f}, and the constant $\alpha_{2,\Gamma}$ is
defined in \eqref{Ac}.
\end{proposition}
\begin{proof}
By the definition of $\bar h$, we have $h \in (0, (4Mr_1)^{2/(2\bar{\theta}-1)}]$. Since $\bar\theta\in(0,1/2)$, this implies $r_1\leqslant h^{\bar{\theta}-1/2}/(4M)$. Hence, for every $|r_h|\in(\sqrt h,r_1]$, \(|r_h|\in[\sqrt h,h^{\bar\theta-1/2}/(4M)]\). Together with $|\hat{r}_h|\in(\sqrt h,2R]$, Lemma \ref{2mom} gives
\begin{equation}\label{Theq1}
\E\left[(|\hat{R}_h|-|r_h|)^2\1_{(|r_h|-\sqrt{h}, |r_h|)}(|\hat{R}_h|)\right]\geqslant \alpha_{1, \Gamma}h.
\end{equation}
Since $f''< 0$ on $(0, r_1)$, by \eqref{cc4} and Lemma \ref{1mom}, we have
\begin{align*}
&\E f(|\hat{R}_h|)-f(|r_h|)\\
\leqslant &f'(|r_h|)(|\hat{r}_h|-|r_h|)+\frac{1}{2}\E\left[(|\hat{R}_h|-|r_h|)^2\1_{(|r_h|-\sqrt{h}, |r_h|)}(|\hat{R}_h|)\right]\sup_{u\in(|r_h|-\sqrt{h}, |r_h|)}f''(u).
\end{align*}
We therefore obtain from (\ref{Theq1}) that
\begin{equation*}
\E f(|\hat{R}_h|)-f(|r_h|)\leqslant f'(|r_h|)(|\hat{r}_h|-|r_h|) + \frac{\alpha_{1, \Gamma} h}{2}\sup_{u\in(|r_h|-\sqrt{h}, |r_h|)}f''(u).
\end{equation*}
Since $f'(u) = \varphi(u)\rho(u)$ on $u\in (0, r_1]$, we have
$$\sup_{u\in(|r_h|-\sqrt{h}, |r_h|)}f''(u)\leqslant \sup_{u\in(|r_h|-\sqrt{h}, |r_h|)}(\varphi'(u)\rho(u))+\sup_{u\in(|r_h|-\sqrt{h}, |r_h|)}(\varphi(u)\rho'(u)).$$
Therefore, 
\begin{align*}
\E f(|\hat{R}_h|)-f(|r_h|)\leqslant& f'(|r_h|)(|\hat{r}_h|-|r_h|) + \frac{\alpha_{1, \Gamma} h}{2} \sup_{u\in(|r_h|-\sqrt{h}, |r_h|)}(\varphi'(u)\rho(u))\\
&+\frac{\alpha_{1, \Gamma} h}{2}\sup_{u\in(|r_h|-\sqrt{h}, |r_h|)}(\varphi(u)\rho'(u)).
\end{align*}
It remains to prove
\begin{equation}\label{Ca1}
f'(|r_h|)(|\hat{r}_h|-|r_h|)+\frac{\alpha_{1, \Gamma} h}{2} \sup_{u\in(|r_h|-\sqrt{h}, |r_h|)}(\varphi'(u)\rho(u))\leqslant 0
\end{equation}
and 
\begin{equation}\label{Ca2}
\frac{\alpha_{1, \Gamma}h}{2}\sup_{u\in(|r_h|-\sqrt{h}, |r_h|)}(\varphi(u)\rho'(u))\leqslant -\alpha_{2, \Gamma}h f(|r_h|).
\end{equation}
				
{\bf Step I. Proof of \eqref{Ca1}.}
Since $\rho$ is decreasing and positive, while $\varphi'<0$, we have
\begin{align*}
\sup_{u\in(|r_h|-\sqrt{h}, |r_h|)}(\varphi'(u)\rho(u)) \leqslant  \rho(|r_h|)\sup_{u\in(|r_h|-\sqrt{h}, |r_h|)}\varphi'(u),
\end{align*}
Since $f'=\varphi\rho$, it is enough to show that
\begin{equation}\label{Theq7}
\sup_{u\in(|r_h|-\sqrt{h}, |r_h|)}\varphi'(u)\leqslant \frac{2(|r_h|-|\hat{r}_h|)}{\alpha_{1, \Gamma} h}\varphi(|r_h|).
\end{equation}
If $|r_h|\in(R, r_1]$, then by \eqref{Sub4eq1}, $|r_h|-|\hat{r}_h|\geqslant (hK/2)|r_h|>0$. Since $\varphi'\leqslant 0$, \eqref{Theq7} follows immediately. It remains to consider the case $|r_h|\in(\sqrt h,R]$. By \eqref{varphi}, \(\varphi'(u)
=
-(2L+1)(u+1)\varphi(u)/\alpha_{1,\Gamma}\). Hence
\begin{equation}\label{Theq3}
\sup_{u\in(|r_h|-\sqrt{h}, |r_h|)}\varphi'(u)
=-\frac{2L+1}{\alpha_{1, \Gamma}}\inf_{u\in(|r_h|-\sqrt{h}, |r_h|)}\left((u+1) \varphi(u)\right).
\end{equation}
Since $\varphi$ is decreasing and $h\leqslant 1$, we have
\begin{equation}\label{Tq2}
\inf_{u\in(|r_h|-\sqrt{h}, |r_h|)}\left((u+1) \varphi(u)\right)\geqslant(|r_h|-\sqrt{h}+1)\varphi(|r_h|)\geqslant |r_h|\varphi(|r_h|).	   
\end{equation}
Substituting (\ref{Tq2}) into (\ref{Theq3}) yields 
\begin{equation}\label{Tq3}
\sup_{u\in(|r_h|-\sqrt{h}, |r_h|)}\varphi'(u)\leqslant -\frac{(2L+1)|r_h|}{\alpha_{1, \Gamma}}\varphi(|r_h|).
\end{equation}
We now estimate $|r_h|-|\hat r_h|$.	
By (\ref{ss}), and using the fact that \(h\leqslant M^{1/(2\bar\theta-1)}\) so that \(Mh^{1-2\bar{\theta}}\leqslant 1\), we get 
\begin{equation*}	
\frac{|r_h|-|\hat{r}_h|}{h}\geqslant \frac{|r_h|}{h}\left(1-\sqrt{1+2hL+Mh^{2-2\bar{\theta}}}\right)\geqslant \frac{1-\sqrt{1+(2L+1)h}}{h}|r_h|.
\end{equation*}		     
One observes that function $x\mapsto (\sqrt{1+(2L+1) x}-1)/x$ is decreasing on $(0, +\infty)$. This implies that for any $x>0$,
\begin{align*}
0\leqslant \frac{\sqrt{1+(2L+1)x}-1}{x}\leqslant \frac{2L+1}{2}.
\end{align*}
Consequently, \((1-\sqrt{1+(2L+1)h})/h \geqslant -(2L+1)/2\). Thus, \(-(2L+1)|r_h|\leqslant 2(|r_h|-|\hat r_h|)/h\).
Combining this with \eqref{Tq3} proves \eqref{Theq7}. Hence
\eqref{Ca1} follows.
				
{\bf Step II. Proof of \eqref{Ca2}.} By the definition of \(f\) given in \eqref{f}, we have $\rho'(u) = -2\alpha_{2, \Gamma}\Phi(u+1)/(\alpha_{1, \Gamma}\varphi(u))$. Therefore, \(\varphi(u)\rho'(u)=-2\alpha_{2,\Gamma}
\Phi(u+1)/\alpha_{1,\Gamma}\). Since \(\Phi\) is increasing,
\begin{align*}
\frac{\alpha_{1, \Gamma}}{2}\sup_{u\in(|r_h|-\sqrt{h},|r_h|)}(\varphi(u)\rho'(u))
= -\alpha_{2, \Gamma}\Phi(|r_h|+1-\sqrt{h})\leqslant-\alpha_{2, \Gamma}\Phi(|r_h|),
\end{align*}
where we use the fact that $\Phi(\cdot)$ is increasing. 
Using $f(|r_h|)\leqslant \Phi(|r_h|)$ as shown in (\ref{f1}), we obtain
$$\frac{\alpha_{1, \Gamma}}{2}\sup_{u\in(|r_h|-\sqrt{h}, |r_h|)}(\varphi(u)\rho'(u))\leqslant -\alpha_{2, \Gamma}f(|r_h|).$$
Multiplying by \(h\) gives \eqref{Ca2}.

Combining \eqref{Ca1} and \eqref{Ca2}, we conclude that \(\mathbb E f(|\hat R_h|)
\leqslant (1-\alpha_{2,\Gamma}h)f(|r_h|)\).				
Finally, by \eqref{dd6}, $\alpha_{2, \Gamma}^{-1}\geqslant 2r_1\Phi(1)/\alpha_{1, \Gamma}$. Moreover, since $M\geqslant \sqrt{\alpha_{1, \Gamma}}/\sqrt{32r_1\Phi(1)}$, we have \[(4M)^{-2}\leqslant \frac{2r_1\Phi(1)}{\alpha_{1, \Gamma}} \leqslant \alpha_{2, \Gamma}^{-1}.\] 
Thus, for \(h\leqslant (4M)^{-2}\), $1-\alpha_{2, \Gamma}h>0$. The proof is complete.	
\end{proof}
			
For \underline{\bf{Case 3. \textit{$|r_h|\in(r_1, +\infty)$}}}, we note that $|\hat{r}_h|\leqslant (1-hK/2)|r_h|<|r_h|$. Thus, by \eqref{cc4} and Lemma \ref{1mom}, we have $\E f(|\hat{R}_h|)-f(|\hat{r}_h|)\leqslant 0$. Hence we only consider $f(|\hat{r}_h|)-f(|r_h|)$. Since we do not know whether $|r_h|(1-hK/2)$ is larger or smaller than $r_1$,  we use a case-by-case discussion for establishing Theorem \ref{Th1}.
\begin{proposition}\label{PR5}
Let Assumptions \ref{A} and \ref{A2} hold. For any $x, y\in\R^d$ and $h\in(0, 2/K]$, consider the coupling $(\hat{X}_h(x,Z), \hat{Y}_h(x,y,Z))$ defined by (\ref{hY}). If $|r_h|\in(r_1, +\infty)$, then
\begin{equation*}
\E f(|\hat{R}_h|)\leqslant (1-K\varphi(r_1)h/8)f(|r_h|),
\end{equation*}
where $|r_h|, |\hat{r}_h|, |\hat{R}_h|$ are defined in (\ref{eq3}), functions $f$ and $\varphi$ are defined by (\ref{f}) and (\ref{varphi}) respectively.
\end{proposition} 
\begin{proof}
If $|r_h|(1-hK/2)>r_1$, since $f$ is increasing, using the fact $\rho(u)\geqslant 1/2$ for $u\in(0,r_1)$, we have
\begin{align*}
f(|\hat{r}_h|)-f(|r_h|)&\leqslant f\left(|r_h|\left(1-hK/2\right)\right)-f(|r_h|)\\
&\leqslant -\frac{1}{2}\int_{|r_h|(1-hK/2)}^{|r_h|}\varphi(s\wedge r_1)\md s= -\frac{hK}{4}\varphi(r_1)|r_h|.
\end{align*}
Therefore, by the fact $f(|r_h|)\leqslant |r_h|$ given in (\ref{f1}), we obtain the desired result. It remains to consider the case \(|r_h|(1-hK/2)\leqslant r_1\).
Set the intermediate point \(|r_h|(1-hK/4)\).
We distinguish two subcases. If $|r_h|(1-hK/4)< r_1$, then using the fact that $\varphi$ is decreasing and $\rho(r)\geqslant 1/2$ for $r\in(0,r_1)$, we have 
\begin{align*}
f(|\hat{r}_h|)-f(|r_h|)&\leqslant f\left(|r_h|\left(1-\frac{hK}{2}\right)\right)-f\left(|r_h|\left(1-\frac{hK}{4}\right)\right)\\
&\leqslant -\frac{1}{2}\int_{|r_h|(1-\frac{hK}{2})}^{|r_h|(1-hK/4)}\varphi(s\wedge r_1)\md s\leqslant -\frac{K\varphi(r_1)}{8}h|r_h|\leqslant -\frac{K\varphi(r_1)}{8}hf(|r_h|).
\end{align*}
If $|r_h|(1-hk/4)\geqslant r_1$, then 
\begin{align*}
f(|\hat{r}_h|)-f(|r_h|)&\leqslant f\left(|r_h|\left(1-\frac{hK}{4}\right)\right)-f(|r_h|)\leqslant -\frac{1}{2}\int_{|r_h|(1-hK/4)}^{|r_h|}\varphi(s\wedge r_1)\md s\\
&= -\frac{K\varphi(r_1)}{8}h|r_h|\leqslant -\frac{K\varphi(r_1)}{8}hf(|r_h|).
\end{align*}	
Hence, combining all estimates, we obtain
\begin{equation*}
\E f(|\hat{R}_h|)-f(|r_h|)\leqslant f(|\hat{r}_h|)-f(|r_h|)\leqslant  -\frac{K\varphi(r_1)}{8}hf(|r_h|)
\end{equation*}
for $|r_h|\in(r_1, +\infty)$.
Finally, since $0<\varphi(r_1)<1$ and $h\leqslant 2/K$, we have $K\varphi(r_1)h/ 8
\leqslant\varphi(r_1)/4<1$. Then, $1-K\varphi(r_1)h/8>0$. The proof is complete.
\end{proof}
			
			Combining Propositions \ref{PR1}, \ref{PR4} and \ref{PR5}, we establish the contractivity of coupling $(\hat{X}_h(x,Z),$ $ \hat{Y}_h(x, y,Z))$, as formalized in Theorem \ref{Th1}.
			
			\subsection{Uniform-in-time error analysis}
			In this section, we give the uniform-in-time strong error bound for the TEM scheme (\ref{TEM}) and derive the $1/2$-order convergence rate for the corresponding numerical invariant measure $\mu_h$ to the underlying invariant measure $\mu$. 
			
			Noting that, after one step of the coupling, the initial conditions for the
			next step become random, we now extend the construction from deterministic
			initial points to random initial variables. More precisely, let
			\(\xi_1\) and \(\xi_2\) be \(\mathbb R^d\)-valued random variables which are
			independent of \((Z,\zeta)\), where \(Z\sim N(\mathbf 0,h\mathbf I_q)\) and \(\zeta\sim {\rm Unif}[0,1]\).
			The random variables \(\xi_1\) and \(\xi_2\) need not be independent of each
			other. Define
			\[
			u_h(\xi_i):=\pi_h(\xi_i)+h b(\pi_h(\xi_i)),
			\qquad i=1,2,
			\]
			and \(\hat r_h(\xi_1,\xi_2):=u_h(\xi_1)-u_h(\xi_2)\).
			Define \(\hat X_h(\xi_1,Z):=u_h(\xi_1)+\Sigma Z\).
			If \(\hat r_h(\xi_1,\xi_2)=\mathbf 0\), set \(\hat Y_h(\xi_1,\xi_2,Z):=u_h(\xi_2)+\Sigma Z\).
			If \(\hat r_h(\xi_1,\xi_2)\neq\mathbf 0\), define
			\begin{equation}\label{XY0}
				\hat Y_h(\xi_1,\xi_2,Z)
				:=
				\begin{cases}
					u_h(\xi_1)+\Sigma Z,
					& \zeta\leqslant v_{\hat r_h(\xi_1,\xi_2)}^m(\Sigma Z),\
					\Sigma Z\in A_{\hat r_h(\xi_1,\xi_2)}^m,\\
					&~~~~~~~~~~~~~~~~~~~~~~~\text{and}~|\hat r_h(\xi_1,\xi_2)|\leqslant 2R,\\[1mm]
					u_h(\xi_2)+\mathcal R_{\hat r_h(\xi_1,\xi_2)}^\Gamma \Sigma Z,
					& \zeta> v_{\hat r_h(\xi_1,\xi_2)}^m(\Sigma Z),\
					\Sigma Z\in A_{\hat r_h(\xi_1,\xi_2)}^m,\\
					&~~~~~~~~~~~~~~~~~~~~~~~\text{and}~|\hat r_h(\xi_1,\xi_2)|\leqslant 2R,\\[1mm]
					u_h(\xi_2)+\Sigma Z,
					& \Sigma Z\notin A_{\hat r_h(\xi_1,\xi_2)}^m
					~\text{or}~
					|\hat r_h(\xi_1,\xi_2)|>2R .
				\end{cases}
			\end{equation}
			Here \(v_{\hat r_h(\xi_1,\xi_2)}^m\) and \(A_{\hat r_h(\xi_1,\xi_2)}^m\)
			are defined exactly as in the deterministic case, with \(x\) and \(y\)
			replaced by \(\xi_1\) and \(\xi_2\), respectively.
			
			\begin{corollary}\label{rvcoup}
				For fixed \(h>0\), let \(\xi_1\) and \(\xi_2\) be \(\mathbb R^d\)-valued
				random variables independent of \((Z,\zeta)\). Then the $2d$-dimensional random vector $(\hat{X}_h(\xi_1,Z), \hat{Y}_h(\xi_1, \xi_2,Z))$ defined by (\ref{XY0}) is a coupling of the laws of
				\(\hat X_h(\xi_1,Z)\) and \(\hat X_h(\xi_2,Z)\).	
			\end{corollary}
			\begin{proof}
				The first marginal is clear from the definition of \(\hat X_h(\xi_1,Z)\). It remains to identify the second marginal. Let \(g:\mathbb R^d\to\mathbb R\)
				be bounded and continuous, and denote \(\mathcal G:=\Sigma(\xi_1,\xi_2)\).
				Since \((\xi_1,\xi_2)\) is independent of \((Z,\zeta)\), conditioning on
				\(\mathcal G\) freezes the values of \(\xi_1\) and \(\xi_2\). Hence, by \cite[Theorem 2.24]{FC} and Lemma \ref{coup}, 
				\[
				\begin{aligned}
					\mathbb E\left[
					g\bigl(\hat Y_h(\xi_1,\xi_2,Z)\bigr)
					\right]
					&=
					\mathbb E\left[
					\mathbb E\left[
					g\bigl(\hat Y_h(\xi_1,\xi_2,Z)\bigr)
					\mid \mathcal G
					\right]
					\right]=
					\mathbb E\left[
					\left.
					\mathbb E\left[
					g\bigl(\hat Y_h(x,y,Z)\bigr)
					\right]
					\right|_{x=\xi_1,\ y=\xi_2}
					\right] \\
					&=
					\mathbb E\left[
					\left.
					\mathbb E\left[
					g\bigl(u_h(y)+\Sigma Z\bigr)
					\right]
					\right|_{y=\xi_2}
					\right]=
					\mathbb E\left[
					g\bigl(\hat X_h(\xi_2,Z)\bigr)
					\right].
				\end{aligned}
				\]
				Therefore, \(\mathcal L\bigl(\hat Y_h(\xi_1,\xi_2,Z)\bigr)
				=\mathcal L\bigl(\hat X_h(\xi_2,Z)\bigr)\).
				Thus \(\bigl(\hat X_h(\xi_1,Z),\hat Y_h(\xi_1,\xi_2,Z)\bigr)\)
				is indeed a coupling of the laws of
				\(\hat X_h(\xi_1,Z)\) and \(\hat X_h(\xi_2,Z)\). The proof is complete.
			\end{proof}
			Similarly, the conclusion of Theorem \ref{Th1} remains valid for the
			coupling defined by \eqref{XY0}. We now establish the uniform-in-time
			strong convergence rate of the TEM scheme \eqref{TEM} for \eqref{E}.
			
			\begin{theorem}\label{Th2}
				Let Assumptions \ref{A1} and \ref{A2} hold. For any
				\(h\in(0,\bar h]\), we have \[\sup_{k\geqslant 0}
				\mathbb E\left|X_k-x_{t_k}\right|
				\leqslant
				C h^{1/2},\]
				where \(\bar h\) is defined in \eqref{bar}, and \(C\) is a positive
				constant independent of \(h\) and \(k\).
			\end{theorem}
			
			\begin{proof}
				We construct a process \((\hat G_k)_{k\geqslant 0}\) coupled with
				\((\hat S_k)_{k\geqslant 1}\) according to the coupling mechanism
				\eqref{XY0}. By Corollary \ref{rvcoup}, if
				\(\hat G_0=\hat X_0\), then \(\mathcal L(\hat G_k)=\mathcal L(\hat X_k)\) for all \(k\geqslant 0\).
				Moreover, the processes \((\hat G_k)_{k\geq 0}\) and
				\((\hat S_k)_{k\geq 1}\) inherit the moment estimates stated in Lemma \ref{Le2}. By \eqref{f1}, Corollary \ref{co2} with \(s\geqslant 1\vee(3\ell/(2\theta))\),
				and Lemma \ref{Le5}, we obtain, for any \(k\geqslant 0\),
				\begin{equation}\label{35}
					\mathbb E f\left(
					\left|\pi_h(\hat G_{k+1})-x_{t_{k+1}}\right|
					\right)
					\leqslant
					C h^{3/2}
					+
					\mathbb E f\left(
					\left|\hat G_{k+1}-\hat S_{k+1}\right|
					\right).
				\end{equation}				
				Using Theorem \ref{Th1}, Corollary \ref{rvcoup}, and Corollary \ref{co2}, we have, for any \(k\geqslant 0\),
				\[\begin{aligned}
					\mathbb E f\left(
					\left|\hat G_{k+1}-\hat S_{k+1}\right|
					\right)
					&\leqslant
					(1-ch)
					\mathbb E f\left(
					\left|\pi_h(\hat G_k)-\pi_h(x_{t_k})\right|
					\right)  \\
					&\leqslant
					(1-ch)
					\mathbb E f\left(
					\left|\pi_h(\hat G_k)-x_{t_k}\right|
					\right)+(1-ch)
					\mathbb E f\left(
					\left|x_{t_k}-\pi_h(x_{t_k})\right|
					\right)  \\
					&\leqslant(1-ch)
					\mathbb E f\left(
					\left|\pi_h(\hat G_k)-x_{t_k}\right|
					\right)
					+C h^{3/2}.
				\end{aligned}\]
				Here we used the subadditivity of \(f\), which follows from the fact
				that \(f\) is increasing, concave, and \(f(0)=0\).
				Substituting the above estimate into \eqref{35} gives
				\[\mathbb E f\left(
					\left|\pi_h(\hat G_{k+1})-x_{t_{k+1}}\right|
					\right)\leqslant
					(1-ch)
					\mathbb E f\left(
					\left|\pi_h(\hat G_k)-x_{t_k}\right|
					\right)
					+C h^{3/2}.\]
				Iterating this inequality yields
				\[\mathbb E f\left(
					\left|\pi_h(\hat G_{k+1})-x_{t_{k+1}}\right|
					\right)
					\leqslant
					(1-ch)^{k+1}
					f\left(
					\left|\pi_h(\hat G_0)-\pi_h(x_0)\right|
					\right)+
					C h^{3/2}
					\sum_{j=0}^{k}(1-ch)^j .\]
				Since the two processes start from the same initial value, the first
				term vanishes. Therefore,
				\[\sup_{k\geqslant 0}
				\mathbb E f\left(
				\left|\pi_h(\hat G_k)-x_{t_k}\right|
				\right)
				\leqslant C h^{1/2}.\]				
				Finally, by \eqref{f1}, using \(X_k=\pi_h(\hat X_k)\) and
				\(\mathcal L(\hat X_k)=\mathcal L(\hat G_k)\), we obtain
			$$\E\left|X_{k}-x_{t_{k}}\right|\leqslant \frac{2}{\varphi(r_1)}\E f\left(\left|X_{k}-x_{t_{k}}\right|\right)= \frac{2}{\varphi(r_1)}\E f(|\pi_h(\hat{G}_{k})-x_{t_{k}}|)\leqslant Ch^{1/2}.$$
				Taking the supremum over \(k\geqslant 0\) gives the desired result.
			\end{proof}
			
			Based on the above result, we obtain the convergence rate of the
			numerical invariant measure to the invariant measure of the underlying
			SDE in the Wasserstein distance.
			
			\begin{theorem}\label{Th3}
				Let Assumptions \ref{A1} and \ref{A2} hold. Then, for any initial value
				\(x_0\in\mathbb R^d\) and any \(h\in(0,\bar h]\), we have
				\[
				\mathcal W_1(\mu_h,\mu)
				\leqslant
				C h^{1/2},
				\]
				where \(\bar h\) is defined in \eqref{bar}, and \(C\) is a positive
				constant independent of \(h\).
			\end{theorem}
			
			\begin{proof}
				By the triangle inequality, for every \(k\geqslant 0\),
				\[\mathcal W_1(\mu_h,\mu)
					\leqslant
					\mathcal W_1\bigl(\mu_h,\mathcal L(X_k)\bigr)
					+
					\mathcal W_1\bigl(\mathcal L(X_k),\mathcal L(x_{t_k})\bigr)+
					\mathcal W_1\bigl(\mathcal L(x_{t_k}),\mu\bigr).\]
				By Lemma \ref{LeX} and Theorem \ref{Co},
				\[\lim_{k\to\infty}
				\mathcal W_1\bigl(\mu_h,\mathcal L(X_k)\bigr)
				=0,\qquad \lim_{k\to\infty}
				\mathcal W_1\bigl(\mathcal L(x_{t_k}),\mu\bigr)
				=0.\]
				Therefore, taking \(\limsup_{k\to\infty}\) in the preceding inequality
				gives
				\[\mathcal W_1(\mu_h,\mu)
				\leqslant
				\limsup_{k\to\infty}
				\mathcal W_1\bigl(\mathcal L(X_k),\mathcal L(x_{t_k})\bigr).\]
				Using Theorem \ref{Th2}, we obtain \(\mathcal W_1\bigl(\mathcal L(X_k),\mathcal L(x_{t_k})\bigr)
				\leqslant\mathbb E\left|X_k-x_{t_k}\right|
				\leqslant C h^{1/2}\).
				Consequently, the proof is complete.
			\end{proof}
			
			\section{Numerical examples}
			This section presents an example with long-range dissipativity and reports numerical experiments based on the truncated Euler--Maruyama scheme.
			\begin{example}\label{EX}
				Consider the planar SDE
				\begin{equation}\label{Ex}
					\md X_t = b(X_t)\md t+\Sigma\md B_t,\qquad 
					X_0=(1,0.5)^{\mathrm T},
				\end{equation}
				where $X_t=(x_t,y_t)^{\mathrm T}$, $B_t=(B_t^1,B_t^2)^{\mathrm T}$ is a standard two-dimensional Brownian motion, and
				\[b(z)=-z|z|+6z-2z\cos(|z|),\qquad z=(x,y)^{\mathrm T}\in\R^2.\]
				In the numerical experiment, we choose a correlated non-degenerate diffusion matrix
				\[
				\Sigma=
				\begin{pmatrix}
					1&0\\
					0.6&0.8
				\end{pmatrix}.
				\]
				The eigenvalues of \(\Sigma\Sigma^{\mathrm T}\) are $1.6$ and $0.4$. Hence $\Sigma\Sigma^{\mathrm T}>0$.
			\end{example}
			We first verify that the drift coefficient in Example \ref{EX} satisfies Assumption \ref{A1}.	
			A direct computation shows that
			\[
			\left\langle z_1-z_2,b(z_1)-b(z_2)\right\rangle
			\leqslant 8|z_1-z_2|^2,\qquad z_1,z_2\in\R^2.
			\]
			Moreover, the dominant term \(-z|z|\) yields dissipation at large distances. More precisely, for all \(z_1,z_2\in\R^2\) satisfying \(|z_1-z_2|>20\), one has
			\[
			\left\langle z_1-z_2,b(z_1)-b(z_2)\right\rangle
			\leqslant -|z_1-z_2|^2 .
			\]
			Consequently, condition \eqref{CaI} holds with \(L=8\), \(K=1\) and \(R=20\). Thus Assumption \ref{A1} is satisfied. By Lemmas \ref{Ju} and \ref{LeX}, system \eqref{Ex} admits a unique global solution and a unique invariant probability measure. The drift in Example \ref{EX} has a natural interpretation as an overdamped diffusion in a confining rough radial potential. The term \(-z|z|\) gives the dominant long-range restoring force, while \(-2z\cos(|z|)\) produces local radial corrugations. This is closely related to Zwanzig's classical model of diffusion in a rough potential \cite{Zwanzig}, where microscopic roughness of the energy landscape affects Brownian diffusion. Thus, Example \ref{EX} may be regarded as a two-dimensional radial analogue of rough-potential diffusion, combining long-range dissipativity with local oscillatory effects.

			It is important to note that the drift coefficient in \eqref{Ex} is not only dissipative at infinity, but also genuinely nonlinear, oscillatory, and of superlinear growth. Therefore, the numerical approximation theories developed in the existing literature, for instance \cite{Bao, Bao1, Bao2, Majka, A.D. N}, do not directly apply to the present example. This motivates the use of the TEM scheme proposed in this paper, which is designed to handle such nonlinear systems with long-range dissipativity while preserving the relevant long-time stability structure.
					
Let \(X_k^h=(x_k^h,y_k^h)^{\mathrm T}\) denote the numerical solution at time $t_k=kh$, and set
\[
r_k^h=|X_k^h|
=
\sqrt{(x_k^h)^2+(y_k^h)^2}.
\]
For the truncation function, we choose \(\phi(u)=4u+8, u\geqslant 0\),
so that \(\phi^{-1}(v)=(v-8)/4\). Taking \(M=16\) and \(\theta=1/4\), the truncation radius becomes \(4h^{-1/4}-2\).
The TEM scheme for \eqref{Ex} is then written explicitly as
\[\left\{\begin{aligned}
	\widehat x_{k+1}^h
	&=
	x_k^h+h(-r_k^h+6-2\cos r_k^h)x_k^h
	+\Delta B_{k+1}^{1},\\
	\widehat y_{k+1}^h
	&=
	y_k^h+h(-r_k^h+6-2\cos r_k^h)y_k^h
	+0.6\Delta B_{k+1}^{1}
	+0.8\Delta B_{k+1}^{2},\\
	x_{k+1}^h
	&=
	\left(1\wedge \frac{4h^{-1/4}-2}{\sqrt{(\widehat x_{k+1}^h)^2+(\widehat y_{k+1}^h)^2}}\right)
	\widehat x_{k+1}^h,\\
	y_{k+1}^h
	&=
	\left(1\wedge \frac{4h^{-1/4}-2}{\sqrt{(\widehat x_{k+1}^h)^2+(\widehat y_{k+1}^h)^2}}\right)
	\widehat y_{k+1}^h.
\end{aligned}\right.\]
Here \(\Delta B_{k+1}^{i} =
B_{t_{k+1}}^i-B_{t_k}^i\), \(i=1,2\), and $\Delta B_{k+1}^{1}$, $\Delta B_{k+1}^{2}$ are independent normal random variables distributed as $N(0,h)$. The multiplicative truncation factor is understood to be equal to one when \((\widehat x_{k+1}^h)^2+(\widehat y_{k+1}^h)^2=0\).

\begin{figure}[htb]
				\centering
				\includegraphics[width=15cm]{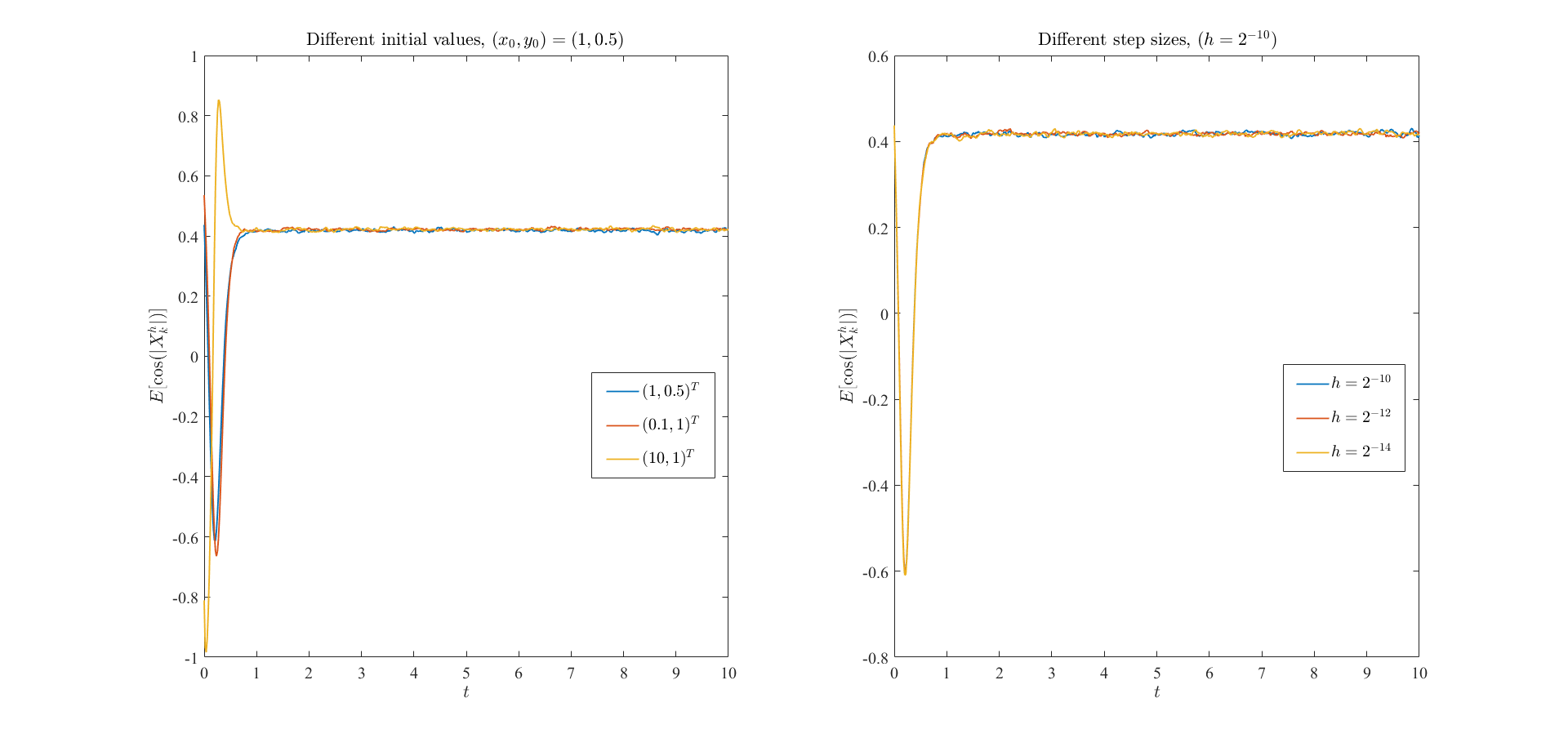}
				\caption{Trajectories of $(\mathbb{E}\cos(|X_k^h|))_{k\geqslant 0}$: Left---fixed step size with different initial conditions; Right---fixed initial condition with different step sizes.}
				\label{fig5}
\end{figure}

The left panel of Figure \ref{fig5} displays Monte Carlo approximations of \(t_k\mapsto \mathbb{E}[\cos(|X_k^h|)]\)
with step size \(h=2^{-10}\), \(2000\) sample paths, and initial values \((x_0,y_0)^{\mathrm T}\in\{(1,0.5),(0.1,1),(10,1)\}\).
The curves approach a common limiting level as time increases, independently of the initial condition. This behaviour is consistent with the convergence-to-equilibrium statement in Theorem \ref{Pro1}. The right panel of Figure \ref{fig5} shows the same quantity for the fixed initial value \((x_0,y_0)^{\mathrm T}=(1,0.5)\) and for step sizes \(h\in\{2^{-10},2^{-12},2^{-14}\}\).
For all the step sizes considered, the numerical trajectories approach an equilibrium level. This provides numerical evidence for the existence of numerical invariant measures and for their convergence as the step size decreases.
			
\begin{figure}[htb]
	\centering
	\includegraphics[width=12cm]{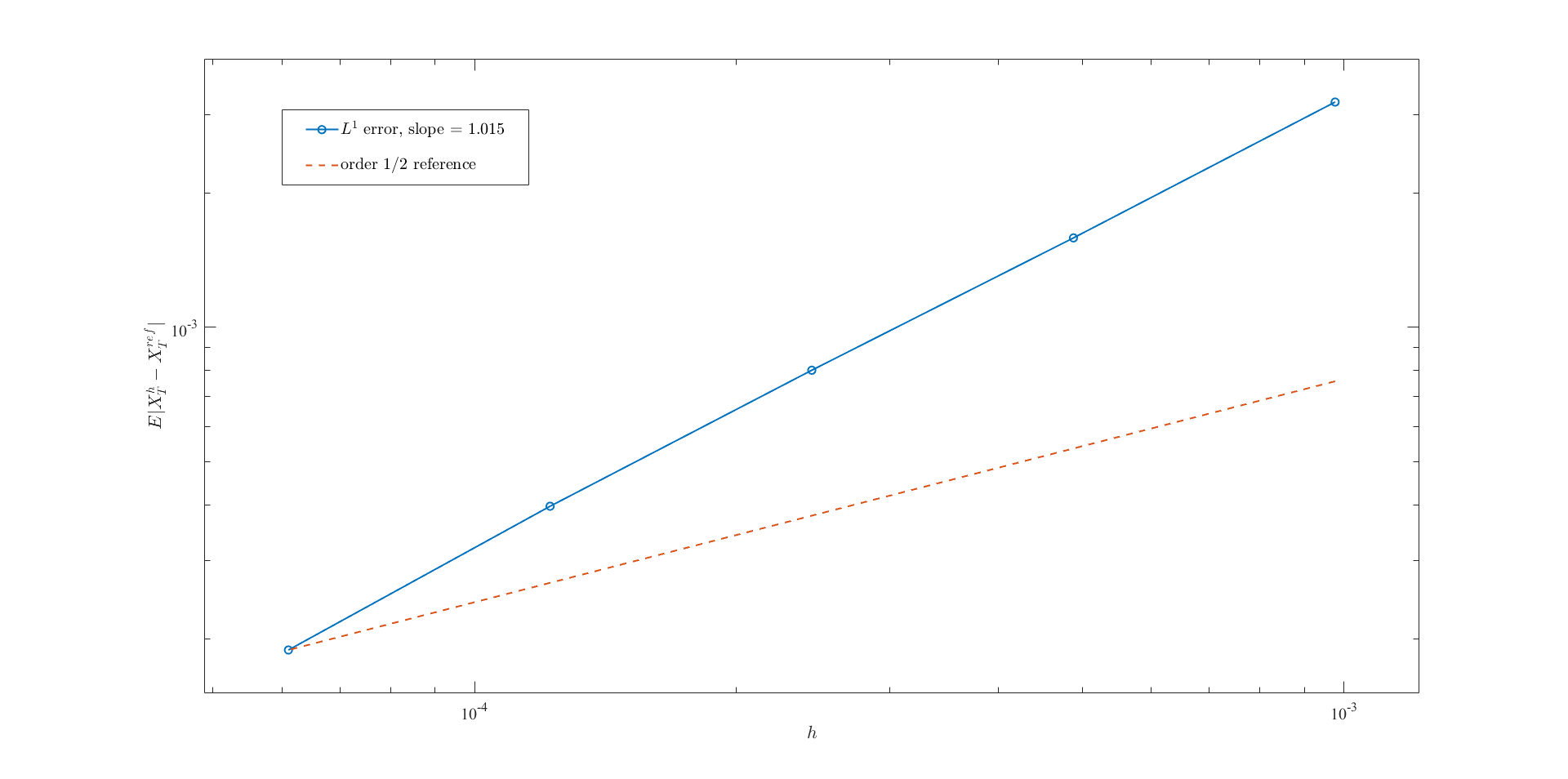}
	\caption{Terminal-time strong error between the TEM approximation and the reference approximation as a function of the step size \(h\).}
	\label{fig2}
\end{figure}

We next illustrate the strong convergence behaviour predicted by Theorem \ref{Th2}. The theorem gives a uniform-in-time strong convergence rate of at least order \(1/2\). Since no closed-form solution of \eqref{Ex} is available, we use the TEM approximation with the smaller step size \(h_{\mathrm{ref}}=2^{-17}\) as the reference solution. We fix \(T=16\), \((x_0,y_0)^{\mathrm T}=(1,0.5)\), and use \(2000\) sample paths. Figure \ref{fig2} plots the terminal-time strong error \(\mathbb{E}|X_T^h-X_T^{\mathrm{ref}}|\) for \(h\in\{2^{-10},2^{-11},2^{-12},2^{-13},2^{-14}\}\).
Here \(X_T^h\) denotes the TEM approximation with step size \(h\), while \(X_T^{\mathrm{ref}}\) denotes the reference approximation with step size \(h_{\mathrm{ref}}\). The error decreases as the step size becomes smaller. The observed slope is close to the first-order reference line; in particular, it is consistent with, and in this example appears sharper than, the rate guaranteed by Theorem \ref{Th2}.

\begin{figure}[htb]
	\centering
	\includegraphics[width=15cm]{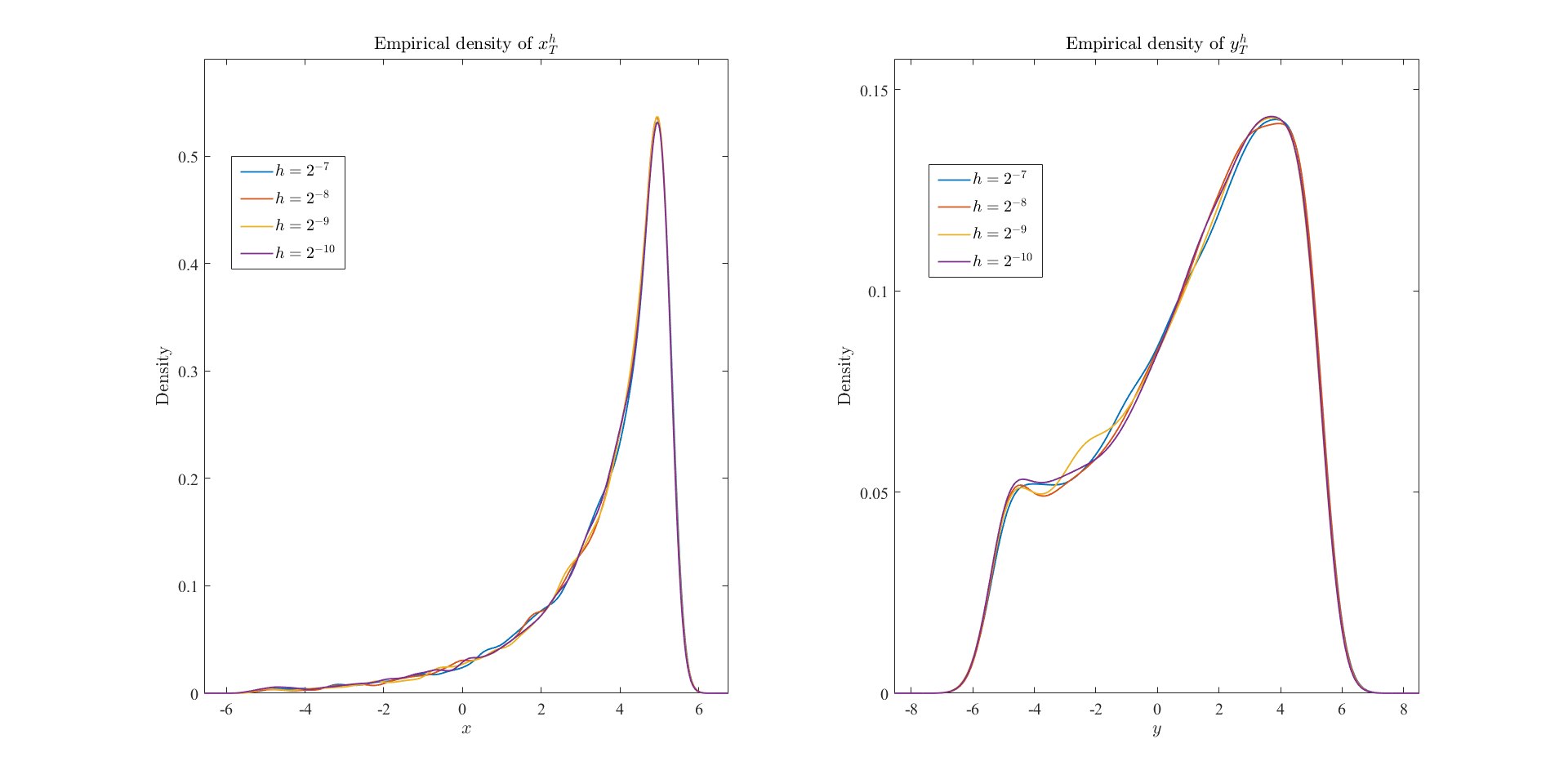}
	\caption{Empirical marginal densities of the TEM approximation at terminal time \(T=16\). Left: empirical density of \(x_T^h\). Right: empirical density of \(y_T^h\).}
	\label{fig3}
\end{figure}

\begin{figure}[htb]
	\centering
	\includegraphics[width=15cm]{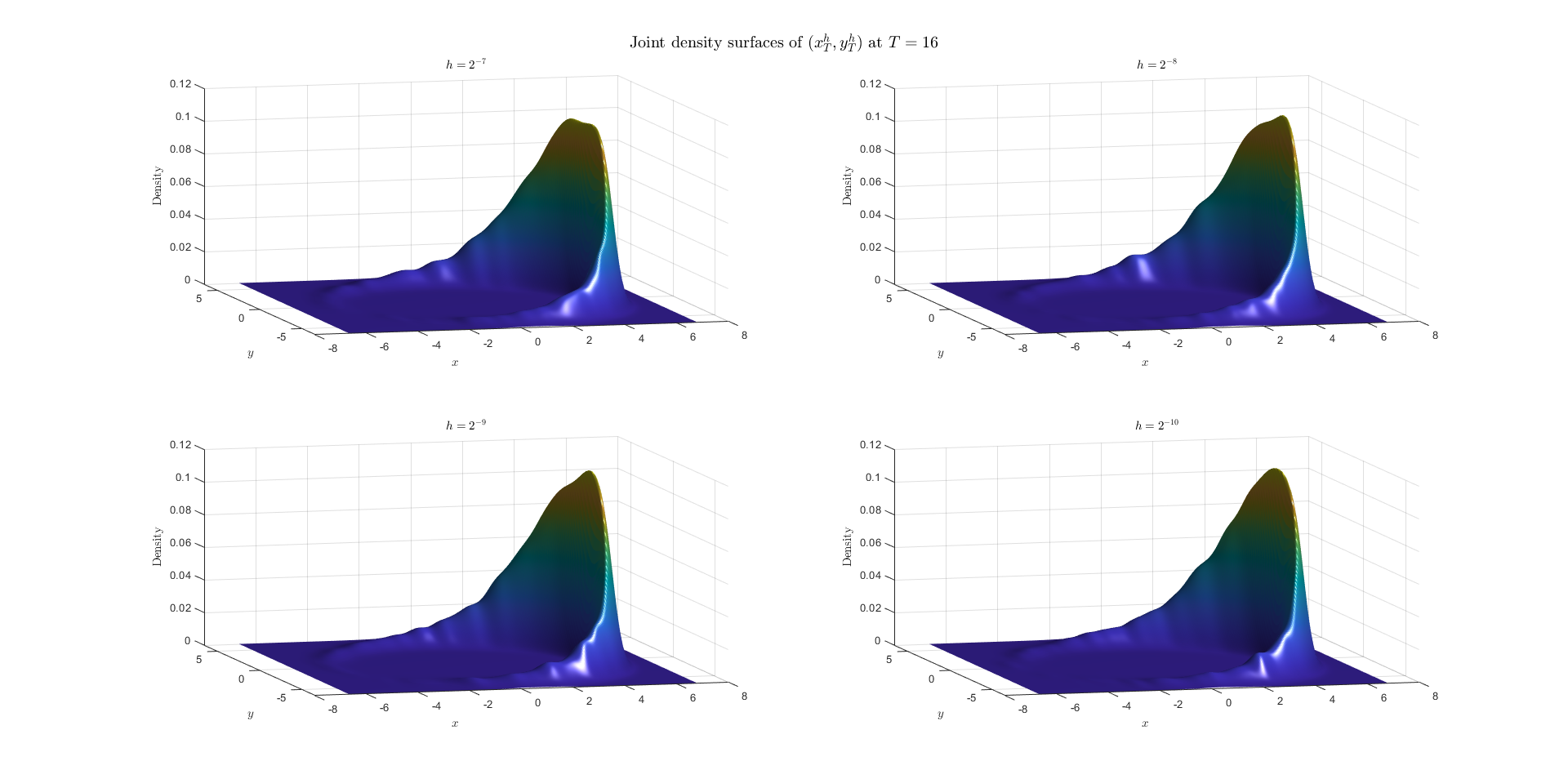}
	\caption{Empirical joint density surfaces of \((x_T^h,y_T^h)\) at terminal time \(T=16\), estimated from \(8000\) sample paths.}
	\label{fig4}
\end{figure}

Finally, we examine the convergence of numerical invariant measures. By Theorem \ref{Th3}, the invariant measure of the TEM scheme converges to the
invariant measure of the exact solution as \(h\to0\). Since an explicit expression for the invariant density of \eqref{Ex} is not available, we use the
empirical law of \(X_T^h\) at a sufficiently large terminal time \(T=16\) as a numerical proxy for the invariant law, and compare the resulting empirical
densities for different step sizes. With initial value
\((x_0,y_0)^{\mathrm T}=(1,0.5)\) and \(8000\) sample paths, Figures \ref{fig3} and \ref{fig4} display, respectively, the empirical marginal densities of
\(x_T^h\) and \(y_T^h\), and the three-dimensional kernel density estimates of the joint law of \((x_T^h,y_T^h)\), for \(h\in\{2^{-7},2^{-8},2^{-9},2^{-10}\}\). The marginal densities corresponding to different step sizes almost overlap, while the joint density surfaces are also visually close to each other and become increasingly stable as \(h\) decreases. These observations indicate that the long-time distribution of the numerical solution is stable with respect to time discretisation, and provide numerical support for the convergence of numerical invariant measures asserted in Theorem \ref{Th3}.
		
\section{Concluding remarks}
We have proposed an explicit truncated Euler--Maruyama scheme for additive-noise SDEs with superlinearly growing drifts that are dissipative only at infinity. The scheme preserves the long-time stability of the underlying diffusion: it admits a unique numerical invariant measure and converges exponentially to it in Wasserstein distance. By means of a fully discrete coupling argument adapted to the TEM transition kernel and to the \(\Gamma^{-1}\)-geometry of the noise, we obtained a uniform-in-time strong error estimate of order \(1/2\). This estimate further yields an order \(1/2\) approximation of the invariant measure in the \(L^1\)-Wasserstein distance. The numerical experiments are consistent with these theoretical results. Extensions to multiplicative or degenerate noise, and to higher-order explicit schemes, remain natural directions for future work.

			%%%%%%%%%%%%%%%%%%%%%%%%%%%%%%%%%%%%%%%%%%%%%%
			%% Support information, if any,             %%
			%% should be provided in the                %%
			%% Acknowledgements section.                %%
			%%%%%%%%%%%%%%%%%%%%%%%%%%%%%%%%%%%%%%%%%%%%%%
			%\begin{acks}[Acknowledgments]
			%	The authors would like to thank the anonymous referees, an Associate Editor and the Editor for their constructive comments that improved the quality of this paper.
			%\end{acks}
			
			%%%%%%%%%%%%%%%%%%%%%%%%%%%%%%%%%%%%%%%%%%%%%%
			%% Funding information, if any,             %%
			%% should be provided in the                %%
			%% funding section.                         %%
			%%%%%%%%%%%%%%%%%%%%%%%%%%%%%%%%%%%%%%%%%%%%%%
			\begin{funding}
				The second author was supported by the  National Natural Science Foundation of China (No. 12371402) and the Tianjin Natural Science Foundation (24JCZDJC00830).
			\end{funding}

\end{document}